\documentclass{tac}

\usepackage{graphicx}%
\usepackage{multirow}%
\usepackage{mathrsfs}%
\usepackage[title]{appendix}%
\usepackage{xcolor}%
\usepackage{textcomp}%
\usepackage{manyfoot}%
\usepackage{booktabs}%
\usepackage{algorithm}%
\usepackage{algorithmicx}%
\usepackage{algpseudocode}%
\usepackage{listings}%
\usepackage{geometry}%
\usepackage{bm}
\usepackage[all,2cell]{xy} \UseAllTwocells \SilentMatrices
\usepackage{xy}
\input diagxy

\usepackage[colorlinks=true]{hyperref}
\hypersetup{allcolors=[rgb]{0.1,0.1,0.4}}

\newcommand{\maps}{\colon}

\newcommand{\Cat}{{\mathscr Cat}}
\newcommand{\Cats}{{\mathscr CAT}}

\newcommand{\Set}{{\mathscr Set}}

\newcommand{\Y}{{\mathscr Y}}
\newcommand{\X}{{\mathscr X}}

\newcommand{\A}{{\mathscr A}}
\newcommand{\B}{{\mathscr B}}
\newcommand{\M}{{\mathscr M}}

\newcommand{\Fib}{{\mathscr Fib}}

\newcommand{\N}{{\mathscr N}}

\newcommand{\Pe}{{\mathscr P}}

\newcommand{\I}{{\mathscr I}}

\newcommand{\E}{{\mathscr E}}

\newcommand{\C}{{\mathscr C}}

\newcommand{\G}{{\mathscr G}}

\newcommand{\F}{{\mathscr F}}

\newcommand{\T}{{\mathscr T}}

\newcommand{\W}{{\mathscr W}}

\newcommand{\D}{{\mathscr D}}

\newcommand{\ce}{{\bm \xi}}
\newcommand{\de}{{\bm d}}
\newcommand{\del}{\bm{\delta}}

\newcommand{\n}{{\bm \nu}}

\author{Igor Bakovi\' c}

\thanks{The author would like to thank Nathanael Arkor, Ivan di Liberti,  Fosco Loregian, Paul-Andr\' e Melli\` es, Axel Osmond and Paul Sejourne for numerous inspiring discussions, their friendship and support over the last few years.  The author would also like to thank to Olivia Caramello, Alain Connes and Laurent Lafforgue as the members of the scientific committee of "Toposes in Mondov\` i" for excepting my talk for a presentation of a part of this program at the conference. The last but not the least, the author would like to express a genuine gratitude to a dear friend and colleague Zoran \v Skoda for his continuous encouragement, guidance and support over the decades without whom this work would have probably never seen the light of the day.}

\address{}

\title{Comma 2-comonad I: Eilenberg-Moore 2-category of colax coalgebras}

\copyrightyear{2025}


\eaddress{ibakovic@gmail.com}

\newtheorem{theorem}{Theorem}

\newtheoremrm{rem}{Remark}

\mathrmdef{Hom}
\mathbfdef{Set}

\begin{document}

\maketitle
\begin{abstract}
In this paper we describe a comma 2-comonad on the 2-category whose objects are functors, 1-cell are colax squares and 2-cells are their transformations. We give a complete description of the Eilenberg-Moore 2-category of colax coalgebras, colax morphisms between them and their transformations and we show how many fundamental constructions in formal category theory like adjoint triples, distributive laws, comprehension structures, Frobenius functors etc. naturally fit in this context.
\end{abstract}

% NOTE: it is good practice to \label all headings (and proclamations) immediately

\section{Introduction}\label{sec-Introduction}
The operation which associates to any pair of categories $\A$ and $\B$ a category $(\A,\B)$ which came to be known by the name "comma category" was introduced by Lawvere in his thesis \cite{Law0} for the purpose of a foundational clarification, in particular of the notion of {\it adjointness}. Since then the comma construction turned out to be fundamental in many areas of category theory from computing Kan extensions to the general calculus of adjoints and limits.  Based on his observation from 1963 that cartesian closed categories serve as a common abstraction of type theory and propositional logic,  Lawvere used the comma construction in his paper \cite{Law3} in which he summed up a stage of the development of the relationship between category theory and proof theory based on the strategy to interpret proofs themselves as structures, following idea which goes back to Kreisel and Kleene.  In that paper he introduced a notion of a \emph{hyperdoctrine} consisting of a category $\T$ of \emph{types} and a functor
\[\Pe \maps \T^{op} \to \Cats \]
which served as a main device in his subsequent influential paper \cite{Law4} to to find an analogue in category theory to the comprehension scheme of set theory which says, essentially, that given a property, there is a set consisting exactly of the elements having that property.  Lawvere translated this statement into an adjunction
\[\xymatrix{(\Set,X) \ar@<-1.2ex>[r]_{}^{\perp} & \Pe(X) \ar@<-1.2ex>[l]_{}} \]
where the comma category on the left consists of all functions into the fixed set $X$. Then he generalized the adjunction to the case of the presheaf hyperdoctrine
\[\xymatrix{(\Cat,\X) \ar@<-1.2ex>[r]_{}^{\perp} & \Pe(\X) \ar@<-1.2ex>[l]_{}} \]
where $\X$ is a category and $\Pe(\X)$ is a category of presheaves on $\X$ and a functor from right to left is given by a Grothendieck construction assigning a (discrete) fibration to any presheaf.  A further step was done by Gray in \cite{Gray2} based on his reformulation \cite{Gray1} of Grothendieck's theory of (co)fibrations \cite{Gr2} in order to discuss what it would mean for an adjunction
\[\xymatrix{(\Cat,\X) \ar@<-1.2ex>[r]_{}^-{\perp} & \Cats[\X^{op},\Cat] \ar@<-1.2ex>[l]_{}} \]
to be an instance of the comprehension scheme. Gray's answer is that it is not, and the reason is that in this context the comprehension scheme is equivalent to asserting
\[{\displaystyle  \lim_{\to }}(\tau_{\B}) = \B \]
where $\tau_{\B} \maps \B \to \bullet$ is a constant functor into the terminal category $\bullet$ in $\Cat$, and this assertion is false unless $\B$ is discrete. Gray realized that the difficulty lied in applying the above notion of a colimit in a traditional 1-dimensional sense which is not suitable for $\Cat$ being intrinsically a 2-category. Therefore he started developing the theory of functors, comma categories, adjointness and limits in a 2-dimensional context culminating in a milestone monograph \cite{Gray} which meant to be the first part of a projected four volumes, but its subsequent parts never appeared (as well as "The calculus of comma categories" which he promised in the bibliography of \cite{Gray2}). In the first paragraph of \cite{Gray} he explains the subject of the formal category theory:
\begin{quote}
The purpose of category theory is to try to describe certain general aspects of the structure of mathematics. Since category theory is also part of mathematics, this categorical type of description should apply to it as well as to other parts of mathematics.
\end{quote}

Already until the publication of \cite{Gray} there has been a considerable amount of work on various aspects of this study by various authors, most notably by B\' enabou \cite{Be1}, Eilenberg and Kelly \cite{EK} and \cite{Kelly4}, Linton \cite{Lin} and Street \cite{Street1} and \cite{Street2} which was directly concerned with, or at least, relevant to it. The theory of adjoint squares which is the subject of the last chapter of \cite{Gray} was previously studied by Palmquist \cite{Palm} and Maranda \cite{Mara} in the context of categories enriched over the closed monoidal category in the sense of \cite{EK}.

One of the purposes of this paper is to give some conceptual clarifications of the theory of adjoint squares using a theory of 2-adjunctions and 2-comonads.  Particularly, there exists a strict 2-adjunction
\[\xymatrix{\Cat^{2}_{c}  \ar@<-1.2ex>[r]_{D}^-{\perp} & \Cat \ar@<-1.2ex>[l]_{I}} \]
where $\Cat^{2}_{c}$ is the 2-category whose objects are functors and 1-cells are colax squares filled by upwards pointing natural transformations with an obvious notion of 2-cells between the latter.  The 2-functor $D \maps \Cat^{2}_{c} \to \Cat$ sends any functor $U \maps \A \to \X$ to its comma category $(\X,U)$ and $I \maps \Cat \to \Cat^{2}_{c}$ is an obvious 2-functor which sends any category $\X$ to identity functor $I_{\X} \maps \X \to \X$. From the perspective of formal category theory it is very important that this 2-adjunction does not require any further restrictions on categories or functors whatsoever. This point of view is very much in the spirit of B\' enabou's philosophy of regarding functors as \emph{generalized fibrations}. As any other, this 2-adjunction generates a 2-comonad $\D$ on $\Cat^{2}_{c}$ and the main contribution of this paper is a description of colax $\D$-coalgebras and investigation of its properties. In the light of the aforementioned line of thought it should not come as a surprise that these complex mathematical structures are lurking behind many important constructions in formal category theory, most notably in the \emph{formal theory of monads} which Street developed in \cite{Street1} followed by its sequel with Lack in \cite{LSt} and based on B\' enabou's realization that the definition of a monad can be made in an arbitrary bicategory (which was actually the main reason why B\' enabou insisted on incorporating \emph{morphisms of bicategories} or \emph{lax functors} as the central part of the theory of bicategories). Although Lack and Street claim in Remark 1.1. 
of \cite{LSt} that it would be possible to develop the whole formal theory of monads in the context of general bicategories,  the extension of their theory from 2-categories is not starightforward as it seems. The reason lies in the concluding remark behind Proposition I.6.4 on page 142 of \cite{Gray} that Yoneda embedding for bicategories is not locally full and Lack and Street's argument using 2-categorical Yoneda lemma cannot be carried out to the case of bicategories.

The important attempt to develop such theory based on generalized operads was done by Chikhladze in \cite{Chik} by means of \emph{generalized multicategories} or \emph{$T$-monoids} who defined as monads within a Kleisli bicategory. Motivated by investigation on functoriality of the Kleisli construction Chikhladze developed a lax version of the formal theory of monads, and studied its connection to bicategories. 

The theory of generalized multicategories itself has its early origins in the work of Burroni \cite{Burr} but the more contemporary approaches mostly use the theory developed by Hermida in \cite{Her1} (for the excellent overview of the state of affairs in the filed the diligent reader may consult Leinster's book \cite{Le}).

The point of departure from bicategories in the development of the theory was Cruttwell and Shulman's paper \cite{CS} who recognized that generalized multicategories have been defined in numerous contexts throughout the literature as the \emph{lax algebras} or \emph{Kleisli monoids} relative to a pseudomonad on a bicategory but they questioned the meanings of these constructions from author to author. They proposed a unified framework for generalized multicategories by means of \emph{virtual double categories} contenting themselves by working with monads on double categories and related structures rather than bicategories. An important role in their theory is played by normal oplax $T$-algebras
where $T$ is a normal monad on a virtual equipment. In the statement of Theorem 9.13. they construct a strict 2-functor whose underlying 1-functor is fully faithful, and which becomes 2-fully-faithful (that is, an isomorphism on hom-categories) when restricted to normal oplax T-algebras, with a remark immediately after the theorem that the restriction to normal oplax algebras in its final statement cannot be dispensed with either.  Based on this result Shulman introduced in \cite{Shu} so called \emph{twisted $G$-actions} in order to define 2-categories with contravariance based on his axiomatic treatment of  duality involutions of categories.  The purpose was to treat both covariant and contravariant morphisms inside 2-categories and it is the author's opinion that this is one of the most important approaches to deal with such questions in categorical structures.

However it is also a conviction of the author that the main reason why there is a proliferation of such unifying theories is caused by insufficient understanding of bicategories, in particular of adjunctions in bicategories and consequently higher adjunctions in Gray-categories and other higher categorical structures. This paper tends to bridge this gap by starting a program which might be called \emph{formal theory of adjunctions} by paraphrasing Street's formal approach to monads. One of the most important discoveries of the author during the last few years is that the associated split fibration 2-monad $\F$ whose pseudo-algebras are Grothendieck fibrations extends to the 2-category $\Cat^{2}_{c}$ from the beginning of this paper and moreover it turns out that this 2-monad is admissible in the sense of Bunge and Funk \cite{BF}.

The comma constructions from the beginning plays a crucial role in several 2-(co)monads on the 2-category $\Cat^{2}_{c}$ and its cousins.  The dual version of the comma 2-comonad which is the main object of the study of this paper was investigated by Pavlovi\'  c in \cite{Pav}. He studied three comonads derived from the comma construction and its induced coalgebras which correspond to the three main concepts of his paper: cofree equivalences, dualities and $\ast$-autonomous categories.He showed that the comonad which yields the $\ast$-autonomous categories is the Chu construction and then he proceeded to show that it is induced by the right adjoint to the inclusion of $\ast$-autonomous categories among autonomous categories, with lax structure-preserving morphisms. Moreover,  Pavlovi\'  c showed that this inclusion turns out to be comonadic: autonomous categories are exactly the Chu-coalgebras. This paper was one of the main motivations for the author to start an investigation of the comma 2-comonad in foundations. 

The other important source of inspiration for the author was Maltsiniotis paper \cite{Maltsiniotis} whose impetus for his work was the observation that the class of fibrations shares many formal properties with the class of $\W$-smooth functors where $\W$ is a basic localizer in the sense of Grothendieck \cite{Gr}. The crucial observation which we will address in the sequels ton this paper is his observation that the class of functors which have a right adjoint shares many properties with the class of $W$-aspheric functors.

However this paper is not so ambitious as its introduction seem to suggest. Its basic purpose is to set up the scenery which will be used in its sequels for a sudden change in order to reflect on the unifying aspects which it seems to allow for the pillars of the formal category theory:
 \begin{itemize}
\item [1)] Representability
\item  [2)] Coherence 
\item [3)] Duality
\end{itemize}
Its structure and content is inevitably very simple, due to the complexity of diagrams and structures which it describes. Like a famous Chinese proverb says: "A journey of a thousand miles begins with a single step."

\section{A comma 2-comonad}
Let $\Cat^{2}_{c}$ denote a 2-category whose objects are functors and 1-cells are colax squares
\[\xymatrix@!=3pc{\A \ar[r]^-{F} \ar[d]_-{S}  & \E \ar[d]^{G} \\
\B \ar[r]_-{V} & \X \ultwocell<\omit>{\phi}}\]
with a natural transformation $\phi$ pointing upwards. For any two 1-cells, a 2-cell
\[\xymatrix@!=0.5pc{\E \ar[dd]_{U} \ar@/^1pc/[rr]^{F} \ar@/_1pc/[rr]_{F'} \rrtwocell<\omit>{\tau} && \A \ar[dd]^{G} \\
&&  \\
\B \ar@/^1pc/@{-->}[rr]^{V} \ar@/_1pc/[rr]_{V'} \rrtwocell<\omit>{\sigma} && \X \uulltwocell<\omit>{\phi'} }\]
is given by the above diagram with an obvious commutativity condition. There is a canonical 2-functor
\begin{equation}
\begin{array}{c}\label{Inc}
I \maps \Cat \to \Cat^{2}_{c}
\end{array}
\end{equation}
which sends a category $\B$ to the identity functor $I_{\B} \maps \B \to \B$.  We have the following:

\begin{theorem}\label{Rightcomma}
There exists a strict 2-adjunction
\begin{equation}\label{ID}
\begin{array}{c}
I \dashv D
\end{array}
\end{equation} 
where $D \maps \Cat^{2}_{c} \to \Cat$ sends any functor $U \maps \A \to \X$ to its comma category $(\X,U)$.
\begin{proof}
The action of $D$ on a 1-cell $(U,\beta,F) \maps P \to G$ in $\Cat^{2}_{c}$ is induced by the universal property of a comma category $(\X,G)$
\[\xymatrix@!=0.2pc{(\B,P) \ar[rrr]^{\de_{0}} \ar[ddd]_{\de_{1}} \ar[dr]^(0.6){D(U,\beta,F)} &&& \E \ar@{-->}[ddd]^{P} \ar[dr]^{F} &\\
& (\X,G) \ar[rrr]^{\de_{0}} \ar[ddd]_{\de_{1}} &&& \A \ar[ddd]^{G} \\
&& \ultwocell<\omit>{\delta_{P}\,\,\,\,\,} &&  \\
\B \ar[dr]_{U} \ar@{==}[rrr]^{} &&& \B \ar@{-->}[dr]_(0.4){U} \ultwocell<\omit>{\delta_{G}\,\,\,\,\,} & \uultwocell<\omit>{\beta}   \\
& \X \ar@{=}[rrr]_{} &&& \X }\]
For any 2-cell $(\theta,\tau) \maps (U,\beta,F) \Rightarrow (V,\gamma,K)$ as in the following diagram 
\[\hspace{-1cm}\xymatrix@!=0.75pc{(\B,P) \ar[rrr]^-{\de_{0}} \ar[ddd]_{\de_{1}} \ar@/_1.25pc/[dr]_(0.3){D(U,\beta,F)} \ar@/^1.25pc/[dr]^(0.7){D(V,\gamma,G)} &&&
\E \ar@{-->}[ddd]^{P} \ar@/_1.25pc/[dr]_(0.3){F} \ar@/^1.25pc/[dr]^{K} &\\
& (\X,G) \ar[ddd]^{\de_{1}} \ar[rrr]^-{\de_{0}} \ultwocell<\omit>{(\theta,\tau),\,\,\,\,\,\,} &&& \A \ar[ddd]^{G} \ultwocell<\omit>{\tau} \\
&&&& \\
\B \ar@{-->}[rrr]_-{I_{\B}} \ar@/_1.25pc/[dr]_{U} \ar@{-->}@/^1.25pc/[dr]^{V} &&& \B \ar@{-->}@/_1.25pc/[dr]_(0.4){U} \ar@{-->}@/^1.25pc/[dr]^{V} & \uultwocell<\omit>{<-2>\beta}  \uultwocell<\omit>{<2>\gamma}\\
& \X \ar[rrr]_-{I_{\X}} \ultwocell<\omit>{\theta} &&& \X \ultwocell<\omit>{\theta}}\]
we have functors $D(U,\beta,F)$ and  $D(V,\gamma,K)$ taking any object $(B,p,E)$ of $(\B,P)$ to 
\[\begin{array}{c}
D(U,\beta,F)(B,p,E):=(U(B),\beta_{E}U(p),F(E))\\
D(V,\gamma,G)(B,p,Q):=(V(B),\gamma_{E}V(p),K(E))
\end{array}\]
respectively. The two functors take any morphism $(u,e) \maps (B,p,E) \to (B',p',E')$ to 
\[\begin{array}{c}
D(U,\beta,F)(u,e):=(U(u),F(e)) \\
D(V,\gamma,K)(u,e):=(V(u),K(e))
\end{array}\]
respectively, from which it follows that the action of $D$ on a 2-cell $(\theta,\tau)$ is a 2-cell $(\theta,\tau) \maps D(F,\beta,U) \Rightarrow D(G,\gamma,V)$ (denoted by the same pair of natural transformations) whose component $(\theta,\tau)_{(B,p,Q)}$ indexed by an object $(B,p,P)$ is a morphism $(\theta_{B},\tau_{E}) \maps (U(B),\beta_{E}U(p),F(E)) \to (V(B),\gamma_{E}V(p),G(E))$ represented by the back side of a diagram
\[\xymatrix@!=0.75pc{U(B) \ar[dd]_{U(p)} \ar[dr]^{U(u)} \ar[rrr]^-{\theta_{B}} &&&
V(B) \ar@{-->}[dd]_(0.6){V(p)} \ar[dr]^{V(u)} &\\
& U(B') \ar[rrr]^{\theta_{B'}} \ar[dd]_(0.3){U(p')} &&& \F \ar[dd]^{V(p')} \\
UP(E) \ar@{-->}[rrr]^{\theta_{P(E)}} \ar[dr]^(0.4){UP(e)}  \ar[dd]_-{\beta_{E}} &&& 
VP(E) \ar@{-->}[dr]_(0.3){VP(e)} \ar@{-->}[dd]^(0.6){\gamma_{E}} &  \\
& UP(E') \ar[rrr]^{\theta_{P(E')}}  \ar[dd]_-{\beta_{E'}}  &&& VP(E') \ar[dd]^{\gamma_{E'}}  \\
GF(E) \ar[dr]_{GF(e)} \ar@{-->}[rrr]^{G(\tau_{E})} &&& GK(E) \ar@{-->}[dr]_{GK(e)} & \\
& GF(E') \ar[rrr]_{G(\tau_{E'})} &&& GK(E')  }\]
The above construction is clearly functorial with respect to vertical composition of (pairs of) natural transformations forming 2-cells in hom-categories of $\Cat^{2}_{c}$ and it is easily established that it is also functorial with respect to a horizontal composition. The counit $\del_{G} \maps ID(G) \to G$ of the 2-adjunction (\ref{ID}) indexed by an object $G$ in $\Cat^{2}_{c}$
\begin{equation}
\begin{array}{c}\label{counit}
\xymatrix@!=2pc{(\X,G) \ar[r]^-{\de_{0}} \ar@{=}[d]_{} & \A \ar[d]^{G} \\
(\X,G) \ar[r]_-{\de_{1}} & \X  \ultwocell<\omit>{\delta_{G}\,\,\,\,\,} }
\end{array}
\end{equation} 
is a 1-cell $(\de_{1},\delta_{G},\de_{0}) \maps I_{(\X,G)} \to G$ given by a diagram (\ref{counit}) and its universal property says that for any 1-cell  $(U,\beta,F) \maps I_{\E} \to G$ in $\Cat^{2}_{c}$ there exists a unique functor $(U,\beta,F) \maps \E \to (\X,G)$ (again denoted by the same pair of natural transformations) such that a diagram 
\[\xymatrix@!=0.25pc{\E \ar[drrrr]^{F} \ar@{=}[ddd]_{} \ar[dr]_{(U,\beta,F)} &&&&\\
& (\X,G) \ar[rrr]^{\de_{0}} \ar@{=}[ddd]_{} &&& \A \ar[ddd]^{G} \\
&&&&  \\
\E \ar[dr]_{(U,\beta,F)} \ar@{-->}[drrrr]^{U} &&& \uulltwocell<\omit>{\beta} \ultwocell<\omit>{\delta_{G}\,\,\,\,\,} &  \\
& (\X,G) \ar[rrr]_{\de_{1}} &&& \X }\]
commutes.  The 1-cell $ID(U,\beta,F):=(D(U,\beta,F),D(U,\beta,F)) \maps ID(P) \to ID(G)$ fits into a commutative diagram
\[\xymatrix@!=0.5pc{(\B,P) \ar[rrr]^{\de_{0}} \ar@{=}[ddd]_{} \ar[dr]^{D(U,\beta,F)} &&& \E \ar@{-->}[ddd]^{P} \ar[dr]^{F} &\\
& (\X,G) \ar[rrr]^{\de_{0}} \ar@{=}[ddd]_{} &&& \A \ar[ddd]^{G} \\
&& \ultwocell<\omit>{\delta_{P}\,\,\,\,\,} &&  \\
(\B,P) \ar[dr]_{D(U,\beta,F)} \ar@{-->}[rrr]^{\de_{1}} &&& \B \ar@{-->}[dr]_(0.4){U} \ultwocell<\omit>{\delta_{G}\,\,\,\,\,} & \uultwocell<\omit>{\beta}   \\
& (\X,G) \ar[rrr]_{\de_{1}} &&& \X }\]
from which we conclude that the counit (\ref{counit}) is strictly natural 2-transformation
\begin{equation}\label{counit2}
\begin{array}{c}
\del \maps ID \to Id_{\Cat^{2}_{c}}.
\end{array}
\end{equation}
The component of the unit $\n \maps Id_{\Cat} \Rightarrow DI$ of the 2-adjunction (\ref{ID}) indexed by a category $\E$ is a functor
\begin{equation}
\begin{array}{c}\label{unit}
\n_{\E} \maps \E \to \E^{2}
\end{array}
\end{equation}
which takes any object $E$ of $\E$ to the identity morphism $1_{E} \maps E \to E$. Its universal property says that for any object $G$ in $\Cat^{2}_{c}$ 
\[\xymatrix@!=0.5pc{&& \E \ar[dll]_-{\n_{\E}} \ar[drr]^-{(U,\beta,F)} && \\
\E^{2} \ar[rrrr]_-{D(U,\beta,F)} &&&& (\X,G) }\]
and any functor $(U,\beta,F) \maps \E \to (\X,G)$ where $\beta \maps U \Rightarrow GF$ is a natural transformation there exists a unique 1-cell $(U,\beta,F) \maps I_{\E} \to G$ in $\Cat^{2}_{c}$ such that the above diagram commutes. The 1-cell $(U,\beta,F)$ sends any object $(E,f,Y)$ in $\E^{2}$ where $f \maps E \to Y$ is a morphism in $\E$ to an object $(U(E),\beta_{Y}f,F(Y))$ in $(\X,G)$.  The two triangular identities
\[\xymatrix@!=0.5pc{&& D \ar@2[dll]_-{\n D} \ar@{=}[drr]^-{} && \\
DID \ar@2[rrrr]_-{D\del} &&&& D}
\hspace{2cm}
\xymatrix@!=0.5pc{&& I \ar@2[dll]_-{I\n} \ar@{=}[drr]^-{} && \\
IDI \ar@2[rrrr]_-{\del I} &&&& I }\]
says that for any two objects $G$ in $\Cat^{2}_{c}$ and $\E$ in $\Cat$ the following diagrams
\[\xymatrix@!=0.5pc{&& (\X,G) \ar@2[dll]_-{\n_{(\X,G)}} \ar@{=}[drr]^-{} && \\
(\X,G)^{2} \ar@2[rrrr]_-{D(\del_{G})} &&&& (\X,G)}
\hspace{2cm}
\xymatrix@!=0.5pc{&& I_{\E} \ar@2[dll]_-{(\n_{\E},\n_{\E})} \ar@{=}[drr]^-{} && \\
I_{\E^{2}} \ar@2[rrrr]_-{\del_{I_{\E}}} &&&& I_{\E} }\]
commute in $\Cat^{2}_{c}$ and $\Cat$ respectively.  From the universal property of $(\X,G)$
\[\xymatrix@!=1.5pc@R1pc@C1pc{(\X,G)^{2} \ar[rrr]^{\de_{0}} \ar[ddd]_{\de_{1}} \ar[dr]^{D(\de_{1},\delta_{G},\de_{0})} &&& (\X,G) \ar@{==}[ddd]^{} \ar[dr]^{\de_{0}} &\\
& (\X,G) \ar[rrr]^{\de_{0}} \ar[ddd]_{\de_{1}} &&& \A \ar[ddd]^{G} \\
&& \ultwocell<\omit>{\delta_{I_{(\X,G)}}\,\,\,\,\,\,\,\,\,\,\,\,\,\,} &&  \\
(\X,G) \ar[dr]_{\de_{1}} \ar@{==}[rrr]^{} &&& (\X,G) \ar@{-->}[dr]_(0.4){\de_{1}} \ultwocell<\omit>{\delta_{G}\,\,\,\,\,} & \uultwocell<\omit>{\delta_{G}}   \\
& \X \ar@{=}[rrr]_{} &&& \X }\]
for any object $[(X,f,A),(x,a),(X',f',A')]$ in $(\X,G)^{2}$ we obtain an object of $(\X,G)$
\begin{equation}\label{Daction}
\begin{array}{c}
D(\del_{G})[(X,f,A),(x,a),(X',f',A')] = \\
(\de_{1}(X,f,A),\delta_{G}(X',f',A')\de_{1}(x,a),\de_{0}(X',f',A'))=(X,f'x,A')
\end{array}
\end{equation} 
together with an action on morphisms in $(\X,G)^{2}$ described by the following diagram
\begin{equation}\label{Dmor}
\begin{array}{c}
\hspace{-1cm}\xymatrix@!=0.1pc{X \ar[rrr]^{h} \ar[ddd]_{f} \ar[dr]^{x} &&& Y \ar@{-->}[ddd]^{g} \ar[dr]^{y} &&& &&&\\
& X' \ar[rrr]^{k} \ar[ddd]_{f'} &&& Y' \ar[ddd]^{g'} &&&& X \ar[rrr]^{h} \ar[ddd]_{f'x=G(a)f} &&& Y \ar[ddd]^{g'y=G(e)g} \\
&&&&&&&&&&&\\
G(A) \ar[dr]_{G(a)} \ar@{-->}[rrr]^{G(m)} &&& G(E) \ar@{-->}[dr]_(0.3){G(e)}  &&  \ar@{|->}[r]^{D(\del_{G})} &&&&&& \\ 
& G(A') \ar[rrr]_{G(n)} &&& G(E') &&&& G(A') \ar[rrr]_{G(n)} &&& G(E') }
\end{array}
\end{equation}
sending $[(h,m),(k,n)] \maps [(X,f,A),(x,a),(X',f',A')] \to [(Y,g,E),(y,e),(Y',g',E')]$ in $(\X,G)^{2}$ to its diagonal $(h,n) \maps (X,f'x,A') \to (Y,g'y,E')$. The above construction defines a functor 
\begin{equation}\label{DD}
\begin{array}{c}
D(\del_{G}) \maps (\X,G)^{2} \to (\X,G)
\end{array}
\end{equation}
and by defining $\D:=ID$ we obtain a natural 2-transformation 
\begin{equation}\label{DE2}
\begin{array}{c}
\D\del \maps \D^{2} \Rightarrow \D
\end{array}
\end{equation} 
whose component indexed by $G$ is a 1-cell $(D(\del_{G}),\iota_{D(\del_{G})},D(\del_{G}))$ in $\Cat^{2}_{c}$. The natural transformation $\delta_{I_{(\X,G)}} \maps \de_{1} \Rightarrow \de_{0}$ is a component of the 1-cell $\del_{\D(G)}$ 
\[\xymatrix@!=0.5pc{(\X,G)^{2} \ar[rrr]^{\de_{0}} \ar@{=}[ddd]_{} \ar[dr]^{D(\de_{1},\delta_{G},\de_{0})} &&& (\X,G) \ar@{==}[ddd]^{} \ar[dr]^{\de_{0}} &\\
& (\X,G) \ar[rrr]^{\de_{0}} \ar@{=}[ddd]_{} &&& \A \ar[ddd]^{G} \\
&& \ultwocell<\omit>{\delta_{I_{(\X,G)}}\,\,\,\,\,\,\,\,\,\,\,\,\,\,} &&  \\
(\X,G)^{2} \ar[dr]_{D(\de_{1},\delta_{G},\de_{0})} \ar@{-->}[rrr]^{\de_{1}} &&& (\X,G) \ar@{-->}[dr]_(0.4){\de_{1}} \ultwocell<\omit>{\delta_{G}\,\,\,\,\,} & \uultwocell<\omit>{\delta_{G}} \\
& (\X,G) \ar[rrr]_{\de_{1}} &&& \X }\]
on the back side of the above diagram which is equivalent to the naturality diagram 
\[\xymatrix@!=3pc{\D^{2}(G) \ar[r]^-{\del_{\D(G)}} \ar[d]_{\D(\del_{G})} & \D(G) \ar[d]^{\del_{G}} \\
\D(G)  \ar[r]_-{\del_{G}} & G  }\]
This clearly establishes the commutativity of one of the triangle identities. The other one follows by precomposing the component of the counit 
\[\xymatrix@!=3pc{\E^{2} \ar[r]^-{\de_{0}} \ar@{=}[d]_{} & \E \ar@{=}[d]^{I_{\E}} \\
\E^{2}  \ar[r]_-{\de_{1}} & \E  \ultwocell<\omit>{\delta_{I_{\E}}\,\,\,\,\,} }\]
indexed by $I_{\E}$ with a 1-cell $(\n_{\E},\n_{\E}) \maps  I_{\E} \to I_{\E^{2}}$.
\end{proof}
\end{theorem}
One of the classical consequences of Theorem \ref{Rightcomma} is the following:

\begin{theorem}\label{comma2comonad} 
There is a strict 2-comonad $(\D,E,C)$ on the 2-category $\Cat^{2}_{c}$ whose underlying 2-functor
\begin{equation}
\begin{array}{c}\label{D}
\D \maps \Cat^{2}_{c} \to \Cat^{2}_{c}
\end{array}
\end{equation}
is a composition $\D:=ID$ of the pair (\ref{Rightcomma}) of adjoint 2-functors. 
\begin{proof}
The comultiplication is a natural 2-transformation 
\begin{equation}\label{C}
\begin{array}{c}
\ce \maps \D \to \D^{2}
\end{array}
\end{equation}
whose component indexed by an object $G$ in $\Cat^{2}_{c}$ is a functor $ \ce_{G} \maps \D(G) \to \D^{2}(G)$ defined by 
\[\ce_{G}:=I(\n_{D(G)}) \maps ID(G) \to IDID(G)\] 
taking an object $(X,f,A)$ in $(\X,G)$ to the identity $(1_{X},1_{A}) \maps (X,f,A) \to (X,f,A)$
\begin{equation}\label{CC}
\begin{array}{c}
\xymatrix@!=3pc{X \ar[r]^-{1_{X}} \ar[d]_-{f} & X \ar[d]^{f} \\
G(A) \ar[r]_-{G(1_{A})} & G(A) }
\end{array}
\end{equation}
in $(\X,G)$ as in the above  diagram. The axioms for a (strict) 2-comonad require that the following diagrams
\[\xymatrix@!=3pc{\D \ar@{=}[r]^-{} \ar@{=}[d]_-{} & \D \ar[d]^{\ce} \ar@{=}[r]^-{} & \D \ar@{=}[d]^-{} \\
\D & \D^{2} \ar[l]^-{\D \del} \ar[r]_-{\del \D} & \D}
\hspace{2cm}
\xymatrix@!=3pc{\D \ar[r]^-{\ce} \ar[d]_-{\ce} & \D^{2} \ar[d]^{\D \ce} \\
\D^{2} \ar[r]_-{\ce\D} & \D^{3} }\]
of natural 2-transformations commute are straightforward consequences of triangular identities. 
\end{proof}
\end{theorem}
We shall start to investigate colax $\D$-coalgebras by recalling their definition:

\begin{definition}\label{colaxDcoalgebra}
A colax $\D$-coalgebra consists of the following data:

\begin{itemize}
\item a 1-cell $\bm{F}_{G}=(\bm{F}_{1},\bm{\varphi}, \bm{F_{0}}) \maps G \to \D(G)$
\begin{equation}\label{varphi}
\begin{array}{c}
\xymatrix@!=3pc{\A \ar[r]^-{\bm{F}_{0}} \ar[d]_-{G} & (\X,G) \ar@{=}[d]^{} \\
\X \ar[r]_-{\bm{F}_{1}} & (\X,G) \ultwocell<\omit>{\bm{\varphi}} }
\end{array}
\end{equation}
in $\Cat^{2}_{c}$
\item a 2-cell $\bm{\zeta} \maps \iota_{G} \Rightarrow \del_{G}\bm{F}_{G}$ in $\Cat^{2}_{c}$ as in a diagram
\begin{equation}\label{zeta}
\begin{array}{c}
\xymatrix@!=3pc{G \ar[d]_{\bm{F}_{G}} \ar@{=}[r]^{} & G \ar@{=}[d]^-{}\\
\D(G) \ar[r]_-{\del_{G}} & G \ultwocell<\omit>{\bm{\zeta}} }
\end{array}
\end{equation}
\item a 2-cell $\bm{\theta} \maps \D(\bm{F}_{G})\bm{F}_{G} \Rightarrow \ce_{G}\bm{F}_{G}$ in $\Cat^{2}_{c}$ as in a diagram
\begin{equation}\label{theta}
\begin{array}{c}
\xymatrix@!=3pc{G \ar[d]_{\bm{F}_{G}} \ar[r]^-{\bm{F}_{G}} & \D(G) \ar[d]^-{\D(\bm{F}_{G})} \\
\D(G) \ar[r]_-{\ce_{G}} & \D^{2}(G) \ultwocell<\omit>{\bm{\theta}} }
\end{array}
\end{equation}
\end{itemize}
such that the following axioms are satisfied:

\begin{itemize}
\item a diagram
\[\xymatrix@!=0.1pc{G \ar[ddd]_{\bm{F}_{G}} \ar[dr]^{\bm{F}_{G}} \ar@{=}[rrr]^-{} &&& G \ar@{-->}[ddd]_(0.6){\bm{F}_{G}} \ar@{=}[dr]^{} &\\
& \D(G) \ar[rrr]^{\del_{G}} \ar[ddd]^(0.4){\D(\bm{F}_{G})} && \ulltwocell<\omit>{\bm{\zeta}} & G \ar[ddd]^(0.4){\bm{F}_{G}} \\
&&&&  \\
\D(G) \ar[dr]_{\ce_{G}} \ar@{==}[rrr]^{} & \uultwocell<\omit>{\bm{\theta}} && \D(G) \ar@{==}[dr]_(0.4){} &\\
& \D^{2}(G) \ar[rrr]_{\del_{\D(G)}} &&& \D(G) }\]
commutes meaning 
\begin{equation}\label{leftidentity}
\begin{array}{c}
\bm{F}_{G}\bm{\zeta} \cdot \del_{\D(G)} \bm{\theta} = \iota_{\bm{F}_{G}}.
\end{array}
\end{equation}

\item a diagram
\[\xymatrix@!=0.1pc{G \ar[ddd]_{\bm{F}_{G}} \ar@{=}[dr]^{} \ar[rrr]^-{\bm{F}_{G}} &&& \D(G) \ar@{-->}[ddd]_(0.5){\D(\bm{F}_{G})} \ar@{=}[dr]^{} &\\
& G \ar[rrr]^{\bm{F}_{G}} \ar[ddd]^(0.4){\bm{F}_{G}} &&& \D(G) \ar@{=}[ddd] \\
&& \ultwocell<\omit>{\bm{\theta}} &&  \\
\D(G) \ar@{=}[dr]_{} \ar@{-->}[rrr]^{\ce_{G}} &&& \D^{2}(G) \ar@{-->}[dr]_(0.4){\D(\ce_{G})} & \uultwocell<\omit>{\D(\bm{\zeta})\,\,\,\,\,\,\,\,\,\,} \\
& \D(G) \ar@{=}[rrr]_{} &&& \D(G) }\]
commutes meaning 
\begin{equation}\label{rightidentity}
\begin{array}{c}
\D(\bm{\zeta}) \bm{F}_{G} \cdot \D(\del_{G}) \bm{\theta} = \iota_{\bm{F}_{G}}.
\end{array}
\end{equation}

\item a diagram
\[\xymatrix@!=0.1pc{G \ar[ddd]_{\bm{F}_{G}} \ar[dr]^{\bm{F}_{G}} \ar[rrr]^-{\bm{F}_{G}} &&& \D(G) \ar@{-->}[ddd]_(0.6){\D(\bm{F}_{G})} \ar[dr]^{\D(\bm{F}_{G})} &\\
& \D(G) \ar[rrr]^{\ce_{G}} \ar[ddd]^(0.4){\D(\bm{F}_{G})} && \ulltwocell<\omit>{\bm{\theta}} & \D^{2}(G)  \ar[ddd]^(0.4){\D^{2}(\bm{F}_{G})} \\
&& \ultwocell<\omit>{\bm{\theta}} &&  \\
\D(G) \ar[dr]_{\ce_{G}} \ar@{-->}[rrr]^{\ce_{G}} & \uultwocell<\omit>{\bm{\theta}} && \D^{2}(G) \ar@{-->}[dr]_(0.4){\D(\ce_{G})} & \uultwocell<\omit>{\D(\bm{\theta})\,\,\,\,\,\,\,\,\,\,} \\
& \D^{2}(G) \ar[rrr]_{\ce_{\D(G)}} &&& \D^{3}(G) }\]
commutes meaning that 
\begin{equation}\label{coassociativity}
\begin{array}{c}
\D^{2}(\bm{F}_{G}) \bm{\theta}  \cdot \ce_{\D(G)} \bm{\theta}= \D(\bm{\theta}) \bm{F}_{G} \cdot \D(\ce_{G}) \bm{\theta}. 
\end{array}
\end{equation}
\end{itemize}
We say that the a colax $\D$-coalgebra $(G,\bm{F}_{G},\bm{\eta},\bm{\theta})$ is normal if $\bm{\zeta}$ in (\ref{zeta}) is identity and we call it a strict $\D$-coalgebra if $\bm{\theta}$ in (\ref{theta}) is identity too. If both $\bm{\zeta}$ and $\bm{\theta}$ are natural isomorphisms we call $(G,\bm{F}_{G},\bm{\eta},\bm{\theta})$ a pseudo $\D$-coalgebra. Finally, if we reverse the direction of 2-cells $\bm{\zeta}$ in (\ref{zeta}) and $\bm{\theta}$ in (\ref{theta}) we talk about lax $\D$-coalgebras.
\end{definition}

\begin{definition}\label{colaxDcoalgebramor}
Let $(G,\bm{F}_{G},\bm{\zeta},\bm{\theta})$ and $(U,\bm{F}_{U},\bm{\zeta},\bm{\theta})$ be two colax $\D$-coalgebras
\[\xymatrix@!=3pc{\A \ar[r]^-{\bm{F}_{0}} \ar[d]_-{G} & (\X,G) \ar@{=}[d]^{} \\
\X \ar[r]_-{\bm{F}_{1}} & (\X,G) \ultwocell<\omit>{\bm{\varphi}} }
\hspace{2cm}
\xymatrix@!=3pc{\E \ar[r]^-{\bm{F}_{0}} \ar[d]_-{U} & (\Y,U) \ar@{=}[d]^{} \\
\Y \ar[r]_-{\bm{F}_{1}} & (\Y,U) \ultwocell<\omit>{\bm{\upsilon}} }\]
A colax morphism $(\bm{K},\bm{\beta}) \maps (G,\bm{F}_{G},\bm{\eta},\bm{\theta}) \to (U,\bm{F}_{U},\bm{\eta},\bm{\theta})$ consists of the following data:
\begin{itemize}
\item a 1-cell $\bm{K}=(K_{1},\kappa,K_{0}) \maps G \to U$ in $\Cat^{2}_{c}$ as in a diagram
\begin{equation}\label{kappa}
\begin{array}{c}
\xymatrix@!=3pc{\A \ar[d]_{G} \ar[r]^{K_{0}} & \E \ar[d]^-{U}\\
\X \ar[r]_-{K_{1}} & \Y \ultwocell<\omit>{\kappa} }
\end{array}
\end{equation}
\item a 2-cell $\bm{\beta} \maps \D(\bm{K})\bm{F}_{G} \Rightarrow \bm{F}_{U}\bm{K}$ in $\Cat^{2}_{c}$ as in a diagram
\begin{equation}\label{beta}
\begin{array}{c}
\xymatrix@!=0.1pc{\A \ar[ddd]_{G} \ar[dr]^{\bm{F}_{0}} \ar[rrr]^-{K_{0}} &&& \E \ar@{-->}[ddd]_(0.6){U} \ar[dr]^{\bm{U}_{0}} &\\
& (\X,G) \ar[rrr]^{\D(\bm{K})} \ar@{=}[ddd] && \ulltwocell<\omit>{\bm{\beta}_{0}\,\,\,\,} & (\Y,U) \ar@{=}[ddd] \\
&& \ultwocell<\omit>{\kappa} &&  \\
\X \ar[dr]_(0.4){\bm{F}_{1}} \ar@{-->}[rrr]^{K_{1}} & \uultwocell<\omit>{\bm{\varphi}} && \Y \ar@{-->}[dr]_(0.4){\bm{U}_{1}} & \uultwocell<\omit>{\bm{\upsilon}} \\
& (\X,G) \ar[rrr]_{\D(\bm{K})} && \ulltwocell<\omit>{\bm{\beta}_{1}\,\,\,\,} & (\Y,U) }
\end{array}
\end{equation}
\end{itemize}
such that the following axioms are satisfied:
\begin{itemize}
\item a diagram 
\begin{equation}\label{colaxlaxmor1}
\begin{array}{c}
\xymatrix@!=0.1pc{G \ar[ddd]_{\bm{F}_{G}} \ar@{=}[dr] \ar[rrr]^-{\bm{K}} &&& U \ar@{-->}[ddd]_(0.6){\bm{F}_{U}} \ar@{=}[dr] &\\
& G \ar[rrr]^{\bm{K}} \ar@{=}[ddd] &&& U \ar@{=}[ddd] \\
&& \ultwocell<\omit>{\bm{\beta}} &&  \\
\D(G) \ar[dr]_(0.4){\del_{G}} \ar@{-->}[rrr]^{\D(\bm{K})} & \uultwocell<\omit>{\bm{\zeta}} && \D(U) \ar@{-->}[dr]_(0.4){\del_{U}} & \uultwocell<\omit>{\bm{\zeta}} \\
& G \ar[rrr]_{\bm{K}} &&& U }
\end{array}
\end{equation}
of 2-cells in $\Cat^{2}_{c}$ commute.

\item a diagram 
\begin{equation}\label{colaxlaxmor2}
\begin{array}{c}
\xymatrix@!=0.1pc{G \ar[ddd]_{\bm{F}_{G}} \ar[dr]^{\bm{F}_{G}} \ar[rrr]^-{\bm{K}} &&& U \ar@{-->}[ddd]_(0.6){\bm{F}_{U}} \ar[dr]^{\bm{F}_{U}} &\\
& \D(G) \ar[rrr]^{\D(\bm{K})} \ar[ddd]^(0.4){\D(\bm{F}_{G})} && \ulltwocell<\omit>{\bm{\beta}} & \D(U)  \ar[ddd]^(0.4){\D(\bm{F}_{U})} \\
&& \ultwocell<\omit>{\bm{\beta}} &&  \\
\D(G) \ar[dr]_(0.4){\ce_{G}} \ar@{-->}[rrr]^{\D(\bm{K})} & \uultwocell<\omit>{\bm{\theta}} && \D(U) \ar@{-->}[dr]_(0.4){\ce_{U}} \ultwocell<\omit>{\D(\bm{\beta})\,\,\,\,\,\,} & \uultwocell<\omit>{\bm{\theta}} \\
& \D^{2}(G) \ar[rrr]_{\D^{2}(\bm{K})} &&& \D^{2}(G) }
\end{array}
\end{equation}
of 2-cells in $\Cat^{2}_{c}$ commute.
\end{itemize}
\end{definition}

\noindent Finally, we define a notion of a transformation between colax $\D$-coalgebra morphisms:

\begin{definition}\label{colaxDcoalgebratrans}
Let $(\bm{K},\bm{\beta}), (\bm{L},\bm{\gamma}) \maps (G,\bm{F}_{G},\bm{\zeta},\bm{\theta}) \to (U,\bm{F}_{U},\bm{\zeta},\bm{\theta})$ be two colax $\D$-coalgebra morphism where $(\bm{L},\bm{\gamma})$ is given by a data in a commutative diagram
\[\xymatrix@!=0.4pc{\A \ar[ddd]_{G} \ar[dr]^{\bm{F}_{0}} \ar[rrr]^-{L_{0}} &&& \E \ar@{-->}[ddd]_(0.6){U} \ar[dr]^{\bm{U}_{0}} &\\
& (\X,G) \ar[rrr]^{\D(\bm{L})} \ar@{=}[ddd] && \ulltwocell<\omit>{\bm{\gamma}_{0}\,\,\,\,} & (\Y,U) \ar@{=}[ddd] \\
&& \ultwocell<\omit>{\lambda} &&  \\
\X \ar[dr]_(0.4){\bm{F}_{1}} \ar@{-->}[rrr]^{L_{1}} & \uultwocell<\omit>{\bm{\varphi}} && \Y \ar@{-->}[dr]_(0.4){\bm{U}_{1}} & \uultwocell<\omit>{\bm{\upsilon}} \\
& (\X,G) \ar[rrr]_{\D(\bm{L})} && \ulltwocell<\omit>{\bm{\gamma}_{1}\,\,\,\,} & (\Y,U) }\]
A colax $\D$-coalgebra transformation $\bm{\tau} \maps (\bm{K},\bm{\beta}) \Rightarrow (\bm{L},\bm{\gamma})$ is a 2-cell $\bm{\tau} \maps \bm{K} \Rightarrow \bm{L}$ in $\Cat^{2}_{c}$ such that the following diagram
\begin{equation}\label{colaxDcoalgebratransdiag}
\begin{array}{c}
\xymatrix@!=0.4pc{G \ar[ddd]_{\bm{F}_{G}} \ar@{=}[dr]^{} \ar[rrr]^{\bm{K}} & \drrtwocell<\omit>{\bm{\tau}} && U \ar@{-->}[ddd]_{\bm{F}_{U}} \ar@{=}[dr]^{} &\\
& G \ar[rrr]^{\bm{L}} \ar[ddd]^(0.4){\bm{F}_{G}} &&& U \ar[ddd]_{\bm{F}_{U}} \\
&& \ultwocell<\omit>{\bm{\theta}} &&  \\
\D(G) \ar@{=}[dr]_{} \ar@{-->}[rrr]^{\D(\bm{K})} & \drrtwocell<\omit>{\,\,\,\,\,\,\,\,\D(\bm{\tau})} && \D(U) \ar@{==}[dr]_(0.4){}  \ultwocell<\omit>{\bm{\theta}} & \\
& \D(G) \ar[rrr]_{\D(\bm{L})} &&& \D(U) }
\end{array}
\end{equation}
commutes.
\end{definition}

\newpage

Let us investigate the conditions that a colax $\D$-coalgebra needs to satisfy.  The 2-cell $\bm{\zeta}$ in (\ref{zeta}) with the two components $\zeta^{0}$ and $\zeta^{1}$ on the top and bottom squares
\begin{equation}\label{unitcoalgebra}
\begin{array}{c}
\xymatrix@!=0.5pc{\A \ar@{=}[rrr]^{} \ar[ddd]_{G} \ar[dr]_{(H,\chi,Q)} &&& \A \ar@{==}[ddd]^{} \ar@{=}[dr]^{} &\\
& (\X,G) \ar[rrr]^{\de_{0}} \ar@{=}[ddd]_{} && \ulltwocell<\omit>{\zeta^{0}} & \A \ar[ddd]^{G} \\
&&&&\\
\X \ar[dr]_{(C,\eta,K)} \ar@{==}[rrr]^{} & \uultwocell<\omit>{(\omega,\varepsilon)\,\,\,\,\,\,} && \X \ar@{==}[dr]_(0.3){}  \ultwocell<\omit>{\delta_{G}\,\,\,\,\,} & \\
& (\X,G) \ar[rrr]_{\de_{1}} && \ulltwocell<\omit>{\zeta^{1}} & \X}
\end{array}
\end{equation}
of (\ref{unitcoalgebra}) respectively, consists of the following data:
\begin{itemize}
\item [1)] The natural transformation $\zeta^{0} \maps Q \Rightarrow I_{\A}$ where the functor $\bm{F}_{0} \maps \A \to (\X,G)$ is completely determined a triple $\bm{F}_{0}=(H,\chi,Q)$ where $H \maps \A \to \X$ and $Q \maps \A \to \A$ are functors and $\chi \maps H \Rightarrow GQ$ a natural transformation 
\item [2)] The natural transformation $\eta^{1} \maps C \Rightarrow I_{\X}$ where the functor $\bm{F}_{1} \maps \X \to (\X,G)$  is completely determined a triple $\bm{F}_{1}=(C,\eta,K)$ where $C \maps \X \to \X$ and $K \maps \X \to \A$ are functors and $\eta \maps C \Rightarrow GK$ a natural transformation

\item [3)] As a consequence of the first two items the 1-cell $\bm{F}_{G}=(\bm{F}_{1},\bm{\varphi}, \bm{F}_{0}) \maps G \to \D(G)$ has the from $((C,\eta,K),(\omega,\varepsilon), (H,\chi,Q))$ as in the following diagram
\begin{equation}\label{colaxcoalg}
\begin{array}{c}
\xymatrix@!=3pc{\A \ar[r]^-{(H,\chi,Q)} \ar[d]_-{G} & (\X,G) \ar@{=}[d]^{} \\
\X \ar[r]_-{(C,\eta,K)} & (\X,G) \ultwocell<\omit>{(\omega,\varepsilon)\,\,\,\,\,\,\,\,\,} }
\end{array}
\end{equation}
where $\omega \maps CG \Rightarrow H$ and $\varepsilon \maps KG \Rightarrow Q$ are natural transformations such that a diagram
\begin{equation}\label{colaxcoh}
\begin{array}{c}
\xymatrix@!=3pc{CG(A) \ar[d]_{\eta_{G(A)}} \ar[r]^-{\omega_{A}} & H(A) \ar[d]^-{\chi_{A}}\\
GKG(A) \ar[r]_-{G(\varepsilon_{A})} & GQ(A)}
\end{array}
\end{equation}
commutes for any object $A$ in $\A$.

\item [4)] $\del \circ \bm{\varphi}$ is natural transformation whose component indexed by $A$ is defined by a commutative diagram
\[\xymatrix@!=0.1pc{\de_{1}(C,\eta,K)G(A) \ar[rrrrr]^-{\de_{1}((\omega,\varepsilon)_{A})} \ar[dd]_{\delta_{(C,\eta,K)(A)}} &&&&&\de_{1}(H,\chi,Q)G(A) \ar[dd]^-{\delta_{(H,\chi,Q)(A)}}\\
&&&&& \\
G\de_{0}(C,\eta,K)G(A) \ar[rrrrr]_-{G\de_{0}((\omega,\varepsilon)_{A})} &&&&& G\de_{0}(H,\chi,Q)(A) }\]
which is identical to (\ref{colaxcoh}).  The last diagram defines a pasting composition of 1-cells
\[\xymatrix@!=2pc{\A \ar[r]^-{(H,\chi,Q)} \ar[d]_-{G} & (\X,G) \ar@{=}[d]^{} \ar[r]^-{\de_{0}} & \A  \ar[d]^{G} \\
\X \ar[r]_-{(C,\eta,K)} & (\X,G) \ultwocell<\omit>{(\omega,\varepsilon)\,\,\,\,\,\,\,\,\,} \ar[r]_-{\de_{1}} & \X \ultwocell<\omit>{\delta_{G}} }\]
which is a (vertical) domain of the 2-cell $\bm{\zeta}$ as in a commutative diagram
\[\xymatrix@!=0.5pc{\A \ar[dd]_{G} \ar@/^1.25pc/[rr]^{I_{\A}} \ar@/_1.25pc/[rr]_{Q} && \A \ar[dd]^{G} \lltwocell<\omit>{\zeta^{0}\,\,\,} \\
&&  \\
\X \ar@/^1.25pc/@{-->}[rr]^{I_{\X}} \ar@/_1.25pc/[rr]_{C} && \X \lltwocell<\omit>{\zeta^{1}\,\,\,} \uulltwocell<\omit>{\del \circ \bm{\varphi}\,\,\,\,\,\,\,} }\]
providing a compatibility condition for the two components of $\bm{\zeta}$ given by the identity
\begin{equation}\label{zetacomp}
\begin{array}{c}
G \zeta^{0} \cdot (\del \circ \bm{\varphi}) = \zeta^{1} G.
\end{array}
\end{equation}
It follows from  (\ref{zetacomp}) that for any object $A$ in $\A$ there is a commutative diagram
\begin{equation}\label{colaxcoh2}
\begin{array}{c}
\xymatrix@!=0.1pc{CG(A) \ar[dd]_{\eta_{G(A)}} \ar[dr]^{\omega_{A}} \ar[rrr]^-{\zeta^{1}_{G(A)}} &&& G(A) \ar@{=}[dr] \ar@{==}[dd]_(0.4){} & \\
& H(A) \ar[dd]^(0.6){\chi_{A}} \ar[rrr]_-{G(\zeta^{0}_{A})\chi_{A}} &&& G(A) \ar@{=}[dd]^{} \\
GKG(A)  \ar[dr]_-{G(\varepsilon_{A})} \ar@{-->}[rrr]^-{G(\zeta^{0}_{A}\varepsilon_{A})} &&& G(A) \ar@{==}[dr]_(0.4){} & \\
& GQ(A) \ar[rrr]_{G(\zeta^{0}_{A})}  &&& G(A)}
\end{array}
\end{equation}
\end{itemize}

The functor $D({\bm{F}_{G}})$ is induced from the universal property of the comma square
\[\xymatrix@!=0.1pc{(\X,G) \ar[rrr]^{\de_{0}} \ar[ddd]_{\de_{1}} \ar[dr]^{D({\bm{F}_{G}})} &&& \A \ar@{-->}[ddd]^{G} \ar[dr]^{(H,\chi,Q)} &\\
& (\X,G)^{2} \ar[rrr]^{\de_{0}} \ar[ddd]_{\de_{1}} &&& (\X,G) \ar@{=}[ddd]^{} \\
&& \ultwocell<\omit>{\delta_{G}\,\,\,\,} &&\\
\X \ar[dr]_(0.3){(C,\eta,K)} \ar@{==}[rrr]^{} &&& \X \ar@{-->}[dr]_(0.3){(C,\eta,K)} \ultwocell<\omit>{\delta_{I_{(\X,G)}}\,\,\,\,\,\,\,\,\,\,\,\,\,\,\,\,\,} & \uultwocell<\omit>{(\omega,\varepsilon)\,\,\,\,\,\,\,\,\,\,\,\,}   \\ 
& (\X,G) \ar@{=}[rrr]_{} &&& (\X,G) }\]
on the front side. The value of $D({\bm{F}_{G}})$ on an object $(X,f,A)$ in $(\X,G)$ is given by
\begin{equation}\label{DFGob}
D({\bm{F}_{G}})(X,f,A)=((C,\eta,K)(X),(\omega,\varepsilon)_{A}(C,\eta,K)(f),(H,\chi,Q)(A))
\end{equation} 
and any morphism $(x,a) \maps (X,f,A) \to (X',f',A')$ in $(\X,G)$ is sent to a morphism $D({\bm{F}_{G}})(x,a)$  in $(\X,G)^{2}$ represented the horizontal boundary of the following diagram 
\begin{equation}\label{DFG}
\begin{array}{c}
\hspace{-1cm}\xymatrix@!=1.5pc@R1pc@C1pc{C(X) \ar[ddd]_{\eta_{X}} \ar[rrr]^-{C(x)} \ar[dr]^{C(f)} &&& C(X') \ar@{-->}[ddd]_{\eta_{X'}} \ar[dr]^{C(f')} &&\\
& CG(A) \ar[ddd]_{\eta_{G(A)}} \ar[dr]^(0.6){\omega_{A}} \ar[rrr]^-{CG(a)} &&& CG(A') \ar@{-->}[ddd]_{\eta_{G(A')}} \ar[dr]^{\omega_{A'}} &\\
&& H(A) \ar[ddd]_{\chi_{A}} \ar[rrr]^{H(a)} &&& H(A') \ar[ddd]^{\chi_{A'}}  \\
GK(X) \ar[dr]_{GK(f)} \ar@{-->}[rrr]^-{GK(x)} &&& GK(X') \ar@{-->}[dr]_(0.4){GK(f')} && \\
& GKG(A) \ar[dr]_{G(\varepsilon_{A})} \ar@{-->}[rrr]_{GKG(a)} &&& GKG(A') \ar@{-->}[dr]_{G(\varepsilon_{A'})} & \\
&& GQ(A)  \ar[rrr]_{GQ(a)} &&& GQ(A') }
\end{array}
\end{equation}
taking the object  $[(C(X),\eta_{X},K(X)), (\omega_{A}C(f),\varepsilon_{A}K(f)), (H(A),\chi_{A},Q(A))]$ on the left to an object $[(C(X'),\eta_{X'},K(X')), (\omega_{A'}C(f'),\varepsilon_{A'}K(f')), (H(A'),\chi_{A'},Q(A'))]$ on the right of (\ref{DFG}). The 2-cell $\bm{\theta}=(\bm{\theta}^{1},\bm{\theta}^{0})$ in (\ref{theta}) is a pair of natural transformations
\[\xymatrix@!=1.5pc@R1pc@C1pc{\A \ar[ddd]_{G} \ar[rrr]^{(H,\chi,Q)} \ar[dr]^{(H,\chi,Q)} &&& (\X,G) \ar@{==}[ddd]^{} \ar[dr]^{D({\bm{F}_{G}})} &\\
& (\X,G) \ar[rrr]_{\n_{(\X,G)}} \ar@{=}[ddd]_{} && \ulltwocell<\omit>{\bm{\theta}^{0}} & (\X,G)^{2} \ar@{=}[ddd]_{}  \\
&& \ultwocell<\omit>{(\omega,\varepsilon)\,\,\,\,\,}  &&  \\
\X \ar[dr]_{(C,\eta,K)} \ar@{-->}[rrr]^{(C,\eta,K)} & \uultwocell<\omit>{(\omega,\varepsilon)\,\,\,} && (\X,G) \ar@{-->}[dr]_(0.4){D({\bm{F}_{G}})} &\\
& (\X,G) \ar[rrr]_{\n_{(\X,G)}} && \ulltwocell<\omit>{\bm{\theta}^{1}} & (\X,G)^{2} }\]
respectively,  which induces the following identity of natural transformations
\begin{equation}\label{theta2}
\begin{array}{c}
\bm{\theta}^{0} \cdot \n_{(\X,G)}(\omega,\varepsilon)= D({\bm{F}_{G}})(\omega,\varepsilon) \cdot \bm{\theta}^{1} G.
\end{array}
\end{equation}
The component $\bm{\theta}^{1}_{X}$ of the natural transformation $\bm{\theta}^{1}$ indexed by an object $X$ in $\X$
\begin{equation}\label{theta1}
\begin{array}{c}
\xymatrix@!=1.5pc@R1pc@C1pc{C(X) \ar[ddd]_{\eta_{X}} \ar[rrr]^-{\theta^{100}_{X}} \ar[ddrr]^{1_{C(X)}} &&& CC(X) \ar@{-->}[ddd]_{\eta_{C(X)}} \ar[dr]^{C(\eta_{X})} &&\\
&&&& CGK(X) \ar@{-->}[ddd]_{\eta_{GK(X)}} \ar[dr]^{\omega_{K(X)}} &\\
&& C(X) \ar[ddd]_{\eta_{X}} \ar[rrr]^{\theta^{101}_{X}} &&& HK(X) \ar[ddd]^{\chi_{K(X)}}  \\
GK(X) \ar[ddrr]_{G(1_{K(X)})} \ar@{-->}[rrr]^-{G(\theta^{110}_{X})} &&& GKC(X) \ar@{-->}[dr]_(0.3){GK(\eta_{X})} && \\
&&&& GKGK(X) \ar@{-->}[dr]_(0.3){G(\varepsilon_{K(X)})} & \\
&& GK(X) \ar[rrr]_{G(\theta^{111}_{X})} &&& GQK(X) }
\end{array}
\end{equation}
is a horizontal boundary of (\ref{theta1}) given by a quadruple $[(\theta^{100}_{X},\theta^{110}_{X},) (\theta^{101}_{X},\theta^{111}_{X})]$. It takes an object $[(C(X),\eta_{X},K(X)),(1_{C(X)},1_{K(X)}),(C(X),\eta_{X},K(X))]$ to an object $[(CC(X),\eta_{C(X)},KC(X)),(\omega_{K(X)}C(\eta_{X}),\varepsilon_{K(X)}K(\eta_{X})),(HK(X),\chi_{K(X)},QK(X))]$.  

The identities
\begin{equation}\label{theta1eq1}
\eta_{C(X)}\theta^{100}_{X}=\G(\theta^{110}_{X})\eta_{X},  
\end{equation}
\begin{equation}\label{theta1eq2}
\chi_{K(X)}\theta^{101}_{X} = G(\theta^{111}_{X})\eta_{X}, 
\end{equation}
\begin{equation}\label{theta1eq3}
\omega_{K(X)} C(\eta_{X}) \theta^{100}_{X} = \theta^{101}_{X}, 
\end{equation}
\begin{equation}\label{theta1eq4}
\varepsilon_{K(X)}K(\eta_{X})\theta^{110}_{X} = \theta^{111}_{X}
\end{equation}
hold by definition of (\ref{theta1}).
The component $\bm{\theta}^{0}_{A}$ of the natural transformation $\bm{\theta}^{0}$
\begin{equation}\label{theta0}
\begin{array}{c}
\xymatrix@!=1.5pc@R1pc@C1pc{H(A) \ar[ddd]_{\chi_{A}} \ar[rrr]^-{\theta^{000}_{A}} \ar[ddrr]^{1_{H(A)}} &&& CH(A) \ar@{-->}[ddd]_{\eta_{H(A)}} \ar[dr]^{C(\chi_{A})} &&\\
&&&& CGQ(A) \ar@{-->}[ddd]^{\eta_{GQ(A)}} \ar[dr]^{\omega_{Q(A)}} &\\
&& H(A) \ar[ddd]_{\chi_{A}} \ar[rrr]^{\theta^{001}_{A}} &&& HQ(A) \ar[ddd]^{\chi_{Q(A)}}  \\
GQ(A) \ar[ddrr]_{G(1_{Q(A)})} \ar@{-->}[rrr]^-{G(\theta^{010}_{A})} &&& GKH(A) \ar@{-->}[dr]_(0.3){GK(\chi_{A})} && \\
&&&& GKGQ(A) \ar@{-->}[dr]_(0.3){G(\varepsilon_{Q(A)})} & \\
&& GQ(A) \ar[rrr]_{G(\theta^{011}_{A})} &&& GQQ(A) }
\end{array}
\end{equation}
is a horizontal boundary of (\ref{theta0}) given by a quadruple $[(\theta^{000}_{A},\theta^{010}_{A}),(\theta^{001}_{A},\theta^{011}_{A})]$.  It takes an object $[(H(A),\chi_{A},Q(A)),(1_{H(A)},1_{Q(A)}),(H(A),\chi_{A},Q(A))]$ on the left to an object $[(H(A),\eta_{H(A)},KH(A)),(\omega_{A}\chi_{A},\varepsilon_{A}K(\chi_{A})),(H(A),\chi_{A},Q(A))]$ on the right (\ref{theta0}) inducing identities 
\begin{equation}\label{theta0eq1}
\eta_{H(A)}\theta^{000}_{A}=G(\theta^{010}_{A})\chi_{A},
\end{equation}
\begin{equation}\label{theta0eq2}
\chi_{Q(A)}\theta^{001}_{A} = G(\theta^{011}_{A})\chi_{A},
\end{equation}
\begin{equation}\label{theta0eq3}
\omega_{Q(A)}C(\chi_{A})\theta^{000}_{A} = \theta^{001}_{A}, 
\end{equation}
\begin{equation}\label{theta0eq4}
\varepsilon_{Q(A)}K(\chi_{A})\theta^{010}_{A} = \theta^{011}_{A}. 
\end{equation}
\newpage
The component $D({\bm{F}_{G}})(\omega,\varepsilon)_{A}$ indexed by $A$ in $\A$ is the right half of a diagram
\begin{equation}\label{thetacomp1}
\begin{array}{c}
\xymatrix@!=1.25pc@R1pc@C1pc{CG(A) \ar[ddd]_{\eta_{G(A)}} \ar[rrr]^-{\theta^{100}_{G(A)}} \ar[ddrr]^{1_{CG(A)}} &&& 
CCG(A) \ar@{-->}[ddd]_{\eta_{CG(A)}} \ar[dr]^{C(\eta_{G(A)})} \ar[rrr]^-{C(\omega_{A})}  &&& 
CH(A) \ar@{-->}[ddd]^{\eta_{H(A)}} \ar[dr]^{C(\chi_{A})} &&\\
&&&& CGKG(A) \ar@{-->}[ddd]_{\eta_{GKG(A)}}  \ar[rrr]^-{CG(\varepsilon_{A})} \ar[dr]^(0.6){\omega_{KG(A)}} &&& CGQ(A) \ar@{-->}[ddd]^{\eta_{GQ(A)}} \ar[dr]^{\omega_{Q(A)}}&\\
&& CG(A) \ar[ddd]_(0.6){\eta_{G(A)}} \ar[rrr]^{\theta^{101}_{G(A)}} &&& HKG(A) \ar[ddd]_{\chi_{KG(A)}} \ar[rrr]^{H(\varepsilon_{A})} &&& HQ(A) \ar[ddd]^{\chi_{Q(A)}}  \\
GKG(A) \ar[ddrr]_{G(1_{KG(A)})} \ar@{-->}[rrr]^-{G(\theta^{110}_{G(A)})} &&& GKCG(A) \ar@{-->}[dr]_(0.3){GK(\eta_{G(A)})} \ar@{--}[rrr]^-{GK(\omega_{A})} &&& GKH(X) \ar@{-->}[dr]_(0.3){GK(\chi_{A})} && \\
&&&& GKGKG(A) \ar@{-->}[dr]_(0.3){G(\varepsilon_{KG(A)})} \ar@{-->}[rrr]_-{GKG(\varepsilon_{A})} &&& 
GKGQ(A) \ar@{-->}[dr]_(0.3){G(\varepsilon_{Q(A)})} & \\
&& GKG(A) \ar[rrr]_{G(\theta^{111}_{G(A)})} &&& GQKG(A) \ar[rrr]_{GQ(\varepsilon_{A})} &&& GQ(A) }
\end{array}
\end{equation}
whose left half is $\bm{\theta}^{1}_{G(A)}$ and the component $\n_{(\X,G)}((\omega,\varepsilon)_{A})$ is a left half of a diagram
\begin{equation}\label{thetacomp2}
\begin{array}{c}
\xymatrix@!=1.25pc@R1pc@C1pc{CG(A) \ar[ddd]_{\eta_{G(A)}} \ar[rrr]^-{\omega_{A}} \ar[ddrr]^{1_{G(A)}} &&& H(A) \ar@{-->}[ddd]_{\chi_{A}} \ar[ddrr]^{1_{H(A)}} \ar[rrr]^-{\theta^{000}_{A}}  &&& 
CH(A) \ar@{-->}[ddd]_{\eta_{H(A)}} \ar[dr]^{C(\chi_{A})} &&\\
&&&&&&& CGQ(A) \ar@{-->}[ddd]^{\eta_{G(A)}} \ar[dr]^{\omega_{Q(A)}} &\\
&& CG(A) \ar[ddd]_{\eta_{G(A)}} \ar[rrr]^{\omega_{A}} &&& H(A) \ar[ddd]_{\chi_{A}} \ar[rrr]^{\theta^{001}_{A}} &&& HQ(A) \ar[ddd]^{\chi_{Q(A)}}  \\
GKG(A) \ar[ddrr]_{G(1_{KG(A)})} \ar@{-->}[rrr]^-{G(\varepsilon_{A})} &&& GQ(A) \ar@{-->}[ddrr]_{G(1_{A})} \ar@{-->}[rrr]^-{G(\theta^{010}_{A})} &&& 
GKH(A) \ar@{-->}[dr]_(0.3){GK(\chi_{A})} && \\
&&&&&&& GKGQ(A) \ar@{-->}[dr]_(0.3){G(\varepsilon_{Q(A)})} & \\
&& GKG(A) \ar[rrr]_{G(\varepsilon_{A})} &&& GQ(A)  \ar[rrr]_{G(\theta^{011}_{A})} &&& GQQ(A) }
\end{array}
\end{equation}
whose right half is the component $\bm{\theta}^{0}_{A}$. Then the identity (\ref{theta2}) means that the boundaries of the above two diagrams are the same,  giving us the following identities
\begin{equation}\label{thetaAcompeq1}
C(\omega_{A})\theta^{100}_{G(A)} =  \theta^{000}_{A} \omega_{A},
\end{equation}
\begin{equation}\label{thetaAcompeq2}
K(\omega_{A})\theta^{110}_{G(A)}=\theta^{010}_{A}\varepsilon_{A},
\end{equation}
\begin{equation}\label{thetaAcompeq3}
H(\varepsilon_{A})\theta^{101}_{G(A)} = \theta^{001}_{A}\omega_{A}, 
\end{equation}
\begin{equation}\label{thetaAcompeq4}
Q(\varepsilon_{A})\theta^{111}_{G(A)} = \theta^{011}_{A}\varepsilon_{A}.
\end{equation}
\newpage
The corresponding diagram for the axiom (\ref{leftidentity}) becomes the commuting diagram
\begin{equation}\label{leftunitCat}
\begin{array}{c}
\xymatrix@!=0.5pc@R1pc@C1.25pc{\A \ar[rrrrrrrr]^-{I_{\A}} \ar[ddddddddd]_{(H,\chi,Q)}
\ar@{-->}[dddrrr]_(0.6){G} \ar[ddrrrr]^{(H,\chi,Q)} &&&&&&&& \A \ar[ddrrrr]^-{I_{\A}} \ar@{-->}[dddl]_{G} \ar@{-->}[ddddddddd]^(0.55){(H,\chi,Q)} &&&& \\
&&&&&&&&&& \\
&&&& (\X,G) \ar[rrrrrrrr]^(0.55){\de_{0}} \ar[ddddddddd]_(0.45){D(\bm{F}_{G})} \ar@{==}[ddr]^{} && \ulltwocell<\omit>{(\omega,\varepsilon)\,\,\,\,\,\,\,\,\,} & 
\uulltwocell<\omit>{\bm{\zeta}_{0}\,\,} &&&&& \A \ar[ddddddddd]^{(H,\chi,Q)} \ar@{-->}[ddlll]_{G} \\
&&& \X \ar@{-->}[dddd]_(0.5){(C,\eta,K)} \ar@{-->}[drr]_(0.3){(C,\eta,K)} \ar@{-->}[rrrr]^{I_{\X}} & \ultwocell<\omit>{(\omega,\varepsilon)\,\,\,\,\,\,\,\,\,\,} &&& 
\X \ar@{-->}[dddd]_{(C,\eta,K)} \ar@{-->}[drr]^{I_{\X}} &&&&& \\
&& \dltwocell<\omit>{(\omega,\varepsilon)\,\,\,\,\,\,\,\,\,\,} &&& (\X,G) \ar@{-->}[dddd]^{D(\bm{F}_{G})} \ar@{-->}[rrrr]_{\de_{1}} && \ulltwocell<\omit>{\bm{\zeta}_{1}\,\,} &
& \X \ar@{-->}[dddd]^{(C,\eta,K)} \uulltwocell<\omit>{\delta_{G}\,\,\,\,\,} &&&\\
&&&&&&&&&&&& \\
&&& \ulltwocell<\omit>{\bm{\theta}_{0}\,\,} && \ulltwocell<\omit>{\bm{\theta}_{1}\,\,} &&& \ultwocell<\omit>{(\omega,\varepsilon)\,\,\,\,\,\,\,\,\,\,} &&&& \\
&&& (\X,G) \ar@{-->}[drr]_(0.4){\n_{(\X,G)}} \ar@{-->}[rrrr]^{I_{(\X,G)}} &&&& 
(\X,G) \ar@{-->}[drr]^{I_{(\X,G)}} \ar@{==}[ddr]^{} &&& \urtwocell<\omit>{\,\,\,\,\,\,\,\,\,\,(\omega,\varepsilon)} &&&\\
&&&&& (\X,G)^{2} \ar@{-->}[rrrr]_{\de_{1}} &&&& (\X,G) \ar@{==}[dddrrr]^-{} &&& \\
(\X,G) \ar@{-->}[rrrrrrrr]^{I_{(\X,G)}} \ar[ddrrrr]_-{\n_{(\X,G)}} \ar@{==}[uurrr]^{} &&&&&&&& (\X,G) \ar@{-->}[ddrrrr]_-{I_{(\X,G)}} &&&& \\
&&&&&&& \urtwocell<\omit>{\,\,\,\,\,\,\,\,\, \,\,\,\,\,\,\delta_{I_{(\X,G)}}} &&&& \\
&&&& (\X,G)^{2} \ar[rrrrrrrr]_{\de_{0}} \ar@{==}[uuur]_-{} &&&&&&&& (\X,G) }
\end{array}
\end{equation}
with two identities
\begin{equation}\label{thetaleq1}
\begin{array}{c}
(C,\eta,K)\bm{\zeta}^{1} \cdot \de_{1}\bm{\theta}^{1} =  \iota_{(C,\eta,K)}
\end{array}
\end{equation}
\begin{equation}\label{thetaleq0}
\begin{array}{c}
(H,\chi,Q)\bm{\zeta}^{0} \cdot \de_{0}\bm{\theta}^{0} =  \iota_{(H,\chi,Q)}
\end{array}
\end{equation}
meaning that for any objects $X$ in $\X$ and $A$ in $\A$ the following two diagrams 
\[\xymatrix@!=0.25pc{C(X) \ar[dd]_{\eta_{X}} \ar[rr]^-{\theta^{100}_{X}} && CC(X) \ar[dd]^{\eta_{C(X)}} \ar[rr]^-{C(\zeta^{1}_{X})} && C(X) \ar[dd]^{\eta_{X}} \\
&&&& \\
GK(X) \ar[rr]_-{G(\theta^{110}_{X})} && GKC(X) \ar[rr]_-{GK(\zeta^{1}_{X})} && GK(X) }
\hspace{1cm}
\xymatrix@!=0.25pc{H(A) \ar[dd]_{\chi_{A}} \ar[rr]^-{\theta^{001}_{A}} && HQ(A) \ar[dd]^{\chi_{Q(A)}} \ar[rr]^-{H(\zeta^{0}_{A})} && H(A) \ar[dd]^{\chi_{A}} \\
&&&& \\
GQ(A) \ar[rr]_-{G(\theta^{011}_{A})} && GQQ(A) \ar[rr]_-{GQ(\zeta^{0}_{A})} && GQ(A) }\] 
are equal to identities $(1_{C(X)},1_{K(X)}) \maps (C(X),\eta_{X},K(X)) \to (C(X),\eta_{X},K(X))$ and $(1_{H(A)},1_{Q(A)}) \maps (H(A),\chi_{A},Q(A)) \to (H(A),\chi_{A},Q(A))$ respectively, implying the identities
\begin{equation}\label{leftuniteqX1}
C(\zeta^{1}_{X})\theta^{100}_{X} = 1_{C(X)},
\end{equation}
\begin{equation}\label{leftuniteqX2}
K(\zeta^{1}_{X})\theta^{110}_{X} = 1_{K(X)},
\end{equation}
\begin{equation}\label{leftuniteqA1}
H(\zeta^{0}_{A})\theta^{001}_{A} = 1_{H(A)},
\end{equation}
\begin{equation}\label{leftuniteqA2}
Q(\zeta^{0}_{A})\theta^{011}_{A} = 1_{Q(A)}.
\end{equation}
The axiom (\ref{rightidentity}) is equivalent to a diagram
\begin{equation}\label{rightunitCat}
\begin{array}{c}
\xymatrix@!=0.5pc@R1pc@C1.25pc{\A \ar[rrrrrrrr]^-{(H,\chi,Q)} \ar[ddddddddd]_{I_{\A}}
\ar@{-->}[dddrrr]_(0.6){G} \ar[ddrrrr]^{I_{\A}} &&&&&&&& (\X,G) \ar[ddrrrr]^-{I_{(\X,G)}} \ar@{==}[dddl]^{} \ar@{-->}[ddddddddd]^(0.5){D(\bm{F}_{G})} &&&& \\
&&&&&&&&&& \\
&&&& \A \ar[rrrrrrrr]^{(H,\chi,Q)} \ar[ddddddddd]_(0.45){D(\bm{F}_{G})} \ar@{-->}[ddr]^(0.6){G} && \ulltwocell<\omit>{(\omega,\varepsilon)\,\,\,\,\,\,\,\,\,} &&&&&& 
(\X,G) \ar[ddddddddd]^{I_{(\X,G)}} \ar@{==}[ddlll]_{} \\
&&& \X \ar@{-->}[dddd]_{(C,\eta,K)} \ar@{-->}[drr]_(0.3){I_{\X}} \ar@{-->}[rrrr]^{(C,\eta,K)} & \ultwocell<\omit>{(\omega,\varepsilon)\,\,\,\,\,\,\,\,\,\,} &&& 
(\X,G) \ar@{-->}[dddd]_{D(\bm{F}_{G})} \ar@{-->}[drr]^{I_{(\X,G)}} &&&&& \\
&& \dltwocell<\omit>{(\omega,\varepsilon)\,\,\,\,\,\,\,\,\,\,} &&& \X \ar@{-->}[dddd]^{(C,\eta,K)} \ar@{-->}[rrrr]_{(C,\eta,K)} &&&& 
(\X,G) \ar@{-->}[dddd]^{I_{(\X,G)}} \uulltwocell<\omit>{(\omega,\varepsilon)\,\,\,\,\,\,\,} &&& \\
&&&&& \uulltwocell<\omit>{\bm{\theta}_{0}\,\,} &&&&&&& \\
&&&&&& \uulltwocell<\omit>{\bm{\theta}_{1}\,\,} &&& \ulltwocell<\omit>{D(\bm{\zeta})\,\,\,\,\,\,\,\,\,\,} && \ulltwocell<\omit>{D(\bm{\zeta})\,\,\,\,\,\,\,\,\,} & \\
&&& (\X,G) \ar@{-->}[drr]_(0.4){I_{(\X,G)}} \ar@{-->}[rrrr]^{\n_{(\X,G)}} && \ultwocell<\omit>{(\omega,\varepsilon)\,\,\,\,\,\,\,\,\,} && (\X,G)^{2} \ar@{-->}[drr]^{D(\del_{G})} \ar@{==}[ddr]^{} &&&&&&\\
&&&&& (\X,G) \ar@{-->}[rrrr]_{I_{(\X,G)}} &&&& (\X,G) \ar@{==}[dddrrr]^-{}  &&& \\
(\X,G) \ar@{-->}[rrrrrrrr]^{\n_{(\X,G)}} \ar[ddrrrr]_-{I_{(\X,G)}} \ar@{==}[uurrr]^{} &&&&&&&& (\X,G)^{2} \ar@{-->}[ddrrrr]_-{D(\del_{G})} &&&& \\
&&&&&&&&&&& \\
&&&& (\X,G) \ar[rrrrrrrr]_{I_{(\X,G)^{2}}} \ar@{==}[uuur]_-{} &&&&&&&& (\X,G) }
\end{array}
\end{equation}
giving
\begin{equation}\label{thetareq1}
\begin{array}{c}
D(\bm{\zeta})(C,\eta,K) \cdot D(\del_{G})\bm{\theta}^{1} =  \iota_{(C,\eta,K)}
\end{array}
\end{equation}
\begin{equation}\label{thetareq0}
\begin{array}{c}
D(\bm{\zeta})(H,\chi,Q) \cdot D(\del_{G})\bm{\theta}^{0} =  \iota_{(H,\chi,Q)}
\end{array}
\end{equation}
Equations (\ref{thetareq1}) and (\ref{thetareq0}) means that a composition of morphisms as in a diagram
\[\hspace{-1cm}\xymatrix@!=1pc@R1pc@C1.25pc{C(X) \ar[dd]_{\eta_{X}} \ar[rrr]^-{\theta^{100}_{X}} &&& CC(X) \ar[dd]^{\eta_{C(X)}G(\varepsilon_{K(X)}K(\eta_{X}))} \ar[rrr]^-{\zeta^{1}_{C(X)}} &&& 
C(X) \ar[dd]^{\eta_{X}} \\
&&&&&& \\
GK(X) \ar[rrr]_-{G(\theta^{111}_{X})} &&& GQK(X) \ar[rrr]_-{G(\zeta^{0}_{K(X)})} &&& GK(X) }
\xymatrix@!=1pc@R1pc@C1.25pc{H(A) \ar[dd]_{\chi_{A}} \ar[rrr]^-{\theta^{000}_{A}} &&& CH(A) \ar[dd]^{\chi_{Q(A)}} \ar[rrr]^-{\zeta^{1}_{H(A)}} &&& H(A) \ar[dd]^{\chi_{A}} \\
&&&&&& \\
GQ(A) \ar[rrr]_-{G(\theta^{011}_{A})} &&& GQQ(A) \ar[rrr]_-{G(\zeta^{0}_{Q(A)})} &&& GQ(A) }\] 
are equal to the identities $(1_{C(X)},1_{K(X)}) \maps (C(X),\eta_{X},K(X)) \to (C(X),\eta_{X},K(X))$ and $(1_{H(A)},1_{Q(A)}) \maps (H(A),\chi_{A},Q(A)) \to (H(A),\chi_{A},Q(A))$ in $(\X,G) $ respectively. Therefore:
\begin{equation}\label{rightuniteqX1}
\zeta^{1}_{C(X)} \theta^{100}_{X} = 1_{C(X)},
\end{equation}
\begin{equation}\label{rightuniteqX2}
\zeta^{0}_{K(X)} \theta^{111}_{X} = 1_{K(X)},
\end{equation}
\begin{equation}\label{rightuniteqA1}
\zeta^{1}_{H(A)} \theta^{000}_{A} = 1_{H(A)},
\end{equation}
\begin{equation}\label{rightuniteqA2}
\zeta^{0}_{Q(A)} \theta^{011}_{A} = 1_{Q(A)}.
\end{equation}
\newpage
The last condition coming out from (\ref{rightunitCat}) is the compatibility of (\ref{thetareq1}) and (\ref{thetareq0}) expressed by the equation
\[\xymatrix@!=1pc{(C,\eta,K)G \ar@{=}[d]_{} \ar@2[rrrr]^-{(\omega,\varepsilon)} &&&& (H,\chi,Q) \ar@{=}[d]_{}  \\
D(\del_{G})\n_{(\X,G)}(C,\eta,K)G \ar@2[d]_{D(\del_{G})\theta^{1}G}  \ar@2[rrrr]^-{D(\del_{G})\n_{(\X,G)}(\omega,\varepsilon)} &&&& 
D(\del_{G})\n_{(\X,G)}(H,\chi,Q) \ar@2[d]^{D(\del_{G})\theta^{0}} \\
D(\del_{G})D(\bm{F}_{G})(C,\eta,K)G  \ar@2[d]_{D(\bm{\zeta})(C,\eta,K)G} \ar@2[rrrr]_{D(\del_{G})D(\bm{F}_{G})(\omega,\varepsilon)} &&&& 
D(\del_{G})D(\bm{F}_{G})(H,\chi,Q)  \ar@2[d]^{D(\bm{\zeta})(H,\chi,Q)}  \\
(C,\eta,K)G \ar@2[rrrr]_-{(\omega,\varepsilon)} &&&& (H,\chi,Q) }\]

\[\xymatrix@!=1.5pc@R1pc@C1pc{CG(A) \ar[ddd]_{\varrho_{G(A)}} \ar[rrr]^-{\theta^{100}_{G(A)}} \ar[ddrr]^{1_{CG(A)}} &&& 
CCG(A) \ar@{-->}[ddd]_{\varrho_{CG(A)}} \ar[dr]^{C(\varrho_{G(A)})} \ar[rrr]^-{\zeta^{1}_{CG(A)}}  &&& 
CG(A) \ar@{-->}[ddd]^{\varrho_{G(A)}} \ar[dr]^{\varrho_{G(A)}} &&\\
&&&& CGKG(A) \ar@{-->}[ddd]_{\varrho_{GKG(A)}} \ar[rrr]^-{\zeta^{1}_{GKG(A)}} \ar[dr]^(0.6){\omega_{KG(A)}} &&& GKG(A) \ar@{-->}[ddd]^{1_{GKG(A)}} \ar[dr]^{1_{GKG(A)}} &\\
&& CG(A) \ar[ddd]_(0.6){\varrho_{G(A)}} \ar[rrr]^{\theta^{101}_{G(A)}} &&& HKG(A) \ar[ddd]_{\chi_{KG(A)}} \ar[rrr]^{H(\varepsilon_{(A)})} &&& GKG(A) \ar[ddd]^{\chi_{Q(A)}}  \\
GKG(A) \ar[ddrr]_{G(1_{KG(A)})} \ar@{-->}[rrr]^-{G(\theta^{110}_{G(A)})} &&& GKCG(A) \ar@{-->}[dr]_(0.3){GK(\varrho_{G(A)})} \ar@{--}[rrr]^-{GK(\zeta_{CG(A)})} &&& GKG(A) \ar@{-->}[dr]_(0.3){G(1_{KG(A)})} && \\
&&&& GKGKG(A) \ar@{-->}[dr]_(0.3){G(\varepsilon_{KG(A)})} \ar@{-->}[rrr]_-{G(\zeta^{0}_{KG(A)}\varepsilon_{KG(A)})} &&& 
GKG(A) \ar@{-->}[dr]_(0.3){G(1_{KG(A)})} & \\
&& GKG(A) \ar[rrr]_{G(\theta^{111}_{KG(A)})} &&& GQKG(A) \ar[rrr]_{G(\zeta^{0}_{KG(A)})} &&& GKG(A) }\]
The horizontal boundary of the following diagram 
\[\xymatrix@!=1.5pc@R1pc@C1pc{C(X) \ar[ddd]_{\eta_{X}} \ar[rrr]^-{\theta^{100}_{X}} \ar[ddrr]^{1_{C(X)}} &&& 
CC(X) \ar@{-->}[ddd]_{\eta_{C(X)}} \ar[dr]^{C(\eta_{X})} \ar[rrr]^{\zeta^{1}_{C(X)}} \ar@{-->}[ddd]_{\eta_{C(X)}} &&& 
C(X) \ar@{-->}[ddd]^{\eta_{X}} \ar[dr]^{\eta_{X}} &&\\
&&&& CGK(X) \ar@{-->}[ddd]_{\eta_{GK(X)}} \ar[dr]^{\omega_{K(X)}}  \ar[rrr]^{\zeta^{1}_{GK(X)}} &&& GK(X) \ar@{-->}[ddd]_{1_{GK(X)}} \ar[dr]^{1_{GK(X)}} &\\
&& C(X) \ar[ddd]_{\eta_{C(X)}} \ar[rrr]^-{\theta^{101}_{X}} &&& HK(X) \ar[ddd]_{\chi_{K(X)}} \ar[rrr]^{G(\zeta^{0}_{K(X)})\chi_{K(X)}} &&& GK(X) \ar[ddd]^{1_{GK(X)}}    \\
GK(X) \ar[ddrr]_{G(1_{K(X)})} \ar@{-->}[rrr]^-{G(\theta^{110}_{X})} &&& GKC(X) \ar@{-->}[dr]_(0.3){GK(\eta_{X})}  \ar@{-->}[rrr]^{GK(\zeta^{1}_{X})} &&& 
GK(X) \ar@{-->}[dr]_(0.3){G(1_{K(X)})}  && \\
&&&& GKGK(X) \ar@{-->}[dr]_(0.3){G(\varepsilon_{K(X)})}  \ar@{-->}[rrr]_{G(\zeta^{0}_{K(X)}\varepsilon_{K(X)})} &&& GK(X) \ar@{-->}[dr]_(0.3){G(1_{K(X)})} & \\
&& GK(X) \ar[rrr]_{G(\theta^{111}_{X})} &&& GQK(X) \ar[rrr]_{G(\zeta^{0}_{K(X)})} &&& GK(X) }\]
the top raw of the back rectangle is identity by (\ref{rightuniteqX1}) and the bottom raw of the back rectangle is identity by (\ref{leftuniteqX2})

\newpage
\noindent The objects of the category $((\X,G)^{2})^{2}$ are morphisms in $(\X,G)^{2}$ which we represent by cubes 
\begin{equation}\label{hyperobject}
\begin{array}{c}
\xymatrix@!=0.1pc{X \ar[rrr]^{i} \ar[ddd]_{f} \ar[dr]^{x} &&& Z\ar@{-->}[ddd]^{k} \ar[dr]^{z} &\\
& Y \ar[rrr]^{j} \ar[ddd]_{g} &&& W \ar[ddd]^{l} \\
&&&&  \\
G(A) \ar[dr]_{G(a)} \ar@{-->}[rrr]^{G(d)} && & G(M) \ar@{-->}[dr]_(0.4){G(m)} &   \\
& G(E) \ar[rrr]_{G(e)} &&& G(N) }
\end{array}
\end{equation}
and typical morphism looks like a hypercube
\begin{equation}\label{hypermorphism}
\begin{array}{c}
\xymatrix@!=1pc@R1pc@C1pc{X \ar[rrrrrrrr]^-{h_{100}} \ar[ddddddddd]_{f}
\ar@{-->}[dddrrr]_(0.6){i} \ar[ddrrrr]^-{x} &&&&&&&& X' \ar[ddrrrr]^-{x'} \ar@{-->}[dddl]_{i'}  \ar@{-->}[ddddddddd]^(0.6){f'} &&&& \\
&&&&&&&&&& \\
&&&& Y \ar[rrrrrrrr]^(0.6){h_{101}} \ar[ddddddddd]_{g} \ar@{-->}[ddr]_{j} &&&&&&&& 
Y' \ar[ddddddddd]^{g'} \ar@{-->}[ddlll]_{j'} \\
&&& Z \ar@{-->}[dddd]_{k} \ar@{-->}[drr]_(0.4){z} \ar@{-->}[rrrr]^{h_{000}} &&&& 
Z' \ar@{-->}[dddd]^{k'} \ar@{-->}[drr]^(0.3){z'} &&&&& \\
&&&&& W \ar@{-->}[dddd]_{l} \ar@{-->}[rrrr]_{h_{001}} &&&& W'\ar@{-->}[dddd]^{l'} &&&\\
&&&&&&&&&&&& \\
&&&&&&&&&&&& \\
&&& G(M) \ar@{-->}[drr]_(0.4){G(m)} \ar@{-->}[rrrr]^(0.6){G(h_{010})} &&&& G(M') \ar@{-->}[drr]^{G(m')} &&&&&&\\
&&&&& G(N) \ar@{-->}[rrrr]_{G(h_{011})} &&&& g(N') &&& \\
G(A) \ar@{-->}[rrrrrrrr]^{G(h_{110})} \ar[ddrrrr]_-{G(a)} \ar@{-->}[uurrr]^{G(d)} &&&&&&&& G(A') \ar@{-->}[ddrrrr]_-{G(a')}  \ar@{-->}[uul]^(0.7){G(d')} &&&& \\
&&&&&&&&&&& \\
&&&& G(E) \ar[rrrrrrrr]_{G(h_{111})} \ar@{-->}[uuur]_-{G(e)} &&&&&&&& G(E')  \ar@{-->}[uuulll]_-{G(e')} }
\end{array}
\end{equation}
with the convention that the eight components $(h_{000},h_{001},\ldots h_{111})$ always go from left to right.
\noindent The functor $DID({\bm{F}_{G}}) \maps (\X,G)^{2} \to ((\X,G)^{2})^{2}$ is induced by the universal property of a comma square
\[\xymatrix@!=0.3pc{(\X,G)^{2} \ar[rrr]^{\de_{0}} \ar[ddd]_{\de_{1}} \ar[dr]^{DID({\bm{F}_{G}})} &&& (\X,G) \ar@{==}[ddd]^{} \ar[dr]^{D({\bm{F}_{G}})} &\\
& ((\X,G)^{2})^{2} \ar[rrr]^{\de_{0}} \ar[ddd]_{\de_{1}} &&& (\X,G)^{2} \ar@{=}[ddd]^{} \\
&& \ultwocell<\omit>{\delta_{I_{(\X,G)}}\,\,\,\,\,\,\,\,\,\,\,\,\,\,\,\,\,} &&  \\
(\X,G) \ar[dr]_{D({\bm{F}_{G}})} \ar@{==}[rrr]^{} && & (\X,G) \ar@{-->}[dr]_(0.4){D({\bm{F}_{G}})} \ultwocell<\omit>{\delta_{I_{(\X,G)^{2}}}\,\,\,\,\,\,\,\,\,\,\,\,\,\,\,\,\,\,\,\,} &   \\
& (\X,G)^{2} \ar@{=}[rrr]_{} &&& (\X,G)^{2} }\]
in the front side of the above diagram. It takes an object $[(X,f,A),(x,a),(Y,g,E)]$ as on the left side of (\ref{hyperobject}) to an object in $(\X,G)^{3}$ defined by a formula
\begin{equation}\label{Dhyperobject}
\begin{array}{c}
DID({\bm{F}_{G}})[(X,f,A),(x,a),(Y,g,E)]=\\
=[[D({\bm{F}_{G}})(X,f,A),D({\bm{F}_{G}})(x,a),D({\bm{F}_{G}})(Y,g,E)]]
\end{array}
\end{equation}
and represented by the following diagram
\begin{equation}\label{DIDaction}
\begin{array}{c}
\xymatrix@!=0.25pc{C(X) \ar[rrr]^{\omega_{A}C(f)} \ar[ddd]_{\eta_{X}} \ar[dr]^{C(x)} &&& H(A) \ar@{-->}[ddd]^{\chi_{A}} \ar[dr]^{H(a)} &\\
& C(Y) \ar[rrr]^{\omega_{E}C(g)} \ar[ddd]_{\eta_{Y}} &&& H(E) \ar[ddd]^{\chi_{E}} \\
&&&&  \\
GK(X) \ar[dr]_{GK(x)} \ar@{-->}[rrr]^{G(\varepsilon_{A}K(f))} &&& GQ(A) \ar@{-->}[dr]_(0.4){GQ(a)} &   \\
& GK(Y) \ar[rrr]_{G(\varepsilon_{E}K(g))} &&& GQ(E) }
\end{array}
\end{equation}
The functor $DID({\bm{F}_{G}})$ takes the morphism in $(\X,G)^{2}$ represented by (\ref{hyperobject}) to a hypercube
\begin{equation}\label{DIDhyperaction}
\begin{array}{c}
\xymatrix@!=1pc@R1pc@C1pc{C(X) \ar[rrrrrrrr]^-{C(i)} \ar[ddddddddd]_{\eta_{X}}
\ar@{-->}[dddrrr]_(0.6){\omega_{A}C(f)} \ar[ddrrrr]^-{C(x)} &&&&&&&& C(Z) \ar[ddrrrr]^-{C(z)} \ar@{-->}[dddl]_{\omega_{M}C(k)}  \ar@{-->}[ddddddddd]^(0.6){\eta_{Z}} &&&& \\
&&&&&&&&&& \\
&&&& C(Y) \ar[rrrrrrrr]^(0.6){C(j)} \ar[ddddddddd]_{\eta_{Y}} \ar@{-->}[ddr]^(0.75){\omega_{E}C(g)} &&&&&&&& 
C(W) \ar[ddddddddd]^{\eta_{W}} \ar@{-->}[ddlll]_{\omega_{N}C(l)} \\
&&& H(A) \ar@{-->}[dddd]_{\chi_{A}} \ar@{-->}[drr]_(0.4){H(a)} \ar@{-->}[rrrr]^{H(d)} &&&& 
H(M) \ar@{-->}[dddd]^{\chi_{M}} \ar@{-->}[drr]^(0.6){H(m)} &&&&& \\
&&&&& H(E) \ar@{-->}[dddd]_{\chi_{E}} \ar@{-->}[rrrr]_{H(e)} &&&& H(N) \ar@{-->}[dddd]^{\chi_{N}} &&&\\
&&&&&&&&&&&& \\
&&&&&&&&&&&& \\
&&& GQ(A) \ar@{-->}[drr]_(0.4){GQ(a)} \ar@{-->}[rrrr]^(0.6){GQ(d)} &&&& GQ(M) \ar@{-->}[drr]^{GQ(m)} &&&&&&\\
&&&&& GQ(E) \ar@{-->}[rrrr]_{GQ(e)} &&&& GQ(N) &&& \\
GK(X) \ar@{-->}[rrrrrrrr]^{GK(i)} \ar[ddrrrr]_-{GK(x)} \ar@{-->}[uurrr]^{G(\varepsilon_{A}K(f))} &&&&&&&& GK(Z) \ar@{-->}[ddrrrr]_-{GK(z)}  \ar@{-->}[uul]^(0.7){G(\varepsilon_{M}K(k))} &&&& \\
&&&&&&&&&&& \\
&&&& GK(Y) \ar[rrrrrrrr]_{GK(j)} \ar@{-->}[uuur]_-{G(\varepsilon_{E}K(g))} &&&&&&&& GK(W)  \ar@{-->}[uuulll]_(0.8){G(\varepsilon_{N}K(l))} }
\end{array}
\end{equation}
and it simply means "apply $D({\bm{F}_{G}})$ everywhere".

\noindent The functor $DI(\n_{(\X,G)}) \maps (\X,G)^{2} \to ((\X,G)^{2})^{2}$ is induced by the same universal property
\[\xymatrix@!=0.3pc{(\X,G)^{2} \ar[rrr]^{\de_{0}} \ar[ddd]_{\de_{1}} \ar[dr]^{DI(\n_{(\X,G)})} &&& (\X,G) \ar@{==}[ddd]^{} \ar[dr]^{\n_{(\X,G)}} &\\
& ((\X,G)^{2})^{2} \ar[rrr]^{\de_{0}} \ar[ddd]_{\de_{1}} &&& (\X,G)^{2} \ar@{=}[ddd]^{} \\
&& \ultwocell<\omit>{\delta_{I_{(\X,G)}}\,\,\,\,\,\,\,\,\,\,\,\,\,\,\,\,\,} &&  \\
(\X,G) \ar[dr]_{\n_{(\X,G)}} \ar@{==}[rrr]^{} && & (\X,G) \ar@{-->}[dr]_(0.4){\n_{(\X,G)}} \ultwocell<\omit>{\delta_{I_{(\X,G)^{2}}}\,\,\,\,\,\,\,\,\,\,\,\,\,\,\,\,\,\,\,\,} &   \\
& (\X,G)^{2} \ar@{=}[rrr]_{} &&& (\X,G)^{2} }\]
and it takes an object $[(X,f,A),(x,a),(Y,g,E)]$ as on the left side of (\ref{hyperobject}) to an object defined by
\begin{equation}\label{nhyperobject}
\begin{array}{c}
DI(\n_{(\X,G)})[(X,f,A),(x,a),(Y,g,E)]=\\
=[[\n_{(\X,G)}(X,f,A),\n_{(\X,G)}(x,a),\n_{(\X,G)}(Y,g,E)]]
\end{array}
\end{equation}
and the morphism in $(\X,G)^{2}$ represented by (\ref{hyperobject}) to a morphism in $(\X,G)^{3}$ represented by a hypercube
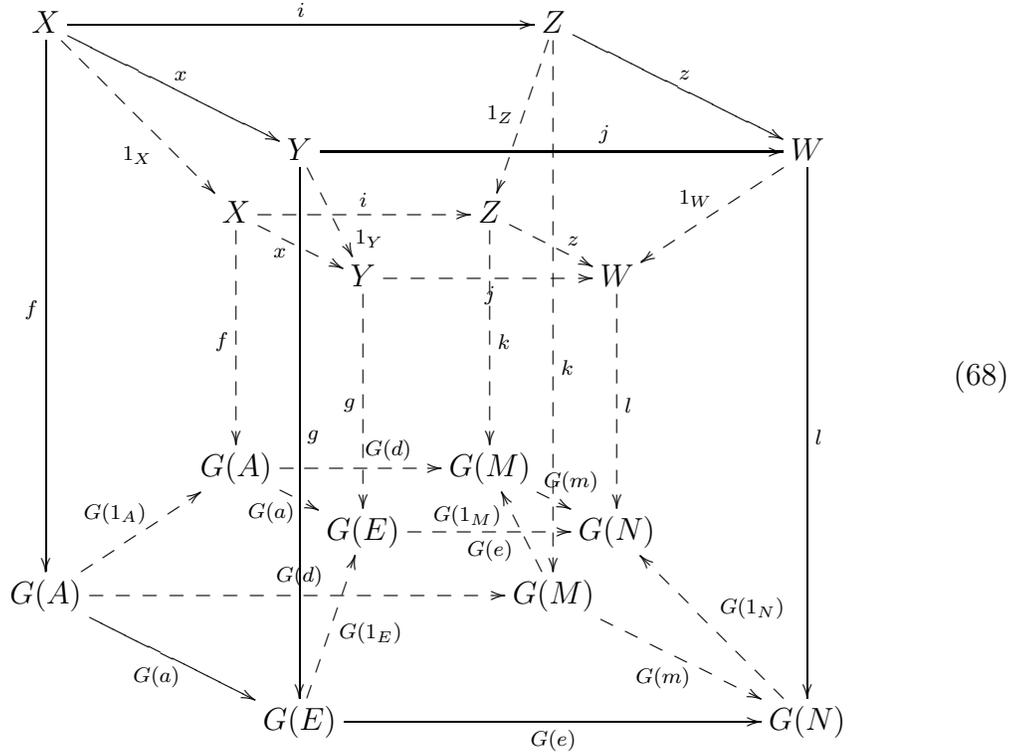
\begin{equation}\label{Dnhyperaction}
\begin{array}{c}
\xymatrix@!=1pc@R1pc@C1pc{X \ar[rrrrrrrr]^-{i} \ar[ddddddddd]_{f}
\ar@{-->}[dddrrr]_(0.6){1_{X}} \ar[ddrrrr]^-{x} &&&&&&&& Z \ar[ddrrrr]^-{z} \ar@{-->}[dddl]_{1_{Z}}  \ar@{-->}[ddddddddd]^(0.6){k} &&&& \\
&&&&&&&&&& \\
&&&& Y \ar[rrrrrrrr]^(0.6){j} \ar[ddddddddd]^{g} \ar@{-->}[ddr]^(0.75){1_{Y}} &&&&&&&& 
W \ar[ddddddddd]^{l} \ar@{-->}[ddlll]_{1_{W}} \\
&&& X \ar@{-->}[dddd]_{f} \ar@{-->}[drr]_(0.4){x} \ar@{-->}[rrrr]^{i} &&&& 
Z \ar@{-->}[dddd]^{k} \ar@{-->}[drr]^(0.6){z} &&&&& \\
&&&&& Y \ar@{-->}[dddd]_{g} \ar@{-->}[rrrr]_{j} &&&& W \ar@{-->}[dddd]^{l} &&&\\
&&&&&&&&&&&& \\
&&&&&&&&&&&& \\
&&& G(A) \ar@{-->}[drr]_(0.4){G(a)} \ar@{-->}[rrrr]^(0.6){G(d)} &&&& G(M) \ar@{-->}[drr]^{G(m)} &&&&&&\\
&&&&& G(E) \ar@{-->}[rrrr]_{G(e)} &&&& G(N) &&& \\
G(A) \ar@{-->}[rrrrrrrr]^{G(d)} \ar[ddrrrr]_-{G(a)} \ar@{-->}[uurrr]^{G(1_{A})} &&&&&&&& G(M) \ar@{-->}[ddrrrr]_-{G(m)} \ar@{-->}[uul]^(0.7){G(1_{M})} &&&& \\
&&&&&&&&&&& \\
&&&& G(E) \ar[rrrrrrrr]_{G(e)} \ar@{-->}[uuur]_-{G(1_{E})} &&&&&&&& G(N) \ar@{-->}[uuulll]_(0.5){G(1_{N})} }
\end{array}
\end{equation}
and it simply means "apply $\n_{(\X,G)}$ everywhere".
\newpage
\noindent The functor $\n_{(\X,G)^{2}} \maps (\X,G)^{2} \to ((\X,G)^{2})^{2}$ takes the morphism (\ref{hyperobject}) in $(\X,G)^{2}$ to a morphism in $(\X,G)^{3}$ given by hypercube
\begin{equation}\label{xihyperaction}
\begin{array}{c}
\xymatrix@!=1pc@R1pc@C1pc{X \ar[rrrrrrrr]^-{i} \ar[ddddddddd]_{f}
\ar@{-->}[dddrrr]_(0.6){x} \ar[ddrrrr]^-{1_{X}} &&&&&&&& Z \ar[ddrrrr]^-{1_{Z}} \ar@{-->}[dddl]_{z}  \ar@{-->}[ddddddddd]^(0.6){k} &&&& \\
&&&&&&&&&& \\
&&&& X \ar[rrrrrrrr]^(0.6){i} \ar[ddddddddd]^{f} \ar@{-->}[ddr]^(0.75){x} &&&&&&&& 
Z \ar[ddddddddd]^{k} \ar@{-->}[ddlll]_{z} \\
&&& Y \ar@{-->}[dddd]_{g} \ar@{-->}[drr]_(0.4){1_{Y}} \ar@{-->}[rrrr]^{j} &&&& 
W \ar@{-->}[dddd]^{l} \ar@{-->}[drr]^(0.6){1_{W}} &&&&& \\
&&&&& Y \ar@{-->}[dddd]_{g} \ar@{-->}[rrrr]_{j} &&&& W \ar@{-->}[dddd]^{l} &&&\\
&&&&&&&&&&&& \\
&&&&&&&&&&&& \\
&&& G(E) \ar@{-->}[drr]_(0.4){G(1_{E})} \ar@{-->}[rrrr]^(0.6){G(e)} &&&& G(N) \ar@{-->}[drr]^{G(1_{N})} &&&&&&\\
&&&&& G(E) \ar@{-->}[rrrr]_{G(e)} &&&& G(N) &&& \\
G(A) \ar@{-->}[rrrrrrrr]^{G(d)} \ar[ddrrrr]_-{G(1_{A})} \ar@{-->}[uurrr]^{G(a)} &&&&&&&& G(M) \ar@{-->}[ddrrrr]_-{G(1_{M})} \ar@{-->}[uul]^(0.7){G(m)} &&&& \\
&&&&&&&&&&& \\
&&&& G(A) \ar[rrrrrrrr]_{G(d)} \ar@{-->}[uuur]_-{G(a)} &&&&&&&& G(M) \ar@{-->}[uuulll]_(0.5){G(m)} }
\end{array}
\end{equation}
When the underlying 2-functor of the 2-comonad $\D$ acts on $\bm{\theta}$ one obtains the natural transformation $\D(\bm{\theta})$ as on the right side of a diagram
\[\hspace{-0.5cm}\xymatrix@!=0.5pc{\A \ar[ddd]_{G} \ar[rrr]^{(H,\chi,Q)} \ar[dr]^{(H,\chi,Q)} &&& (\X,G) \ar@{==}[ddd]^{} \ar[dr]^{D({\bm{F}_{G}})} &\\
& (\X,G) \ar[rrr]_{\n_{(\X,G)}} \ar@{=}[ddd]_{} && \ulltwocell<\omit>{\bm{\theta}^{0}} & (\X,G)^{2} \ar@{=}[ddd]_{}  \\
&& \ultwocell<\omit>{(\omega,\varepsilon)\,\,\,\,\,}  &&  \\
\X \ar[dr]_{(C,\eta,K)} \ar@{-->}[rrr]^{(C,\eta,K)} & \uultwocell<\omit>{(\omega,\varepsilon)\,\,\,} && (\X,G) \ar@{-->}[dr]_(0.4){D({\bm{F}_{G}})} &\\
& (\X,G) \ar[rrr]_{\n_{(\X,G)}} && \ulltwocell<\omit>{\bm{\theta}^{1}} & (\X,G)^{2} }
\xymatrix@!=0.1pc{ && \\
&& \\
\ar@{..>}[rr]^{\D} && \\
&&}
\xymatrix@!=0.5pc{(\X,G) \ar@{=}[ddd]_{} \ar[rrr]^{D(\bm{F}_{G})} \ar[dr]^{D(\bm{F}_{G})} &&& (\X,G)^{2} \ar@{==}[ddd]^{} \ar[dr]^{DID({\bm{F}_{G}})} &\\
& (\X,G)^{2} \ar[rrr]_{D(\n_{(\X,G)})} \ar@{=}[ddd]_{} && \ulltwocell<\omit>{D(\bm{\theta}),\,\,\,\,\,\,} & (\X,G)^{3} \ar@{=}[ddd]_{}  \\
&&&&  \\
(\X,G) \ar[dr]_{D(\bm{F}_{G})} \ar@{-->}[rrr]^{D(\bm{F}_{G})} &&& (\X,G)^{2} \ar@{-->}[dr]_(0.4){DID({\bm{F}_{G}})} &\\
& (\X,G)^{2} \ar[rrr]_{D(\n_{(\X,G)})} && \ulltwocell<\omit>{D(\bm{\theta}),\,\,\,\,\,\,\,\,} & (\X,G)^{3} }\]
whose component $D(\bm{\theta})_{(X,f,A)}$ indexed by an object $(X,f,A)$ in $(\X,G)$ is a diagram
\begin{equation}\label{Dtheta}
\begin{array}{c}
\hspace{-1cm}\xymatrix@!=0.6pc@R1pc@C1pc{C(X) \ar[rrrrrrrr]^-{\theta^{100}_{X}} \ar[ddddddddd]_{\eta_{X}} \ar@{-->}[dddrrr]_(0.6){1_{C(X)}} \ar[ddrrrr]^-{\omega_{A}C(f)} &&&&&&&& 
CC(X) \ar[ddrrrr]^-{C(\omega_{A}C(f))} \ar@{-->}[dddl]_{\omega_{K(X)}C(\eta_{X})}  \ar@{-->}[ddddddddd]^(0.6){\eta_{C(X)}} &&&& \\
&&&&&&&&&& \\
&&&& H(A) \ar[rrrrrrrr]^(0.6){\theta^{000}_{A}} \ar[ddddddddd]_{\chi_{A}} \ar@{-->}[ddr]^(0.75){1_{H(A)}} &&&&&&&& 
CH(A) \ar[ddddddddd]^{\eta_{H(A)}} \ar@{-->}[ddlll]_{\omega_{Q(A)}C(\chi_{A})} \\
&&& C(X) \ar@{-->}[dddd]_{\eta_{X}} \ar@{-->}[drr]_(0.4){\omega_{A}C(f)} \ar@{-->}[rrrr]^{\theta^{101}_{X}} &&&& 
HK(X) \ar@{-->}[dddd]_{\chi_{K(X)}} \ar@{-->}[drr]^(0.5){H(\varepsilon_{A}K(f))} &&&&& \\
&&&&& H(A) \ar@{-->}[dddd]_{\chi_{A}} \ar@{-->}[rrrr]_{\theta^{001}_{A}} &&&& HQ(A) \ar@{-->}[dddd]^{\chi_{Q(A)}} &&&\\
&&&&&&&&&&&& \\
&&&&&&&&&&&& \\
&&& GK(X) \ar@{-->}[drr]_(0.4){G(\varepsilon_{A}K(f))} \ar@{-->}[rrrr]^-{G(\theta^{111}_{X})} &&&& GQK(X) \ar@{-->}[drr]^{GQ(\varepsilon_{A}K(f))} &&&&&&\\
&&&&& GQ(A) \ar@{-->}[rrrr]_{G(\theta^{011}_{A})} &&&& GQQ(A) &&& \\
GK(X) \ar@{-->}[rrrrrrrr]^{G(\theta^{110}_{X})} \ar[ddrrrr]_-{G(\varepsilon_{A}K(f))} \ar@{-->}[uurrr]^{G(1_{K(X)})} &&&&&&&& 
GKC(X) \ar@{-->}[ddrrrr]_-{GK(\omega_{A}C(f))}  \ar@{-->}[uul]^(0.8){G(\varepsilon_{K(X)}K(\eta_{X}))} &&&& \\
&&&&&&&&&&& \\
&&&& GQ(A) \ar[rrrrrrrr]_{G(\theta^{010}_{A})} \ar@{-->}[uuur]_-{G(1_{Q(A)})} &&&&&&&& GKH(A) \ar@{-->}[uuulll]_(0.7){G(\varepsilon_{Q(A)}K(\chi_{A}))} }
\end{array}
\end{equation}
The coassociativity axiom (\ref{coassociativity}) for a colax $\D$-coalgebra means that a diagram
\[\xymatrix@!=0.6pc@R1pc@C1pc{\A \ar[rrrrrrrr]^-{(H,\chi,I_{\A})} \ar[ddddddddd]_{(H,\chi,I_{\A})}
\ar@{-->}[dddrrr]_(0.6){G} \ar[ddrrrr]^{(H,\chi,I_{\A})} &&&&&&&& (\X,G) \ar[ddrrrr]^-{D(\bm{F}_{G})} \ar@{==}[dddl]^{} \ar@{-->}[ddddddddd]^(0.65){D(\bm{F}_{G})} &&&& \\
&&&&&&&&&& \\
&&&& (\X,G) \ar[rrrrrrrr]^{\n_{(\X,G)}} \ar[ddddddddd]^(0.45){D(\bm{F}_{G})} \ar@{==}[ddr]^{} && \ulltwocell<\omit>{(\omega,\varepsilon)\,\,\,\,\,\,\,\,\,} & \uulltwocell<\omit>{\bm{\theta}_{0}\,\,} &&&&& (\X,G)^{2} \ar[ddddddddd]^{DID(\bm{F}_{G})} \ar@{==}[ddlll]_{} \\
&&& \X \ar@{-->}[dddd]_(0.5){(C,\eta,K)} \ar@{-->}[drr]_(0.3){(C,\eta,K)} \ar@{-->}[rrrr]^{(C,\eta,K)} & \ultwocell<\omit>{(\omega,\varepsilon)\,\,\,\,\,\,\,\,\,\,} &&& 
(\X,G) \ar@{-->}[dddd]_{D(\bm{F}_{G})} \ar@{-->}[drr]^{D(\bm{F}_{G})}&&&&& \\
&& \dltwocell<\omit>{(\omega,\varepsilon)\,\,\,\,\,\,\,\,\,\,} &&& (\X,G) \ar@{-->}[dddd]^{D(\bm{F}_{G})} \ar@{-->}[rrrr]_{D({\bm{F}_{G}})} && \ulltwocell<\omit>{\bm{\theta}_{1}\,\,} && (\X,G)^{2} \ar@{-->}[dddd]^{DID(\bm{F}_{G})}&&&\\
&&&&& \uulltwocell<\omit>{\bm{\theta}_{0}\,\,} &&&&&&& \\
&&& \ulltwocell<\omit>{\bm{\theta}_{0}\,\,} && \ulltwocell<\omit>{\bm{\theta}_{1}\,\,} & \uulltwocell<\omit>{\bm{\theta}_{1}\,\,} &&& \ulltwocell<\omit>{D(\bm{\theta})\,\,\,\,\,\,\,\,\,\,} && \ulltwocell<\omit>{D(\bm{\theta})\,\,\,\,\,\,\,\,\,} & \\
&&& (\X,G) \ar@{-->}[drr]_(0.4){\n_{(\X,G)}} \ar@{-->}[rrrr]^{\n_{(\X,G)}} &&&& (\X,G)^{2} \ar@{-->}[drr]^{D(\n_{(\X,G)})} \ar@{==}[ddr]^{} &&&&&&\\
&&&&& (\X,G)^{2} \ar@{-->}[rrrr]_{\n_{(\X,G)^{2}}} &&&& ((\X,G)^{2})^{2} \ar@{==}[dddrrr]^-{}  &&& \\
(\X,G) \ar@{-->}[rrrrrrrr]^{\n_{(\X,G)}} \ar[ddrrrr]_-{\n_{(\X,G)}} \ar@{==}[uurrr]^{} &&&&&&&& (\X,G)^{2} \ar@{-->}[ddrrrr]_-{D(\n_{(\X,G)})} &&&& \\
&&&&&&&&&&& \\
&&&& (\X,G)^{2} \ar[rrrrrrrr]_{\n_{(\X,G)^{2}}} \ar@{==}[uuur]_-{} &&&&&&&& ((\X,G)^{2})^{2} }\]
commutes. Then the identity (\ref{coassociativity}) splits in two equations
\begin{equation}\label{DIDcoassociativity0}
\begin{array}{c}
DID(\bm{F}_{G}) \bm{\theta}_{0}  \cdot \n_{DID(G)} \bm{\theta}_{0}= D(\bm{\theta})_{(H,\chi,I_{\A})} \cdot D(\n_{D(G)}) \bm{\theta}_{0}
\end{array}
\end{equation}
\begin{equation}\label{DIDcoassociativity1}
\begin{array}{c}
DID(\bm{F}_{G}) \bm{\theta}_{1}  \cdot \n_{DID(G)} \bm{\theta}_{1}= D(\bm{\theta})_{(C,\eta,K)} \cdot D(\n_{D(G)}) \bm{\theta}_{1}
\end{array}
\end{equation}

The components $DID(\bm{F}_{G})(\bm{\theta}_{X})$ and $DID(\bm{F}_{G})(\bm{\theta}_{A})$ of the natural transformation $DID(\bm{F}_{G}) \bm{\theta}$ indexed by objects $X$ in $\X$ and $A$ in $\A$ are obtained by applying the functor $DID(\bm{F}_{G})$ to the components (\ref{theta1}) and (\ref{theta0}) respectively.  The eight terms of the identities (\ref{DIDcoassociativity1}) and (\ref{DIDcoassociativity0}) are as follows:

\newpage

\begin{itemize}
\item $\n_{(\X,G)^{2}}(\bm{\theta}_{X})$
\begin{equation}\label{xihyperactionX}
\begin{array}{c}
\xymatrix@!=0.6pc@R1pc@C1.25pc{C(X) \ar[rrrrrrrr]^-{\theta^{100}_{X}} \ar[ddddddddd]_{\eta_{X}}
\ar@{-->}[dddrrr]_(0.6){1_{C(X)}} \ar[ddrrrr]^-{1_{C(X)}} &&&&&&&& CC(X) \ar[ddrrrr]^-{1_{CC(X)}} \ar@{-->}[dddl]_{\omega_{K(X)}C(\eta_{X})} \ar@{-->}[ddddddddd]^(0.6){\eta_{C(X)}} &&&& \\
&&&&&&&&&& \\
&&&& C(X) \ar[rrrrrrrr]^(0.6){\theta^{100}_{X}} \ar[ddddddddd]^{\eta_{X}} \ar@{-->}[ddr]^(0.75){1_{C(X)}} &&&&&&&& 
CC(X) \ar[ddddddddd]^{\eta_{C(X)}} \ar@{-->}[ddlll]_{\omega_{K(X)}C(\eta_{X})} \\
&&& C(X) \ar@{-->}[dddd]_{\eta_{X}} \ar@{-->}[drr]_(0.4){1_{C(X)}} \ar@{-->}[rrrr]^{\theta^{101}_{X}} &&&& 
HK(X) \ar@{-->}[dddd]^{\chi_{K(X)}} \ar@{-->}[drr]^(0.6){1_{HK(X)}} &&&&& \\
&&&&& X \ar@{-->}[dddd]_{\eta_{X}} \ar@{-->}[rrrr]_{\theta^{101}_{X}} &&&& HK(X) \ar@{-->}[dddd]^{\chi_{K(X)}} &&&\\
&&&&&&&&&&&& \\
&&&&&&&&&&&& \\
&&& GK(X) \ar@{-->}[drr]_(0.4){G(1_{K(X)})} \ar@{-->}[rrrr]^-{G(\theta^{111}_{X})} &&&& GQK(X) \ar@{-->}[drr]^{G(1_{QK(X)})} &&&&&&\\
&&&&& GK(X) \ar@{-->}[rrrr]_{G(\theta^{111}_{X})} &&&& GQK(X) &&& \\
GK(X) \ar@{-->}[rrrrrrrr]^{G(\theta^{110}_{X})} \ar[ddrrrr]_-{G(1_{K(X)})} \ar@{-->}[uurrr]^{G(1_{K(X)})} &&&&&&&& 
GKC(X) \ar@{-->}[ddrrrr]_-{G(1_{KC(X)})} \ar@{-->}[uul]^(0.8){G(\varepsilon_{K(X)}K(\eta_{X}))} &&&&\\
&&&&&&&&&&& \\
&&&& GK(X) \ar[rrrrrrrr]_{G(\theta^{110}_{X})} \ar@{-->}[uuur]_-{G(1_{K(X)})} &&&&&&&& GKC(X) \ar@{-->}[uuulll]_(0.5){G(\varepsilon_{K(X)}K(\eta_{X}))} }
\end{array}
\end{equation}

\item $\n_{(\X,G)^{2}}(\bm{\theta}_{A})$
\begin{equation}\label{xihyperactionA}
\begin{array}{c}
\xymatrix@!=0.6pc@R1pc@C1.25pc{H(A) \ar[rrrrrrrr]^-{\theta^{000}_{A}} \ar[ddddddddd]_{\chi_{A}}
\ar@{-->}[dddrrr]_(0.6){1_{H(A)}} \ar[ddrrrr]^-{1_{H(A)}} &&&&&&&& CH(A) \ar[ddrrrr]^-{1_{CH(A)}} \ar@{-->}[dddl]_{\omega_{Q(A)}C(\chi_{A})} \ar@{-->}[ddddddddd]^(0.6){\eta_{H(A)}} &&&& \\
&&&&&&&&&& \\
&&&& H(A) \ar[rrrrrrrr]^(0.6){\theta^{000}_{A}} \ar[ddddddddd]^{\chi_{A}} \ar@{-->}[ddr]^(0.75){1_{H(A)}} &&&&&&&& 
CH(A) \ar[ddddddddd]^{\eta_{H(A)}} \ar@{-->}[ddlll]_{\omega_{Q(A)}C(\chi_{A})} \\
&&& H(A) \ar@{-->}[dddd]_{\chi_{A}} \ar@{-->}[drr]_(0.4){1_{H(A)}} \ar@{-->}[rrrr]^{\theta^{001}_{A}} &&&& 
HQ(A) \ar@{-->}[dddd]^{\chi_{Q(A)}} \ar@{-->}[drr]^(0.6){1_{HQ(A)}} &&&&& \\
&&&&& H(A) \ar@{-->}[dddd]_{\chi_{A}} \ar@{-->}[rrrr]_{\theta^{001}_{A}} &&&& HQ(A) \ar@{-->}[dddd]^{\chi_{Q(A)}} &&&\\
&&&&&&&&&&&& \\
&&&&&&&&&&&& \\
&&& GQ(A) \ar@{-->}[drr]_(0.4){G(1_{Q(A)})} \ar@{-->}[rrrr]^-{G(\theta^{011}_{A})} &&&& GQQ(A) \ar@{-->}[drr]^{G(1_{QQ(A)})} &&&&&&\\
&&&&& GQ(A) \ar@{-->}[rrrr]_-{G(\theta^{011}_{A})} &&&& GQQ(A) &&& \\
GQ(A) \ar@{-->}[rrrrrrrr]^{G(\theta^{010}_{A})} \ar[ddrrrr]_-{G(1_{Q(A)})} \ar@{-->}[uurrr]^{G(1_{Q(A)})} &&&&&&&& 
GKH(A) \ar@{-->}[ddrrrr]_-{G(1_{KH(A)})} \ar@{-->}[uul]^(0.8){G(\varepsilon_{Q(A)}K(\chi_{A}))} &&&&\\
&&&&&&&&&&& \\
&&&& GQ(A) \ar[rrrrrrrr]_{G(\theta^{010}_{A})} \ar@{-->}[uuur]_-{G(1_{Q(A)})} &&&&&&&& GKH(A) \ar@{-->}[uuulll]_(0.5){G(\varepsilon_{Q(A)}K(\chi_{A}))} }
\end{array}
\end{equation}

\item $DID(\bm{F}_{G})(\bm{\theta}_{X})$
\begin{equation}\label{DIDthetaX}
\begin{array}{c}
\hspace{-1cm}\xymatrix@!=0.6pc@R1pc@C1.25pc{CC(X) \ar[rrrrrrrr]^-{C(\theta^{100}_{X})} \ar[ddddddddd]_{\eta_{C(X)}} \ar@{-->}[dddrrr]_(0.6){\omega_{K(X)}C(\eta_{X})} 
\ar[ddrrrr]^-{C(1_{C(X)})} &&&&&&&& 
CCC(X) \ar[ddrrrr]^-{C(\omega_{K(X)}C(\eta_{X}))} \ar@{-->}[dddl]_{\omega_{K(X)}C(\eta_{C(X)})}  \ar@{-->}[ddddddddd]^(0.6){\eta_{CC(X)}} &&&& \\
&&&&&&&&&& \\
&&&& CC(X) \ar[rrrrrrrr]^(0.6){C(\theta^{101}_{X})} \ar[ddddddddd]_{\eta_{C(X)}} \ar@{-->}[ddr]^(0.8){\omega_{K(X)}C(\eta_{X})} &&&&&&&& 
CHK(X) \ar[ddddddddd]^{\eta_{HK(X)}} \ar@{-->}[ddlll]_{\omega_{K(X)}C(\chi_{K(X)})} \\
&&& HK(X) \ar@{-->}[dddd]_{\chi_{K(X)}} \ar@{-->}[drr]_(0.3){H(1_{K(X)})} \ar@{-->}[rrrr]^{H(\theta^{110}_{X})} &&&& 
HKC(X) \ar@{-->}[dddd]_{\chi_{KC(X)}} \ar@{-->}[drr]^(0.5){H(\varepsilon_{K(X)}K(\eta_{X}))} &&&&& \\
&&&&& HK(X) \ar@{-->}[dddd]_{\chi_{K(X)}} \ar@{-->}[rrrr]_{H(\theta^{111}_{X})} &&&& HQK(X) \ar@{-->}[dddd]^{\chi_{QK(X)}} &&&\\
&&&&&&&&&&&& \\
&&&&&&&&&&&& \\
&&& GQK(X) \ar@{-->}[drr]_(0.4){GQ(1_{K(X)})} \ar@{-->}[rrrr]^-{GQ(\theta^{110}_{X})} &&&& GQKC(X) \ar@{-->}[drr]^{GQ(\varepsilon_{K(X)}K(\eta_{X}))} &&&&&&\\
&&&&& GQK(X) \ar@{-->}[rrrr]_-{GQ(\theta^{111}_{X})} &&&& GQQK(X) &&& \\
GKC(X) \ar@{-->}[rrrrrrrr]^{GK(\theta^{100}_{X})} \ar[ddrrrr]_-{GK(1_{C(X)})} \ar@{-->}[uurrr]^{G(\varepsilon_{K(X)}K(\eta_{X}))} &&&&&&&& 
GKCC(X) \ar@{-->}[ddrrrr]_-{GK(\omega_{K(X)}C(\eta_{X}))}  \ar@{-->}[uul]^(0.7){G(\varepsilon_{K(X)}K(\eta_{X}))} &&&& \\
&&&&&&&&&&& \\
&&&& GKC(X) \ar[rrrrrrrr]_{GK(\theta^{101}_{X})} \ar@{-->}[uuur]_-{G(\varepsilon_{K(X)}K(\eta_{X}))} &&&&&&&& GKHK(X)  \ar@{-->}[uuulll]_(0.7){G(\varepsilon_{K(X)}K(\chi_{K(X)})} }\end{array}
\end{equation}

\item $DID(\bm{F}_{G})(\bm{\theta}_{A})$

\begin{equation}\label{DIDthetaA}
\begin{array}{c}
\xymatrix@!=0.6pc@R1pc@C1.25pc{CH(A) \ar[rrrrrrrr]^-{C(\theta^{000}_{A})} \ar[ddddddddd]_{\eta_{H(A)}} \ar@{-->}[dddrrr]_(0.6){\omega_{Q(A)}C(\chi_{A})} 
\ar[ddrrrr]^-{C(1_{H(A)})} &&&&&&&& 
CCH(A) \ar[ddrrrr]^-{C(\omega_{Q(A)}C(\chi_{A}))} \ar@{-->}[dddl]_{\omega_{KH(A)}C(\eta_{H(A)})}  \ar@{-->}[ddddddddd]^(0.6){\eta_{H(A)}} &&&& \\
&&&&&&&&&& \\
&&&& CH(A) \ar[rrrrrrrr]^(0.6){C(\theta^{001}_{A})} \ar[ddddddddd]_{\eta_{H(A)}} \ar@{-->}[ddr]^(0.8){\omega_{Q(A)}C(\chi_{A})} &&&&&&&& 
CHQ(A) \ar[ddddddddd]^{\eta_{H(A)}} \ar@{-->}[ddlll]_{\omega_{QQ(A)}C(\chi_{Q(A)})} \\
&&& HQ(A) \ar@{-->}[dddd]_{\chi_{Q(A)}} \ar@{-->}[drr]_(0.3){H(1_{Q(A)})} \ar@{-->}[rrrr]^-{H(\theta^{010}_{A})} &&&& 
HGKH(A) \ar@{-->}[dddd]_{\chi_{KH(A)}} \ar@{-->}[drr]^(0.5){H(\varepsilon_{Q(A)}K(\chi_{A}))} &&&&& \\
&&&&& HQ(A) \ar@{-->}[dddd]_{\chi_{Q(A)}} \ar@{-->}[rrrr]_{H(\theta^{011}_{A})} &&&& HQQ(A) \ar@{-->}[dddd]^{\chi_{QQ(A)}} &&&\\
&&&&&&&&&&&& \\
&&&&&&&&&&&& \\
&&& GQQ(A) \ar@{-->}[drr]_(0.4){GQ(1_{A})} \ar@{-->}[rrrr]^-{GQ(\theta^{010}_{A})} &&&& GQKH(A) \ar@{-->}[drr]^{GQ(\varepsilon_{Q(A)}K(\chi_{A}))} &&&&&&\\
&&&&& GQQ(A) \ar@{-->}[rrrr]_-{GQ(\theta^{011}_{A})} &&&& GQQQ(A) &&& \\
GKH(A) \ar@{-->}[rrrrrrrr]^{GK(\theta^{000}_{A})} \ar[ddrrrr]_-{GK(1_{H(A)})} \ar@{-->}[uurrr]^{G(\varepsilon_{A}K(\chi_{A}))} &&&&&&&& 
GKCH(A) \ar@{-->}[ddrrrr]_-{GK(\omega_{Q(A)}C(\chi_{A}))}  \ar@{-->}[uul]^(0.8){G(\varepsilon_{KH(A)}K(\eta_{H(A)}))} &&&& \\
&&&&&&&&&&& \\
&&&& GKH(A) \ar[rrrrrrrr]_{GK(\theta^{001}_{A})} \ar@{-->}[uuur]_-{G(\varepsilon_{A}K(\chi_{A}))} &&&&&&&& GKHQ(A)  \ar@{-->}[uuulll]_(0.7){G(\varepsilon_{QQ(A)}K(\chi_{Q(A)}))} }
\end{array}
\end{equation}

\item $D(\n_{(\X,G)})(\bm{\theta}_{X})$
\begin{equation}\label{DnIDX}
\begin{array}{c}
\xymatrix@!=0.6pc@R1pc@C1.25pc{C(X) \ar[rrrrrrrr]^-{\theta^{100}_{X}} \ar[ddddddddd]_{\eta_{X}}\ar@{-->}[dddrrr]_(0.6){1_{C(X)}} \ar[ddrrrr]^-{1_{C(X)}} &&&&&&&& 
CC(X) \ar[ddrrrr]^-{\omega_{K(X)}C(\eta_{X})} \ar@{-->}[dddl]_{1_{CC(X)}} \ar@{-->}[ddddddddd]^(0.6){\eta_{C(X)}} &&&&\\
&&&&&&&&&& \\
&&&& C(X) \ar[rrrrrrrr]^(0.6){\theta^{101}_{X}} \ar[ddddddddd]^{\eta_{X}} \ar@{-->}[ddr]^(0.75){1_{X}} &&&&&&&& 
HK(X) \ar[ddddddddd]^{\chi_{K(X)}} \ar@{-->}[ddlll]_{1_{HK(X)}} \\
&&& C(X) \ar@{-->}[dddd]_{\eta_{X}} \ar@{-->}[drr]_(0.4){1_{C(X)}} \ar@{-->}[rrrr]^{\theta^{100}_{X}} &&&& 
CC(X) \ar@{-->}[dddd]^{\eta_{C(X)}} \ar@{-->}[drr]^(0.6){\omega_{K(X)}C(\eta_{X})} &&&&& \\
&&&&& C(X) \ar@{-->}[dddd]_{\eta_{X}} \ar@{-->}[rrrr]_{\theta^{101}_{X}} &&&& HK(X) \ar@{-->}[dddd]^{\chi_{K(X)}} &&&\\
&&&&&&&&&&&& \\
&&&&&&&&&&&& \\
&&& GK(X) \ar@{-->}[drr]_(0.4){G(1_{K(X)})} \ar@{-->}[rrrr]^-{G(\theta^{110}_{X})} &&&& GKC(X) \ar@{-->}[drr]^{G(\varepsilon_{K(X)})K(\eta_{X}))} &&&&&&\\
&&&&& GK(X) \ar@{-->}[rrrr]_-{G(\theta^{111}_{X})} &&&& GQK(X) &&& \\
GK(X) \ar@{-->}[rrrrrrrr]^{G(\theta^{110}_{X})} \ar[ddrrrr]_-{G(1_{K(X)})} \ar@{-->}[uurrr]^{G(1_{K(X)})} &&&&&&&& 
GKC(X) \ar@{-->}[ddrrrr]_(0.4){G(\varepsilon_{K(X)})K(\eta_{X}))} \ar@{-->}[uul]^(0.8){G(1_{K(X)})} &&&& \\
&&&&&&&&&&& \\
&&&& GK(X) \ar[rrrrrrrr]_{G(\theta^{111}_{X})} \ar@{-->}[uuur]_-{G(1_{K(X)})} &&&&&&&& GQK(X) \ar@{-->}[uuulll]_(0.7){G(1_{QK(X)})} }
\end{array}
\end{equation}

\item $D(\n_{(\X,G)})(\bm{\theta}_{A})$
\begin{equation}\label{DnIDA}
\begin{array}{c}
\xymatrix@!=0.6pc@R1pc@C1.25pc{H(A) \ar[rrrrrrrr]^-{\theta^{000}_{A}} \ar[ddddddddd]_{\chi_{A}}
\ar@{-->}[dddrrr]_(0.6){1_{H(A)}} \ar[ddrrrr]^-{1_{H(A)}} &&&&&&&& 
CH(A) \ar[ddrrrr]^-{\omega_{Q(A)}C(\chi_{A})} \ar@{-->}[dddl]_{1_{CH(A)}}  \ar@{-->}[ddddddddd]^(0.6){\eta_{H(A)}} &&&& \\
&&&&&&&&&& \\
&&&& H(A) \ar[rrrrrrrr]^(0.6){\theta^{001}_{A}} \ar[ddddddddd]^{\chi_{A}} \ar@{-->}[ddr]^(0.75){1_{H(A)}} &&&&&&&& 
HQ(A) \ar[ddddddddd]^{\chi_{A}} \ar@{-->}[ddlll]_{1_{HQ(A)}} \\
&&& H(A) \ar@{-->}[dddd]_{\chi_{A}} \ar@{-->}[drr]_(0.4){1_{H(A)}} \ar@{-->}[rrrr]^{\theta^{000}_{A}} &&&& 
CH(A) \ar@{-->}[dddd]_{\eta_{H(A)}} \ar@{-->}[drr]^(0.6){\omega_{Q(A)}C(\chi_{A})} &&&&& \\
&&&&& H(A) \ar@{-->}[dddd]_{\chi_{A}} \ar@{-->}[rrrr]_{\theta^{001}_{A}} &&&& HQ(A) \ar@{-->}[dddd]^{\chi_{Q(A)}} &&&\\
&&&&&&&&&&&& \\
&&&&&&&&&&&& \\
&&& GQ(A) \ar@{-->}[drr]_(0.4){G(1_{Q(A)})} \ar@{-->}[rrrr]^-{G(\theta^{010}_{A})} &&&& GKH(A) \ar@{-->}[drr]^{G(\varepsilon_{Q(A)}K(\chi_{A}))} &&&&&&\\
&&&&& GQ(A) \ar@{-->}[rrrr]_-{G(\theta^{011}_{A})} &&&& GQQ(A) &&& \\
GQ(A) \ar@{-->}[rrrrrrrr]^{G(\theta^{010}_{A})} \ar[ddrrrr]_-{G(1_{A})} \ar@{-->}[uurrr]^{G(1_{A})} &&&&&&&& 
GKH(A) \ar@{-->}[ddrrrr]_-{G(\varepsilon_{Q(A)}K(\chi_{A}))} \ar@{-->}[uul]^(0.8){G(1_{KH(A)})} &&&& \\
&&&&&&&&&&& \\
&&&& GQ(A) \ar[rrrrrrrr]_{G(\theta^{011}_{A})} \ar@{-->}[uuur]_-{G(1_{Q(A)})} &&&&&&&& GQQ(A) \ar@{-->}[uuulll]_(0.7){G(1_{QQ(A)})} }
\end{array}
\end{equation}

\item $D(\bm{\theta})_{(C(X),\eta_{X},K(X))}$
\begin{equation}\label{Dthetarho}
\begin{array}{c}
\hspace{-1cm}\xymatrix@!=0.6pc@R1pc@C1.25pc{CC(X) \ar[rrrrrrrr]^-{\theta^{100}_{C(X)}} \ar[ddddddddd]_{\eta_{C(X)}} \ar@{-->}[dddrrr]_(0.6){1_{CC(X)}} \ar[ddrrrr]^-{\omega_{K(X)}C(\eta_{X})} &&&&&&&& 
CCC(X) \ar[ddrrrr]^-{C(\omega_{K(X)}C(\eta_{X}))} \ar@{-->}[dddl]_{\omega_{K(X)}C(\eta_{C(X)})} \ar@{-->}[ddddddddd]^(0.6){\eta_{CC(X)}} &&&& \\
&&&&&&&&&& \\
&&&& HK(X) \ar[rrrrrrrr]^(0.6){\theta^{000}_{K(X)}} \ar[ddddddddd]_{\chi_{K(X)}} \ar@{-->}[ddr]^(0.75){1_{HK(X)}} &&&&&&&& 
CHK(X) \ar[ddddddddd]^{\eta_{HK(X)}} \ar@{-->}[ddlll]_{\omega_{K(X)}C(\chi_{K(X)})} \\
&&& CC(X) \ar@{-->}[dddd]_{\eta_{X}} \ar@{-->}[drr]_(0.4){\omega_{K(X)}C(\eta_{X})} \ar@{-->}[rrrr]^{\theta^{101}_{C(X)}} &&&& 
HKC(X) \ar@{-->}[dddd]_{\chi_{KC(X)}} \ar@{-->}[drr]^(0.5){H(\varepsilon_{K(X)}K(\eta_{X}))} &&&&& \\
&&&&& HK(X) \ar@{-->}[dddd]_{\chi_{K(X)}} \ar@{-->}[rrrr]_{\theta^{001}_{K(X)}} &&&& HQK(X) \ar@{-->}[dddd]^{\chi_{QK(X)}} &&&\\
&&&&&&&&&&&& \\
&&&&&&&&&&&& \\
&&& GKC(X) \ar@{-->}[drr]_(0.4){G(\varepsilon_{K(X)}K(\eta_{X}))} \ar@{-->}[rrrr]^-{G(\theta^{111}_{C(X)})} &&&& GQKC(X) \ar@{-->}[drr]^(0.7){GQ(\varepsilon_{K(X)}K(\eta_{X}))} &&&&&&\\
&&&&& GQK(X) \ar@{-->}[rrrr]_-{G(\theta^{011}_{K(X)})} &&&& GQQK(X) &&& \\
GKC(X) \ar@{-->}[rrrrrrrr]^{G(\theta^{110}_{C(X)})} \ar[ddrrrr]_-{G(\varepsilon_{K(X)}\eta_{X})} \ar@{-->}[uurrr]^{G(1_{KC(X)})} &&&&&&&& 
GKCC(X) \ar@{-->}[ddrrrr]_-{GK(\omega_{K(X)}C(\eta_{X}))}  \ar@{-->}[uul]^(0.8){G(\varepsilon_{K(X)}K(\eta_{X}))} &&&& \\
&&&&&&&&&&& \\
&&&& GQK(X) \ar[rrrrrrrr]_{G(\theta^{010}_{K(X)})} \ar@{-->}[uuur]_-{G(1_{QK(X)})} &&&&&&&& GKHK(X) \ar@{-->}[uuulll]_(0.6){G(\varepsilon_{K(X)}K(\chi_{K(X)}))} }
\end{array}
\end{equation}

\item $D(\bm{\theta})_{(H(A),\chi_{A},Q(A))}$
\begin{equation}\label{Dthetachi}
\begin{array}{c}
\xymatrix@!=0.6pc@R1pc@C1.25pc{CH(A) \ar[rrrrrrrr]^-{\theta^{100}_{H(A)}} \ar[ddddddddd]_{\eta_{H(A)}} \ar@{-->}[dddrrr]_(0.6){1_{CH(A)}} 
\ar[ddrrrr]^-{\omega_{Q(A)}C(\chi_{A})} &&&&&&&& 
CCH(A) \ar[ddrrrr]^-{C(\omega_{Q(A)}C(\chi_{A}))} \ar@{-->}[dddl]_{\omega_{KH(A)}\eta_{H(A)}}  \ar@{-->}[ddddddddd]^(0.6){\eta_{H(A)}} &&&& \\
&&&&&&&&&& \\
&&&& HQ(A) \ar[rrrrrrrr]^(0.6){\theta^{000}_{Q(A)}} \ar[ddddddddd]_{\chi_{Q(A)}} \ar@{-->}[ddr]^(0.75){1_{HQ(A)}} &&&&&&&& 
CHQ(A) \ar[ddddddddd]^{\eta_{H(A)}} \ar@{-->}[ddlll]_{\omega_{Q(A)}C(\chi_{Q(A)})} \\
&&& CH(A) \ar@{-->}[dddd]_{\eta_{H(A)}} \ar@{-->}[drr]_(0.4){\omega_{Q(A)}C(\chi_{A})} \ar@{-->}[rrrr]^-{\theta^{101}_{H(A)}} &&&& 
HKH(A) \ar@{-->}[dddd]_{\chi_{KH(A)}} \ar@{-->}[drr]^(0.5){H(\varepsilon_{Q(A)}K(\chi_{Q(A)}))} &&&&& \\
&&&&& HQ(A) \ar@{-->}[dddd]_{\chi_{A}} \ar@{-->}[rrrr]_-{\theta^{001}_{Q(A)}} &&&& HQQ(A) \ar@{-->}[dddd]^{\chi_{QQ(A)}} &&&\\
&&&&&&&&&&&& \\
&&&&&&&&&&&& \\
&&& GKH(A) \ar@{-->}[drr]_(0.4){G(\varepsilon_{A}K(\chi_{A}))} \ar@{-->}[rrrr]^-{G(\theta^{111}_{H(A)})} &&&& GKH(A) \ar@{-->}[drr]^{G(\varepsilon_{A}K(\chi_{A}))} &&&&&&\\
&&&&& GQQ(A) \ar@{-->}[rrrr]_-{G(\theta^{011}_{Q(A)})} &&&& GQQQ(A) &&& \\
GKH(A) \ar@{-->}[rrrrrrrr]^{G(\theta^{110}_{H(A)})} \ar[ddrrrr]_-{G(\varepsilon_{Q(A)}K(\chi_{A}))} \ar@{-->}[uurrr]^{G(1_{KH(A)})} &&&&&&&& 
GKCH(A) \ar@{-->}[ddrrrr]_-{GK(\omega_{Q(A)}C(\chi_{A}))} \ar@{-->}[uul]^(0.8){G(\varepsilon_{KH(A)}K(\eta_{H(A)}))} &&&& \\
&&&&&&&&&&& \\
&&&& GQQ(A) \ar[rrrrrrrr]_{G(\theta^{010}_{Q(A)})} \ar@{-->}[uuur]_-{G(1_{QQ(A)})} &&&&&&&& GKHQ(A) \ar@{-->}[uuulll]_(0.8){G(\varepsilon_{Q(A)}K(\chi_{Q(A)}))} }
\end{array}
\end{equation}

\end{itemize}

\begin{itemize}

\item The composition of (\ref{xihyperactionX}) and (\ref{DIDthetaX}) 
\begin{equation}\label{DIDthetaXxihyperactionX}
\begin{array}{c}
\xymatrix@!=0.6pc@R1pc@C1.25pc{C(X) \ar[rrrrrrrr]^-{C(\theta^{100}_{X})\theta^{100}_{X}} \ar[ddddddddd]_{\eta_{X}} \ar@{-->}[dddrrr]_(0.6){1_{C(X)}} \ar[ddrrrr]^-{1_{C(X)}} &&&&&&&& 
CCC(X) \ar[ddrrrr]^-{C(\omega_{K(X)}C(\eta_{X}))} \ar@{-->}[dddl]_{\omega_{K(X)}C(\eta_{C(X)})} \ar@{-->}[ddddddddd]^(0.6){\eta_{CC(X)}} &&&& \\
&&&&&&&&&& \\
&&&& C(X) \ar[rrrrrrrr]^(0.6){C(\theta^{101}_{X})\theta^{100}_{X}} \ar[ddddddddd]_{\eta_{X}} \ar@{-->}[ddr]^(0.75){1_{C(X)}} &&&&&&&& 
CHK(X) \ar[ddddddddd]^{\eta_{HK(X)}} \ar@{-->}[ddlll]_{\omega_{K(X)}C(\chi_{K(X)})} \\
&&& C(X) \ar@{-->}[dddd]_{\eta_{X}} \ar@{-->}[drr]_(0.4){\omega_{K(X)}C(\eta_{X})} \ar@{-->}[rrrr]^{H(\theta^{110}_{X})\theta^{101}_{X}} &&&& 
HKC(X) \ar@{-->}[dddd]_{\chi_{KC(X)}} \ar@{-->}[drr]^(0.5){H(\varepsilon_{K(X)}K(\eta_{X}))} &&&&& \\
&&&&& C(X) \ar@{-->}[dddd]_{\eta_{X}} \ar@{-->}[rrrr]_{H(\theta^{111}_{X})\theta^{101}_{X}} &&&& HQK(X) \ar@{-->}[dddd]^{\chi_{QK(X)}} &&&\\
&&&&&&&&&&&& \\
&&&&&&&&&&&& \\
&&& GKC(X) \ar@{-->}[drr]_(0.4){G(\varepsilon_{K(X)}K(\eta_{X}))} \ar@{-->}[rrrr]^-{GQ(\theta^{110}_{X})G(\theta^{111}_{X})} &&&& GQKC(X) \ar@{-->}[drr]^(0.7){GQ(\varepsilon_{K(X)}K(\eta_{X}))} &&&&&&\\
&&&&& GQK(X) \ar@{-->}[rrrr]_-{GQ(\theta^{111}_{X})G(\theta^{111}_{X})} &&&& GQQK(X) &&& \\
GK(X) \ar@{-->}[rrrrrrrr]^{GK(\theta^{100}_{X})G(\theta^{110}_{X})} \ar[ddrrrr]_-{G(1_{K(X)})} \ar@{-->}[uurrr]^{G(1_{K(X)})} &&&&&&&& 
GKCC(X) \ar@{-->}[ddrrrr]_-{GK(\omega_{K(X)}C(\eta_{X}))}  \ar@{-->}[uul]^(0.8){G(\varepsilon_{K(X)}K(\eta_{X}))} &&&& \\
&&&&&&&&&&& \\
&&&& GK(X) \ar[rrrrrrrr]_{GK(\theta^{101}_{X})G(\theta^{110}_{X})} \ar@{-->}[uuur]_-{G(1_{K(X)})} &&&&&&&& GKHK(X) \ar@{-->}[uuulll]_(0.6){G(\varepsilon_{K(X)}K(\chi_{K(X)}))} } 
\end{array}
\end{equation}

\item The composition of (\ref{DnIDX}) and (\ref{Dthetarho})
\begin{equation}\label{DthetarhoDnIDX}
\begin{array}{c}
\xymatrix@!=0.6pc@R1pc@C1.25pc{C(X) \ar[rrrrrrrr]^-{\theta^{100}_{C(X)}\theta^{100}_{X}} \ar[ddddddddd]_{\eta_{X}} \ar@{-->}[dddrrr]_(0.6){1_{C(X)}} \ar[ddrrrr]^-{1_{C(X)}} &&&&&&&& 
CCC(X) \ar[ddrrrr]^-{C(\omega_{K(X)}C(\eta_{X}))} \ar@{-->}[dddl]_{\omega_{K(X)}C(\eta_{C(X)})} \ar@{-->}[ddddddddd]^(0.6){\eta_{CC(X)}} &&&& \\
&&&&&&&&&& \\
&&&& C(X) \ar[rrrrrrrr]^(0.6){\theta^{000}_{K(X)}\theta^{101}_{X}} \ar[ddddddddd]_{\eta_{X}} \ar@{-->}[ddr]^(0.75){1_{C(X)}} &&&&&&&& 
CHK(X) \ar[ddddddddd]^{\eta_{HK(X)}} \ar@{-->}[ddlll]_{\omega_{K(X)}C(\chi_{K(X)})} \\
&&& C(X) \ar@{-->}[dddd]_{\eta_{X}} \ar@{-->}[drr]_(0.4){\omega_{K(X)}C(\eta_{X})} \ar@{-->}[rrrr]^{\theta^{101}_{C(X)}\theta^{100}_{X}} &&&& 
HKC(X) \ar@{-->}[dddd]_{\chi_{KC(X)}} \ar@{-->}[drr]^(0.5){H(\varepsilon_{K(X)}K(\eta_{X}))} &&&&& \\
&&&&& C(X) \ar@{-->}[dddd]_{\eta_{X}} \ar@{-->}[rrrr]_{\theta^{001}_{K(X)}\theta^{101}_{X}} &&&& HQK(X) \ar@{-->}[dddd]^{\chi_{QK(X)}} &&&\\
&&&&&&&&&&&& \\
&&&&&&&&&&&& \\
&&& GKC(X) \ar@{-->}[drr]_(0.4){G(\varepsilon_{K(X)}K(\eta_{X}))} \ar@{-->}[rrrr]^-{G(\theta^{111}_{C(X)}\theta^{110}_{X})} &&&& GQKC(X) \ar@{-->}[drr]^(0.7){GQ(\varepsilon_{K(X)}K(\eta_{X}))} &&&&&&\\
&&&&& GQK(X) \ar@{-->}[rrrr]_-{G(\theta^{011}_{K(X)}\theta^{111}_{X})} &&&& GQQK(X) &&& \\
GK(X) \ar@{-->}[rrrrrrrr]^{G(\theta^{110}_{C(X)}\theta^{110}_{X})} \ar[ddrrrr]_-{G(1_{K(X)})} \ar@{-->}[uurrr]^{G(1_{K(X)})} &&&&&&&& 
GKCC(X) \ar@{-->}[ddrrrr]_-{GK(\omega_{K(X)}C(\eta_{X}))}  \ar@{-->}[uul]^(0.8){G(\varepsilon_{K(X)}K(\eta_{X}))} &&&& \\
&&&&&&&&&&& \\
&&&& GK(X) \ar[rrrrrrrr]_{G(\theta^{010}_{K(X)})\theta^{111}_{X})} \ar@{-->}[uuur]_-{G(1_{K(X)})} &&&&&&&& GKHK(X) \ar@{-->}[uuulll]_(0.6){G(\varepsilon_{K(X)}K(\chi_{K(X)}))} } 
\end{array}
\end{equation}

\end{itemize}
The identity (\ref{DIDcoassociativity1}) says that (\ref{DthetarhoDnIDX}) and  (\ref{DIDthetaXxihyperactionX}) are equal. This gives us the identities:
\begin{equation}\label{thetaX000}
C(\theta^{100}_{X}) \theta^{100}_{X} =  \theta^{100}_{C(X)} \theta^{100}_{X},
\end{equation}
\begin{equation}\label{thetaX001}
C(\theta^{101}_{X}) \theta^{100}_{X} =  \theta^{000}_{K(X)} \theta^{101}_{X},
\end{equation}
\begin{equation}\label{thetaX010}
H(\theta^{110}_{X}) \theta^{101}_{X} =  \theta^{101}_{C(X)} \theta^{100}_{X},
\end{equation}
\begin{equation}\label{thetaX011}
H(\theta^{111}_{X}) \theta^{101}_{X} =  \theta^{001}_{K(X)} \theta^{101}_{X},
\end{equation}
\begin{equation}\label{thetaX100}
Q(\theta^{110}_{X}) \theta^{111}_{X} =  \theta^{111}_{C(X)} \theta^{110}_{X},
\end{equation}
\begin{equation}\label{thetaX101}
Q(\theta^{111}_{X}) \theta^{111}_{X} =  \theta^{011}_{K(X)} \theta^{111}_{X},
\end{equation}
\begin{equation}\label{thetaX110}
K(\theta^{100}_{X})\theta^{110}_{X} =  \theta^{110}_{C(X)} \theta^{110}_{X},
\end{equation}
\begin{equation}\label{thetaX111}
K(\theta^{101}_{X}) \theta^{110}_{X} =  \theta^{010}_{K(X)} \theta^{111}_{X}
\end{equation}

\begin{itemize}

\item The composition of (\ref{xihyperactionA}) and (\ref{DIDthetaA})
\begin{equation}\label{DthetachiDnIDA}
\begin{array}{c}
\xymatrix@!=0.6pc@R1pc@C1.25pc{H(A) \ar[rrrrrrrr]^-{C(\theta^{000}_{A})\theta^{000}_{A}} \ar[ddddddddd]_{\chi_{A}} \ar@{-->}[dddrrr]_(0.6){1_{H(A)}} 
\ar[ddrrrr]^-{1_{H(A)}} &&&&&&&& 
CCH(A) \ar[ddrrrr]^-{C(\omega_{Q(A)}C(\chi_{A}))} \ar@{-->}[dddl]_{\omega_{KH(A)}\eta_{H(A)}}  \ar@{-->}[ddddddddd]^(0.6){\eta_{H(A)}} &&&& \\
&&&&&&&&&& \\
&&&& H(A) \ar[rrrrrrrr]^(0.6){C(\theta^{001}_{A})\theta^{000}_{A}} \ar[ddddddddd]_{\chi_{A}} \ar@{-->}[ddr]^(0.75){1_{H(A)}} &&&&&&&& 
CHQ(A) \ar[ddddddddd]^{\eta_{H(A)}} \ar@{-->}[ddlll]_{\omega_{Q(A)}C(\chi_{Q(A)})} \\
&&& H(A) \ar@{-->}[dddd]_{\chi_{A}} \ar@{-->}[drr]_(0.4){1_{H(A)}}   \ar@{-->}[rrrr]^-{H(\theta^{010}_{A})\theta^{001}_{A}} &&&& 
HKH(A) \ar@{-->}[dddd]_{\chi_{KH(A)}} \ar@{-->}[drr]^(0.5){H(\varepsilon_{Q(A)}K(\chi_{Q(A)}))} &&&&& \\
&&&&& H(A) \ar@{-->}[dddd]_{\chi_{A}} \ar@{-->}[rrrr]_-{H(\theta^{011}_{A})\theta^{001}_{A}} &&&& HQQ(A) \ar@{-->}[dddd]^{\chi_{QQ(A)}} &&&\\
&&&&&&&&&&&& \\
&&&&&&&&&&&& \\
&&& GQ(A) \ar@{-->}[drr]_(0.4){G(1_{Q(A)})} \ar@{-->}[rrrr]^-{G(Q(\theta^{010}_{A})\theta^{011}_{A})} &&&& GKH(A) \ar@{-->}[drr]^{G(\varepsilon_{A}K(\chi_{A}))} &&&&&&\\
&&&&& GQ(A) \ar@{-->}[rrrr]_-{G(Q(\theta^{011}_{A})\theta^{011}_{A})} &&&& GQQQ(A) &&& \\
GQ(A) \ar@{-->}[rrrrrrrr]^{G(K(\theta^{000}_{A})\theta^{110}_{A})} \ar[ddrrrr]_-{G(1_{Q(A)})} \ar@{-->}[uurrr]^{G(1_{Q(A)})} &&&&&&&& 
GKCH(A) \ar@{-->}[ddrrrr]_-{GK(\omega_{Q(A)}C(\chi_{A}))} \ar@{-->}[uul]^(0.8){G(\varepsilon_{KH(A)}K(\eta_{H(A)}))} &&&& \\
&&&&&&&&&&& \\
&&&& GQ(A) \ar[rrrrrrrr]_{G(K(\theta^{001}_{A})\theta^{010}_{A})} \ar@{-->}[uuur]_-{G(1_{Q(A)})} &&&&&&&& GKHQ(A) \ar@{-->}[uuulll]_(0.8){G(\varepsilon_{Q(A)}K(\chi_{Q(A)}))} }
\end{array}
\end{equation}

\item The composition of (\ref{DnIDA}) and (\ref{Dthetachi})
\begin{equation}\label{DIDthetaAxihyperactionA}
\begin{array}{c}
\xymatrix@!=0.6pc@R1pc@C1.25pc{H(A) \ar[rrrrrrrr]^-{\theta^{100}_{H(A)}\theta^{000}_{A}} \ar[ddddddddd]_{\chi_{A}} \ar@{-->}[dddrrr]_(0.6){1_{H(A)}} 
\ar[ddrrrr]^-{1_{H(A)}} &&&&&&&& 
CCH(A) \ar[ddrrrr]^-{C(\omega_{Q(A)}C(\chi_{A}))} \ar@{-->}[dddl]_{\omega_{KH(A)}\eta_{H(A)}}  \ar@{-->}[ddddddddd]^(0.6){\eta_{H(A)}} &&&& \\
&&&&&&&&&& \\
&&&& H(A) \ar[rrrrrrrr]^(0.6){\theta^{000}_{Q(A)}\theta^{001}_{A}} \ar[ddddddddd]_{\chi_{A}} \ar@{-->}[ddr]^(0.75){1_{H(A)}} &&&&&&&& 
CHQ(A) \ar[ddddddddd]^{\eta_{H(A)}} \ar@{-->}[ddlll]_{\omega_{Q(A)}C(\chi_{Q(A)})} \\
&&& H(A) \ar@{-->}[dddd]_{\chi_{A}} \ar@{-->}[drr]_(0.4){1_{H(A)}}   \ar@{-->}[rrrr]^-{\theta^{101}_{H(A)}\theta^{000}_{A}} &&&& 
HKH(A) \ar@{-->}[dddd]_{\chi_{KH(A)}} \ar@{-->}[drr]^(0.5){H(\varepsilon_{Q(A)}K(\chi_{Q(A)}))} &&&&& \\
&&&&& H(A) \ar@{-->}[dddd]_{\chi_{A}} \ar@{-->}[rrrr]_-{\theta^{001}_{Q(A)}\theta^{001}_{A}} &&&& HQQ(A) \ar@{-->}[dddd]^{\chi_{QQ(A)}} &&&\\
&&&&&&&&&&&& \\
&&&&&&&&&&&& \\
&&& GQ(A) \ar@{-->}[drr]_(0.4){G(1_{Q(A)})} \ar@{-->}[rrrr]^-{G(\theta^{111}_{H(A)}\theta^{010}_{A})} &&&& GKH(A) \ar@{-->}[drr]^{G(\varepsilon_{A}K(\chi_{A}))} &&&&&&\\
&&&&& GQ(A) \ar@{-->}[rrrr]_-{G(\theta^{011}_{Q(A)}\theta^{011}_{A})} &&&& GQQQ(A) &&& \\
GQ(A) \ar@{-->}[rrrrrrrr]^{G(\theta^{110}_{H(A)}\theta^{010}_{A})} \ar[ddrrrr]_-{G(1_{Q(A)})} \ar@{-->}[uurrr]^{G(1_{Q(A)})} &&&&&&&& 
GKCH(A) \ar@{-->}[ddrrrr]_-{GK(\omega_{Q(A)}C(\chi_{A}))} \ar@{-->}[uul]^(0.8){G(\varepsilon_{KH(A)}K(\eta_{H(A)}))} &&&& \\
&&&&&&&&&&& \\
&&&& GQ(A) \ar[rrrrrrrr]_{G(\theta^{010}_{Q(A)}\theta^{011}_{A})} \ar@{-->}[uuur]_-{G(1_{Q(A)})} &&&&&&&& GKHQ(A) \ar@{-->}[uuulll]_(0.8){G(\varepsilon_{Q(A)}K(\chi_{Q(A)}))} }
\end{array}
\end{equation}
\end{itemize}
The identity (\ref{DIDcoassociativity0}) says that (\ref{DthetachiDnIDA}) and (\ref{DIDthetaAxihyperactionA}) are equal. This gives us the identities:
\begin{equation}\label{thetaA000}
C(\theta^{000}_{A}) \theta^{000}_{A} =  \theta^{100}_{H(A)} \theta^{000}_{A},
\end{equation}
\begin{equation}\label{thetaA001}
C(\theta^{001}_{A}) \theta^{000}_{A} =  \theta^{000}_{Q(A)} \theta^{001}_{A},
\end{equation}
\begin{equation}\label{thetaA010}
H(\theta^{010}_{A}) \theta^{001}_{A} =  \theta^{101}_{H(A)} \theta^{000}_{A},
\end{equation}
\begin{equation}\label{thetaA011}
H(\theta^{011}_{A}) \theta^{001}_{A} =  \theta^{001}_{Q(A)} \theta^{001}_{A},
\end{equation}
\begin{equation}\label{thetaA100}
Q(\theta^{010}_{A}) \theta^{011}_{A} =  \theta^{111}_{H(A)} \theta^{010}_{A},
\end{equation}
\begin{equation}\label{thetaA101}
Q(\theta^{011}_{A}) \theta^{011}_{A} =  \theta^{011}_{Q(A)} \theta^{011}_{A},
\end{equation}
\begin{equation}\label{thetaA110}
K(\theta^{000}_{A})\theta^{010}_{A} =  \theta^{110}_{H(A)} \theta^{010}_{A},
\end{equation}
\begin{equation}\label{thetaA111}
K(\theta^{001}_{A}) \theta^{010}_{A} =  \theta^{010}_{Q(A)} \theta^{011}_{A}.
\end{equation}
In order to summarize what all these equations mean let us recall definitions of the Eilenberg-Moore and Kleisli categories of coalgebras over a comonad:

\begin{definition}\label{EM}
Let $(C,\tau,\zeta)$ be a comonad on a category $\X$.  We say that a pair $(X,x)$ where $x \maps X \to C(X)$ is a $C$-coalgebra if it satisfies the following two axioms:
\[\xymatrix@!=3pc{X \ar[r]^-{x} \ar[d]_-{x} & C(X) \ar[d]^{\tau_{X}} \\
C(X) \ar[r]_-{C(x)} & CC(X)  }
\hspace{1cm}
\xymatrix@!=3pc{X \ar[r]^-{1_{X}} \ar[d]_-{x} & X \ar[d]^{1_{X}} \\
C(X) \ar[r]_-{\zeta_{X}} & X }\]
A morphism of colagebras $f \maps (X,x) \to (X',x')$ is a morphism $x \maps X \to X'$ satisfying
\[\xymatrix@!=3pc{X \ar[r]^-{f} \ar[d]_-{x} & X' \ar[d]^{x'} \\
C(X) \ar[r]_-{C(f)} & C(X') }\]
$\C$-oalgebras and their morphisms form a category $\X^{C}$ which we call the Eilenberg-Moore category of the comonad $\C$.
\end{definition}
The following is well known consequence of Definition \ref{EM}:
\begin{theorem}\label{EMadjunction}
Any comonad $(C,\tau,\zeta)$ on a category $\X$ induces an adjunction
\[\xymatrix@!=3pc{\X^{C} \ar@<1ex>[r]^-{U^{C}} & \X \ar@<1ex>[l]^-{F^{C}}_{\perp} }\]
with the forgetful functor $U^{C}$ which sends any $C$-coalgebra  $(X,x)$ to $X$ left adjoint to the free functor $F^{C}$ which sends any object $X$ in $\X$ to the free coalgebra $(C(X),\tau_{X})$.
\begin{proof}
For any colagebra $(X,x)$ the component of the unit $\eta \maps I_{\X^{C}} \Rightarrow F^{C}U^{C}$ of the adjunction is given by the first diagram of Definition \ref{EM} and $\zeta_{X}$ is the component of the counit $\zeta \maps U^{C}F^{C} \Rightarrow I_{\X}$ indexed by the object $X$ in $\X$. The verification of triangle identities is straightforward.
\end{proof}
\end{theorem}

Now we state without a proof the existence of the Kleisli category of a comonad:
\begin{theorem}\label{Kleisli}
Let $(C,\tau,\zeta)$ be a comonad on a category $\X$.  Then the Kleisli category $\X_{C}$ has the same objects as $\X$ and a hom-set $\X_{C}(X,Y)$ of Kleisli morphisms from an object $X$ to $Y$ is by definition $\X(C(X),Y)$. The composition of two  Kleisli morphisms $f \maps C(X) \to Y$ and $g \maps C(Y) \to Z$ is given by 
\[\xymatrix@!=1pc{C(X) \ar[rr]^-{\theta^{100}_{X}} && CC(X) \ar[rr]^-{C(f)} && C(Y)  \ar[rr]^-{g} && Z}\]
\end{theorem}

\begin{theorem}\label{Kleisliadjunction}
Any comonad $(C,\tau,\zeta)$ on a category $\X$ induces an adjunction
\[\xymatrix@!=3pc{\X_{C} \ar@<1ex>[r]^-{U_{C}} & \X \ar@<1ex>[l]^-{F_{C}}_{\perp} }\]
where $F_{C}$ is identity on objects functor which sends any morphism $f \maps X \to Y$ in $\X$ to a Kleisli morphsim
\begin{equation}\label{Kleislimorphism}
\begin{array}{c}
\xymatrix@!=3pc{C(X) \ar[r]^-{\zeta^{1}_{X}} & X  \ar[r]^-{f} & Y}
\end{array}
\end{equation}
and $U_{C}$ acts as $C$ on objects and takes any Kleisli morphism $\widetilde{f} \maps C(X) \to Y$ from $X$ to $Y$ in $\X_{C}$ to 
\begin{equation}\label{Kleislimorphism2}
\begin{array}{c}
\xymatrix@!=3pc{C(X) \ar[r]^-{\theta^{100}_{X}} & CC(X) \ar[r]^-{C(\widetilde{f})} & C(Y).}
\end{array}
\end{equation}
\begin{proof}
The component of the counit $\zeta \maps U_{C}F_{C} \Rightarrow I_{\X}$ indexed by an object $X$ in $\X$ is a morphism $\zeta^{1}_{X}
\maps C(X) \to X$ and the component of the unit $\eta \maps I_{\X_{C}} \Rightarrow F_{C}U_{C}$ of the adjunction is a Kleisli morphism $1_{C(X)}
\maps C(X) \to C(X)$ from $X$ to $C(X)$ in $\X_{C}$. The verification of triangle identities is straightforward.
\end{proof}
\end{theorem}
The equations which any colax $\D$-coalgebra needs to satisfy correspond to the following commutative diagrams
\[\hspace{-0.75cm}\xymatrix@!=3pc{C(X) \ar[d]_{\theta^{100}_{X}} \ar[r]^-{1_{C(X)}} & C(X) \ar[d]^-{1_{C(X)}} \\
CK(X) \ar[r]_-{C(\zeta^{1}_{X})} & C(X)}
\xymatrix@!=3pc{K(X) \ar[d]_{\theta^{110}_{X}} \ar[r]^-{1_{K(X)}} & K(X) \ar[d]^-{1_{K(X)}} \\
KC(X) \ar[r]_-{K(\zeta^{1}_{X})} & K(X)}
\xymatrix@!=3pc{H(A) \ar[d]_{\theta^{001}_{A}} \ar[r]^-{1_{H(A)}} & H(A) \ar[d]^-{1_{H(A)}} \\
HQ(A) \ar[r]_-{H(\zeta^{0}_{A})} & H(A)}
\xymatrix@!=3pc{Q(A) \ar[d]_{\theta^{011}_{A}} \ar[r]^-{1_{Q(A)}} & Q(A) \ar[d]^-{1_{Q(A)}} \\
QQ(A) \ar[r]_-{Q(\zeta^{0}_{A})} & Q(A)}\]

\[\hspace{-0.75cm}\xymatrix@!=3pc{C(X) \ar[d]_{\theta^{100}_{X}} \ar[r]^-{1_{C(X)}} & C(X) \ar[d]^-{1_{C(X)}} \\
CC(X) \ar[r]_-{\zeta^{1}_{C(X)}} & C(X)}
\xymatrix@!=3pc{K(X) \ar[d]_{\theta^{111}_{X}} \ar[r]^-{1_{K(X)}} & K(X) \ar[d]^-{1_{K(X)}} \\
QK(X) \ar[r]_-{\zeta^{0}_{K(X)}} & K(X)}
\xymatrix@!=3pc{H(A) \ar[d]_{\theta^{000}_{A}} \ar[r]^-{1_{H(A)}} & H(A) \ar[d]^-{1_{H(A)}} \\
CH(A) \ar[r]_-{\zeta^{1}_{H(A)}} & H(A)}
\xymatrix@!=3pc{Q(A) \ar[d]_{\theta^{011}_{A}} \ar[r]^-{1_{Q(A)}} & Q(A) \ar[d]^-{1_{Q(A)}} \\
QQ(A) \ar[r]_-{\zeta^{0}_{Q(A)}} & Q(A)}\]

\[\hspace{-0.75cm}\xymatrix@!=3pc{C(X) \ar[d]_{\theta^{100}_{X}} \ar[r]^-{\theta^{100}_{X}} & CC(X) \ar[d]^-{\theta^{100}_{C(X)}} \\
CC(X) \ar[r]_-{C(\theta^{100}_{X})} & CCC(X)}
\xymatrix@!=3pc{C(X) \ar[d]_{\theta^{100}_{X}} \ar[r]^-{\theta^{101}_{X}} & HK(X) \ar[d]^-{\theta^{000}_{K(X)}} \\
CC(X) \ar[r]_-{C(\theta^{101}_{X})} & CHK(X)}
\xymatrix@!=3pc{C(X) \ar[d]_{\theta^{101}_{X}} \ar[r]^-{\theta^{100}_{X}} & CC(X) \ar[d]^-{\theta^{101}_{C(X)}} \\
HK(X) \ar[r]_-{H(\theta^{110}_{X})} & HKC(X)}
\xymatrix@!=3pc{C(X) \ar[d]_{\theta^{101}_{X}} \ar[r]^-{\theta^{101}_{X}} & HK(X) \ar[d]^-{\theta^{001}_{K(X)}} \\
HK(X) \ar[r]_-{H(\theta^{111}_{X})} & HQK(X)}\]

\[\hspace{-0.75cm}\xymatrix@!=3pc{K(X) \ar[d]_{\theta^{111}_{X}} \ar[r]^-{\theta^{110}_{X}} & KC(X) \ar[d]^-{\theta^{111}_{C(X)}} \\
QK(X) \ar[r]_-{Q(\theta^{110}_{X})} & QKC(X)}
\xymatrix@!=3pc{K(X) \ar[d]_{\theta^{111}_{X}} \ar[r]^-{\theta^{111}_{X}} & QK(X) \ar[d]^-{\theta^{011}_{K(X)}} \\
QK(X) \ar[r]_-{Q(\theta^{111}_{X})} & QQK(X)}
\xymatrix@!=3pc{K(X) \ar[d]_{\theta^{110}_{X}} \ar[r]^-{\theta^{110}_{X}} & KC(X) \ar[d]^-{\theta^{110}_{C(X)}} \\
KC(X) \ar[r]_-{K(\theta^{100}_{X})} & KCC(X)}
\xymatrix@!=3pc{K(X) \ar[d]_{\theta^{110}_{X}} \ar[r]^-{\theta^{111}_{X}} & QK(X) \ar[d]^-{\theta^{010}_{K(X)}} \\
KC(X) \ar[r]_-{K(\theta^{101}_{X})} & KHK(X)}\]

\[\hspace{-0.75cm}\xymatrix@!=3pc{H(A) \ar[d]_{\theta^{000}_{A}} \ar[r]^-{\theta^{000}_{A}} & CH(A) \ar[d]^-{\theta^{100}_{H(A)}} \\
CH(A) \ar[r]_-{C(\theta^{000}_{A})} & CCH(A)}
\xymatrix@!=3pc{H(A) \ar[d]_{\theta^{000}_{A}} \ar[r]^-{\theta^{001}_{A}} & HQ(A) \ar[d]^-{\theta^{000}_{Q(A)}} \\
CH(A) \ar[r]_-{C(\theta^{001}_{A})} & CHQ(A)}
\xymatrix@!=3pc{H(A) \ar[d]_{\theta^{001}_{A}} \ar[r]^-{\theta^{000}_{A}} & CH(A) \ar[d]^-{\theta^{101}_{H(A)}} \\
HQ(A) \ar[r]_-{H(\theta^{010}_{A})} & HKH(A)}
\xymatrix@!=3pc{H(A) \ar[d]_{\theta^{001}_{A}} \ar[r]^-{\theta^{001}_{A}} & HQ(A) \ar[d]^-{\theta^{001}_{Q(A)}} \\
HQ(A) \ar[r]_-{H(\theta^{011}_{A})} & HQQ(A)}\]

\[\hspace{-0.75cm}\xymatrix@!=3pc{Q(A) \ar[d]_{\theta^{011}_{A}} \ar[r]^-{\theta^{010}_{A}} & KH(A) \ar[d]^-{\theta^{111}_{H(A)}} \\
QQ(A) \ar[r]_-{Q(\theta^{010}_{A})} & QKH(A)}
\xymatrix@!=3pc{Q(A) \ar[d]_{\theta^{011}_{A}} \ar[r]^-{\theta^{011}_{A}} & QQ(A) \ar[d]^-{\theta^{011}_{Q(A)}} \\
QQ(A) \ar[r]_-{Q(\theta^{011}_{A})} & QQQ(A)}
\xymatrix@!=3pc{Q(A) \ar[d]_{\theta^{010}_{A}} \ar[r]^-{\theta^{010}_{A}} & KH(A) \ar[d]^-{\theta^{110}_{H(A)}} \\
KH(A) \ar[r]_-{K(\theta^{000}_{A})} & KCH(A)}
\xymatrix@!=3pc{Q(A) \ar[d]_{\theta^{010}_{A}} \ar[r]^-{\theta^{011}_{A}} & QQ(A) \ar[d]^-{\theta^{010}_{Q(A)}} \\
KH(A) \ar[r]_-{K(\theta^{001}_{A})} & KHQ(A)}\]

Now we give an interpretation of the diagrams defining a colax $\D$-coalgebra:
\begin{itemize}
\item $(C,\theta^{100},\zeta^{1})$ is a comonad on $\X$ by (\ref{leftuniteqX1}),  (\ref{rightuniteqX1}) and (\ref{thetaX000}) 
\[\xymatrix@!=2pc{C(X) \ar[d]_{\theta^{100}_{X}} \ar[r]^-{1_{C(X)}} & C(X) \ar[d]^-{1_{C(X)}} \\
CC(X) \ar[r]_-{C(\zeta^{1}_{X})} & C(X)}
\hspace{0.5cm}
\xymatrix@!=2pc{C(X) \ar[d]_{\theta^{100}_{X}} \ar[r]^-{1_{C(X)}} & C(X) \ar[d]^-{1_{C(X)}} \\
CC(X) \ar[r]_-{\zeta^{1}_{C(X)}} & C(X)}
\hspace{0.5cm}
\xymatrix@!=2pc{C(X) \ar[d]_{\theta^{100}_{X}} \ar[r]^-{\theta^{100}_{X}} & CC(X) \ar[d]^-{\theta^{100}_{C(X)}} \\
CC(X) \ar[r]_-{C(\theta^{100}_{X})} & CCC(X)}\]

\item $(Q,\theta^{011},\zeta^{0})$ is a comonad on $\A$ by (\ref{leftuniteqA2}),  (\ref{rightuniteqA2}) and (\ref{thetaA101}) 
\[\xymatrix@!=2.5pc{Q(A) \ar[d]_{\theta^{011}_{A}} \ar[r]^-{1_{Q(A)}} & Q(A) \ar[d]^-{1_{Q(A)}} \\
QQ(A) \ar[r]_-{Q(\zeta_{A})} & Q(A)}
\hspace{0.5cm}
\xymatrix@!=2.5pc{Q(A) \ar[d]_{\theta^{011}_{A}} \ar[r]^-{1_{Q(A)}} & Q(A) \ar[d]^-{1_{Q(A)}} \\
QQ(A) \ar[r]_-{\zeta^{0}_{Q(A)}} & Q(A)}
\hspace{0.5cm}
\xymatrix@!=2.5pc{Q(A) \ar[d]_{\theta^{011}_{A}} \ar[r]^-{\theta^{011}_{A}} & QQ(A) \ar[d]^-{\theta^{011}_{Q(A)}} \\
QQ(A) \ar[r]_-{Q(\theta^{011}_{A})} & QQQ(A)}\]

\item $(K(X),\theta^{111}_{X})$ is $Q$-coalgebra for every object $X$ in $\X$ by (\ref{rightuniteqX2}) and (\ref{thetaX101}) 
\[\hspace{0.5cm}
\xymatrix@!=2.5pc{K(X) \ar[d]_{\theta^{111}_{X}} \ar[r]^-{1_{K(X)}} & K(X) \ar[d]^-{1_{K(X)}} \\
QK(X) \ar[r]_-{\zeta^{0}_{K(X)}} & K(X)}
\hspace{0.5cm}
\xymatrix@!=2.5pc{K(X) \ar[d]_{\theta^{111}_{X}} \ar[r]^-{\theta^{111}_{X}} & QK(X) \ar[d]^-{\theta^{011}_{K(X)}} \\
QK(X) \ar[r]_-{Q(\theta^{111}_{X})} & QQK(X)}\]

\item $(H(A),\theta^{000}_{A})$ is a $C$-coalgebra for every object $A$ in $\A$ by (\ref{rightuniteqA1}) and (\ref{thetaA000})
\[\xymatrix@!=2.5pc{H(A) \ar[d]_{\theta^{000}_{A}} \ar[r]^-{1_{H(A)}} & H(A) \ar[d]^-{1_{H(A)}} \\
CH(A) \ar[r]_-{\zeta^{1}_{H(A)}} & H(A)}
\hspace{0.5cm}
\xymatrix@!=2.5pc{H(A) \ar[d]_{\theta^{000}_{A}} \ar[r]^-{\theta^{000}_{A}} & CH(A) \ar[d]^-{\theta^{100}_{H(A)}} \\
CH(A) \ar[r]_-{C(\theta^{000}_{A})} & CCH(A)}\]

\item $\theta^{101}_{X} \maps (C(X),\theta^{100}_{X}) \to (HK(X),\theta^{000}_{K(X)})$ is a morphism of $C$-coalgebras by (\ref{thetaX001}) 
\[\xymatrix@!=2.5pc{C(X) \ar[d]_{\theta^{100}_{X}} \ar[r]^-{\theta^{101}_{X}} & HK(X) \ar[d]^-{\theta^{000}_{K(X)}} \\
CC(X) \ar[r]_-{C(\theta^{101}_{X})} & CHK(X)}\]
(from the free $C$-coalgebra on $C(X)$) for every object $X$ in $\X$,

\item $\theta^{001}_{X} \maps (H(A),\theta^{000}_{A}) \to (HQ(A),\theta^{000}_{Q(A)})$ is a morphism of $C$-coalgebras by (\ref{thetaA001}) 
\[\xymatrix@!=2.5pc{H(A) \ar[d]_{\theta^{000}_{A}} \ar[r]^-{\theta^{001}_{A}} & HQ(A) \ar[d]^-{\theta^{000}_{Q(A)}} \\
CH(A) \ar[r]_-{C(\theta^{001}_{A})} & CHQ(A)}\]

\item $\theta^{110}_{X} \maps (K(X),\theta^{111}_{X}) \to (KC(X),\theta^{111}_{C(X)})$ is a morphism of $Q$-coalgebras by (\ref{thetaX100})
\[\xymatrix@!=2.5pc{K(X) \ar[d]_{\theta^{111}_{X}} \ar[r]^-{\theta^{110}_{X}} & KC(X) \ar[d]^-{\theta^{111}_{C(X)}} \\
QK(X) \ar[r]_-{Q(\theta^{110}_{X})} & QKC(X)}\]

\item $\theta^{010}_{A} \maps (Q(A),\theta^{011}_{A}) \to (KH(A),\theta^{111}_{H(A)})$  is a morphism of $Q$-coalgebras by  (\ref{thetaA100}) 
\[\xymatrix@!=2.5pc{Q(A) \ar[d]_{\theta^{011}_{A}} \ar[r]^-{\theta^{010}_{A}} & KH(A) \ar[d]^-{\theta^{111}_{H(A)}} \\
QQ(A) \ar[r]_-{Q(\theta^{010}_{A})} & QKH(A)}\]
(from the free $Q$-coalgebra on $Q(A)$) for every object $A$ in $\A$.

\end{itemize}

Now we give a purely combinatorial interpretation of the rest of diagrams in the comma category $(\X,G)$.
%%%%%%%%%%%%%from here
\begin{proposition}\label{Xprop}
For any colax $\D$-coalgebra $(G,\bm{F}_{G},\bm{\zeta},\bm{\theta})$ the following diagram
\[\hspace{-1cm}\xymatrix@!=1.15pc@R1pc@C1pc{C(X) \ar[dddd]_{\varrho_{X}} \ar[rrrr]^-{\theta^{100}_{X}} \ar[dr]^(0.7){\theta^{100}_{X}} &&&& CC(X) \ar@{-->}[dddd]^(0.6){\varrho_{C(X)}} \ar[dr]^(0.7){\theta^{100}_{C(X)}} \ar[rrrr]^-{C(\varrho_{X})} &&&& CGK(X) \ar[rrrr]^-{\omega_{K(X)}} \ar[dr]^(0.7){\theta^{100}_{GK(X)}} \ar@{-->}[dddd]^(0.6){\varrho_{GK(X)}} &&&& 
HK(X) \ar@{-->}[dddd]^(0.6){\chi_{K(X)}} \ar[dr]^(0.7){\theta^{000}_{K(X)}} &&&\\
& CC(X) \ar[dddd]^{\varrho_{C(X)}} \ar[rrrr]^-{C(\theta^{100}_{X})} \ar[dr]^(0.7){C(\varrho_{X})} &&&& CCC(X) \ar@{-->}[dddd]^(0.6){\varrho_{CC(X)}} \ar[dr]^(0.7){C(\varrho_{C(X)})} 
\ar[rrrr]^-{CC(\varrho_{X})} &&&& CCGK(X) \ar@{-->}[dddd]^(0.6){\varrho_{GKC(X)}} \ar[rrrr]^-{C(\omega_{K(X)})} \ar[dr]^(0.7){C(\varrho_{GK(X)})} &&&& 
CHK(X)  \ar@{-->}[dddd]^(0.6){\varrho_{HK(X)}} \ar[dr]^(0.7){C(\chi_{K(X)})} &&\\
&& CGK(X) \ar[dddd]_(0.6){\varrho_{GK(X)}} \ar[dr]^(0.7){\omega_{K(X)}}  \ar[rrrr]^-{CG(\theta^{110}_{X})} &&&& CGKC(X) \ar@{-->}[dddd]_(0.6){\varrho_{GKC(X)}} \ar[dr]^(0.6){\omega_{KC(X)}} \ar[rrrr]^-{CGK(\varrho_{X})} &&&& CGKGK(X) \ar@{-->}[dddd]^{\varrho_{C(X)}} \ar[rrrr]^-{CG(\varepsilon_{K(X)})} \ar[dr]^(0.7){\omega_{KGK(X)}} &&&& 
CGQK(X)  \ar@{-->}[dddd]^{\varrho_{GQK(X)}} \ar[dr]^(0.7){\omega_{QK(X)}} &\\
&&& HK(X) \ar[dddd]_(0.4){\chi_{K(X)}} \ar[rrrr]^{H(\theta^{110}_{X})} &&&& HKC(X) \ar[dddd]_(0.45){\chi_{KC(X)}} \ar[rrrr]^-{HK(\varrho_{X})} &&&& 
HKGK(X) \ar[rrrr]^-{H(\varepsilon_{K(X)})} \ar[dddd]_(0.45){\chi_{KGK(X)}} &&&& HQK(X) \ar[dddd]^{\chi_{QK(X)}} \\
GK(X) \ar[dr]_(0.3){G(\theta^{110}_{X})} \ar@{-->}[rrrr]^-{G(\theta^{110}_{X})} &&&& 
GKC(X) \ar@{-->}[rrrr]^-{GK(\varrho_{X})} \ar@{-->}[dr]_(0.2){G(\theta^{110}_{C(X)})} &&&& 
GKGK(X) \ar@{-->}[rrrr]^-{G(\varepsilon_{K(X)})} \ar@{-->}[dr]_(0.3){G(\theta^{110}_{GK(X)})} &&&& GQK(X) \ar@{-->}[dr]_(0.2){G(\theta^{010}_{K(X)})} &&&\\
& GKC(X) \ar[dr]_(0.2){GK(\varrho_{X})} \ar@{-->}[rrrr]_-{GK(\theta^{100}_{X})} &&&& 
GKCC(X) \ar@{-->}[dr]_(0.3){GK(\varrho_{C(X)})} \ar@{-->}[rrrr]_-{GKC(\varrho_{X})} &&&& 
GKCGK(X) \ar@{-->}[rrrr]_-{GK(\omega_{K(X)})}  \ar@{-->}[dr]_(0.3){GK(\varrho_{GK(X)})} &&&& GKHK(X) \ar@{-->}[dr]_(0.3){GK(\chi_{K(X)})} && \\
&& GKGK(X) \ar@{-->}[dr]_(0.3){G(\varepsilon_{K(X)})} \ar@{-->}[rrrr]_-{GKG(\theta^{110}_{X})} &&&& GKGKC(X) \ar@{-->}[dr]_(0.3){G(\varepsilon_{KC(X)})} \ar@{-->}[rrrr]_-{GKGK(\varrho_{C(X)})} &&&& GKGKGK(X) \ar@{-->}[rrrr]_-{GKG(\varepsilon_{K(X)})} \ar@{-->}[dr]_(0.3){G(\varepsilon_{KGK(X)})} &&&& GQK(X) \ar@{-->}[dr]_(0.3){G(\varepsilon_{QK(X)})} & \\
&&& GQK(X) \ar[rrrr]_{GQ(\theta^{110}_{X})} &&&& GQKC(X) \ar[rrrr]_-{GQK(\varrho_{C(X)})} &&&& GQKGK(X) \ar[rrrr]_-{GQ(\varepsilon_{K(X)})} &&&& GQQK(X) }\]
commutes.
\begin{proof}
All the cubes in the above diagram commute by definition or naturalness except the first and the last ones in the first raw. The top and bottom squares of the former commute by (\ref{thetaX000}) and (\ref{thetaX110}) respectively and the top and bottom squares of the latter commute by (\ref{thetaAcompeq1}) and (\ref{thetaAcompeq2}) respectively.
By vertically composing the columns we obtain three cubes
\[\xymatrix@!=1.5pc@R0.5pc@C.75pc{C(X) \ar[dddd]_{\varrho_{X}} \ar[dr]^(0.8){\theta^{101}_{X}}  \ar[rrrr]^-{\theta^{100}_{X}} &&&& 
CC(X) \ar@{-->}[dddd]_{\varrho_{C(X)}} \ar[dr]^(0.8){\theta^{101}_{C(X)}} \ar[rrrr]^-{C(\varrho_{X})} &&&& 
CGK(X) \ar@{-->}[dddd]^{\varrho_{GK(X)}} \ar[rrrr]^-{\omega_{K(X)}} \ar[dr]^(0.8){\theta^{101}_{GK(X)}} &&&& 
HK(X) \ar@{-->}[dddd]^{\chi_{K(X)}} \ar[dr]^(0.8){\theta^{001}_{K(X)}} &\\
& HK(X) \ar[dddd]_(0.4){\chi_{K(X)}} \ar[rrrr]^{H(\theta^{110}_{X})} &&&& HKC(X) \ar[dddd]^{\chi_{KC(X)}} \ar[rrrr]^-{HK(\varrho_{X})} &&&& 
HKGK(X) \ar[rrrr]^-{H(\varepsilon_{K(X)})} \ar[dddd]^{\chi_{KGK(X)}} &&&& HQK(X) \ar[dddd]^{\chi_{QK(X)}} \\
& &&&& &&&& &&&& \\
& &&&&  &&&& &&&& \\
GK(X) \ar@{-->}[dr]_(0.2){G(\theta^{111}_{X})} \ar@{-->}[rrrr]^-{G(\theta^{110}_{X})} &&&& 
GKC(X) \ar@{-->}[dr]_(0.2){G(\theta^{111}_{C(X)})} \ar@{-->}[rrrr]^-{GK(\varrho_{X})} &&&& 
GKGK(X) \ar@{-->}[rrrr]^-{G(\varepsilon_{K(X)})} \ar@{-->}[dr]_(0.2){G(\theta^{111}_{GK(X)})} &&&& GQK(X) \ar@{-->}[dr]_(0.2){G(\theta^{011}_{K(X)})} & \\
& GQK(X) \ar[rrrr]_{GQ(\theta^{110}_{X})} &&&& GQKC(X) \ar[rrrr]_-{GQK(\varrho_{C(X)})} &&&& GQKGK(X) \ar[rrrr]_-{GQ(\varepsilon_{K(X)})} &&&& GQQK(X) }\]
in which the top and bottom squares of the first one commute by (\ref{thetaX010}) and (\ref{thetaX100}) respectively, and the top and bottom squares of the last one commute by (\ref{thetaAcompeq3}) and (\ref{thetaAcompeq4}) respectively.  On the other side by if we horizontally compose the raws of the first diagram we obtain  a diagram
\[\hspace{-1cm}\xymatrix@!=1.15pc@R1pc@C1pc{C(X) \ar[dddd]_{\varrho_{X}} \ar[rrrr]^-{\theta^{101}_{X}} \ar[dr]^(0.7){\theta^{100}_{X}} &&&& 
HK(X) \ar@{-->}[dddd]^(0.6){\chi_{K(X)}} \ar[dr]^(0.7){\theta^{000}_{K(X)}}  &&&\\
& CC(X) \ar[dddd]^{\varrho_{C(X)}} \ar[rrrr]^-{C(\theta^{101}_{X})} \ar[dr]^(0.7){C(\varrho_{X})} &&&& CHK(X)  \ar@{-->}[dddd]^(0.6){\varrho_{HK(X)}} \ar[dr]^(0.7){C(\chi_{K(X)})} &&\\
&& CGK(X) \ar[dddd]_(0.6){\varrho_{GK(X)}} \ar[dr]^(0.7){\omega_{K(X)}}  \ar[rrrr]^-{CG(\theta^{111}_{X})} &&&& 
CGQK(X)  \ar@{-->}[dddd]^{\varrho_{GQK(X)}} \ar[dr]^(0.7){\omega_{QK(X)}} &\\
&&& HK(X) \ar[dddd]_(0.4){\chi_{K(X)}} \ar[rrrr]^{H(\theta^{111}_{X})} &&&& HQK(X) \ar[dddd]^{\chi_{QK(X)}} \\
GK(X) \ar[dr]_(0.3){G(\theta^{110}_{X})} \ar@{-->}[rrrr]^-{G(\theta^{111}_{X})} &&&& 
GQK(X) \ar@{-->}[dr]_(0.2){G(\theta^{010}_{K(X)})} &&&\\
& GKC(X) \ar[dr]_(0.2){GK(\theta^{010}_{X})} \ar@{-->}[rrrr]_-{GK(\theta^{101}_{X})} &&&& 
GKHK(X) \ar@{-->}[dr]_(0.3){GK(\chi_{K(X)})} && \\
&& GKGK(X) \ar@{-->}[dr]_(0.3){G(\varepsilon_{K(X)})} \ar@{-->}[rrrr]_-{GKG(\theta^{111}_{X})} &&&& GKGQK(X) \ar@{-->}[dr]_(0.3){G(\varepsilon_{QK(X)})}  & \\
&&& GQK(X) \ar[rrrr]_{GQ(\theta^{111}_{X})} &&&& GQQK(X) }\]
in which the top and bottom squares of the first one commute by (\ref{thetaX001}) and (\ref{thetaX111}) respectively. Then by horizontally and vertically composing the former and the latter diagram respectively, we arrive at the cube
\[\hspace{-1cm}\xymatrix@!=1.2pc@R1pc@C1pc{C(X) \ar[dddd]_(0.6){\varrho_{X}} \ar[dr]^(0.7){\theta^{101}_{X}}  \ar[rrrr]^-{\theta^{101}_{X}} &&&& 
HK(X) \ar@{-->}[dddd]^{\chi_{K(X)}} \ar[dr]^(0.7){\theta^{001}_{K(X)}} &\\
& HK(X) \ar[dddd]_(0.4){\chi_{K(X)}} \ar[rrrr]^{H(\theta^{111}_{X})} &&&& HQK(X) \ar[dddd]^{\chi_{QK(X)}} \\
& &&&& \\
& &&&& \\
GK(X) \ar@{-->}[dr]_(0.3){G(\theta^{111}_{X})} \ar@{-->}[rrrr]^-{G(\theta^{111}_{X})}  
&&&& GQK(X) \ar@{-->}[dr]_(0.3){G(\theta^{011}_{K(X)})} & \\
& GQK(X) \ar[rrrr]_{GQ(\theta^{111}_{X})} &&&& GQQK(X) }\]
whose top and bottom squares commute by (\ref{thetaX011}) and (\ref{thetaX101}). 

\end{proof}
\end{proposition}

\begin{proposition}\label{Aprop}
For any colax $\D$-coalgebra $(G,\bm{F}_{G},\bm{\zeta},\bm{\theta})$ the following diagram
\[\hspace{-1.2cm}\xymatrix@!=1.15pc@R1pc@C1pc{H(A) \ar[dddd]_{\chi_{A}} \ar[rrrr]^-{\theta^{000}_{A}} \ar[dr]^(0.7){\theta^{000}_{A}} &&&& 
CH(A) \ar@{-->}[dddd]^(0.6){\eta_{H(A)}} \ar[dr]^(0.7){\theta^{100}_{H(A)}} \ar[rrrr]^-{C(\chi_{A})} &&&& 
CGQ(A) \ar[rrrr]^-{\omega_{Q(A)}} \ar[dr]^(0.7){\theta^{100}_{GQ(A)}} \ar@{-->}[dddd]^(0.6){\eta_{GQ(A)}} &&&& 
HQ(A) \ar@{-->}[dddd]^(0.6){\chi_{Q(A)}} \ar[dr]^(0.7){\theta^{000}_{Q(A)}} &&&\\
& CH(A) \ar[dddd]_{\eta_{H(A)}} \ar[rrrr]^-{C(\theta^{000}_{A})} \ar[dr]^(0.7){C(\chi_{A})} &&&& 
CCH(A) \ar@{-->}[dddd]_(0.6){\eta_{CH(A)}} \ar[dr]^(0.7){C(\eta_{H(A)})} 
\ar[rrrr]^-{CC(\chi_{A})} &&&& CCGQ(A)  \ar@{-->}[dddd]_(0.6){\chi_{KCG(A)}} \ar[rrrr]^-{C(\omega_{Q(A)})} \ar[dr]^(0.7){C(\eta_{GQ(A)})} &&&& 
CHQ(A) \ar@{-->}[dddd]^(0.6){\eta_{HQ(A)}} \ar[dr]^(0.7){C(\chi_{Q(A)})} &&\\
&& CGQ(A) \ar[dddd]_(0.6){\eta_{GQ(A)}} \ar[dr]^(0.7){\omega_{Q(A)}} \ar[rrrr]^-{CG(\theta^{010}_{A})} &&&& 
CGKH(A) \ar@{-->}[dddd]_(0.6){\chi_{KGKC(X)}} \ar[dr]^(0.6){\omega_{KH(A)}} \ar[rrrr]^-{CGK(\chi_{A})} &&&& 
CGKGQ(A) \ar@{-->}[dddd]_(0.6){\eta_{GKGQ(A)}} \ar[rrrr]^-{CG(\varepsilon_{Q(A)})} \ar[dr]^(0.7){\omega_{KGQ(A)}} &&&& 
CGQQ(A)  \ar@{-->}[dddd]^{\chi_{GQ(A)}} \ar[dr]^(0.7){\omega_{QQ(A)}} &\\
&&& HQ(A) \ar[dddd]_(0.45){\chi_{Q(A)}} \ar[rrrr]^{H(\theta^{010}_{A})} &&&& HKH(A) \ar[dddd]_(0.3){\chi_{KH(A)}} \ar[rrrr]^-{HK(\chi_{A})} &&&& 
HKGQ(A) \ar[rrrr]^-{H(\varepsilon_{Q(A)})} \ar[dddd]_(0.4){\chi_{KGQ(A)}} &&&& HQQ(A) \ar[dddd]^{\chi_{QQ(A)}} \\
GQ(A) \ar[dr]_(0.3){G(\theta^{010}_{A})} \ar@{-->}[rrrr]^-{G(\theta^{010}_{A})} &&&& 
GKH(A) \ar@{-->}[rrrr]^-{GK(\chi_{A})} \ar@{-->}[dr]_(0.2){G(\theta^{110}_{H(A)})} &&&& 
GKGQ(A) \ar@{-->}[rrrr]^-{G(\varepsilon_{Q(A)})} \ar@{-->}[dr]_(0.3){G(\theta^{110}_{GQ(A)})} &&&& GQQ(A) \ar@{-->}[dr]_(0.2){G(\theta^{010}_{Q(A)})} &&&\\
& GKH(A) \ar[dr]_(0.2){GK(\chi_{A})} \ar@{-->}[rrrr]_-{GK(\theta^{000}_{A})} &&&& 
GKCH(A) \ar@{-->}[dr]_(0.3){GK(\eta_{H(A)})} \ar@{-->}[rrrr]_-{GKC(\chi_{A})} &&&& 
GKCGQ(A) \ar@{-->}[rrrr]_-{GK(\omega_{Q(A)})}  \ar@{-->}[dr]_(0.3){GK(\eta_{GQ(A)})} &&&& GKHQ(A) \ar@{-->}[dr]_(0.3){GK(\chi_{Q(A)})} && \\
&& GKGQ(A) \ar@{-->}[dr]_(0.3){G(\varepsilon_{Q(A)})} \ar@{-->}[rrrr]_-{GKGQ(\theta^{010}_{A})} &&&& 
GKGKH(A) \ar@{-->}[dr]_(0.3){G(\varepsilon_{KH(X)})} \ar@{-->}[rrrr]_-{GKGK(\chi_{A})} &&&& 
GKGKGQ(A) \ar@{-->}[rrrr]_-{GKG(\varepsilon_{Q(A)})} \ar@{-->}[dr]_(0.3){G(\varepsilon_{KGQ(A)})} &&&& GKGQQ(A) \ar@{-->}[dr]_(0.3){G(\varepsilon_{QQ(A)})} & \\
&&& GQQ(A) \ar[rrrr]_{GQ(\theta^{010}_{A})} &&&& GQKH(X) \ar[rrrr]_-{GQK(\chi_{A})} &&&& GQKGQ(A) \ar[rrrr]_-{GQ(\varepsilon_{Q(A)})} &&&& GQQQ(A) }\]
commutes.
\begin{proof}
All the cubes of the above diagram commute by definition or naturalness except the first and the last ones in the first raw. The top and bottom squares of the former commute by (\ref{thetaA000}) and (\ref{thetaA110}) respectively and the top and bottom squares of the latter commute by (\ref{thetaAcompeq1}) and (\ref{thetaAcompeq2}) respectively.
By vertically composing the columns we obtain the three cubes
\[\hspace{-1.2cm}\xymatrix@!=1.5pc@R0.5pc@C0.75pc{H(A) \ar[dddd]_{\chi_{A}} \ar[rrrr]^-{\theta^{000}_{A}} \ar[dr]^(0.8){\theta^{001}_{A}} &&&& 
CH(A) \ar@{-->}[dddd]_(0.5){\eta_{H(A)}} \ar[dr]^(0.8){\theta^{101}_{H(A)}} \ar[rrrr]^-{C(\chi_{A})} &&&& 
CGQ(A) \ar[rrrr]^-{\omega_{Q(A)}} \ar[dr]^(0.8){\theta^{101}_{GQ(A)}} \ar@{-->}[dddd]_(0.5){\eta_{GQ(A)}} &&&& 
HQ(A) \ar@{-->}[dddd]^(0.6){\chi_{Q(A)}} \ar[dr]^(0.8){\theta^{001}_{Q(A)}} &\\
& HQ(A) \ar[dddd]_(0.45){\chi_{Q(A)}} \ar[rrrr]^{H(\theta^{010}_{A})} &&&& HKH(A) \ar[dddd]_(0.5){\chi_{KH(A)}} \ar[rrrr]^-{HK(\chi_{A})} &&&& 
HKGQ(A) \ar[rrrr]^-{H(\varepsilon_{Q(A)})} \ar[dddd]^(0.5){\chi_{KGQ(A)}} &&&& HQQ(A) \ar[dddd]^{\chi_{QQ(A)}} \\
&& &&&& &&& \\
&& &&&& &&& \\
GQ(A) \ar@{-->}[dr]_(0.2){G(\theta^{011}_{A})} \ar@{-->}[rrrr]^-{G(\theta^{010}_{A})} &&&& 
GKH(A) \ar@{-->}[dr]_(0.2){G(\theta^{111}_{H(A)})} \ar@{-->}[rrrr]^-{GK(\chi_{A})} &&&& 
GKGQ(A) \ar@{-->}[rrrr]^-{G(\varepsilon_{Q(A)})} \ar@{-->}[dr]_(0.2){G(\theta^{111}_{GQ(A)})} &&&& GQQ(A) \ar@{-->}[dr]_(0.2){G(\theta^{011}_{Q(A)})} & \\
& GQQ(A) \ar[rrrr]_{GQ(\theta^{010}_{A})} &&&& GQKH(X) \ar[rrrr]_-{GQK(\chi_{A})} &&&& GQKGQ(A) \ar[rrrr]_-{GQ(\varepsilon_{Q(A)})} &&&& GQQQ(A) }\]
The top and bottom squares of the first cube commute by (\ref{thetaA010}) and (\ref{thetaA100}) respectively, and the top and bottom squares of the last one commute by (\ref{thetaAcompeq3}) and (\ref{thetaAcompeq4}) respectively.  On the other side if we horizontally compose the raws of the first diagram we get
\[\hspace{-1.2cm}\xymatrix@!=1.15pc@R1pc@C1pc{H(A) \ar[dddd]_{\chi_{A}} \ar[rrrr]^-{\theta^{001}_{A}} \ar[dr]^(0.7){\theta^{000}_{A}} &&&&
HQ(A) \ar@{-->}[dddd]^(0.6){\chi_{Q(A)}} \ar[dr]^(0.7){\theta^{000}_{Q(A)}} &&&\\
& CH(A) \ar[dddd]_{\eta_{H(A)}} \ar[rrrr]^-{C(\theta^{001}_{A})} \ar[dr]^(0.7){C(\chi_{A})} &&&& 
CHQ(A) \ar@{-->}[dddd]^(0.6){\eta_{HQ(A)}} \ar[dr]^(0.7){C(\chi_{Q(A)})} &&\\
&& CGQ(A) \ar[dddd]_(0.6){\eta_{GQ(A)}} \ar[dr]^(0.7){\omega_{Q(A)}} \ar[rrrr]^-{CG(\theta^{011}_{A})} &&&& 
CGQQ(A)  \ar@{-->}[dddd]^{\chi_{GQ(A)}} \ar[dr]^(0.7){\omega_{QQ(A)}} &\\
&&& HQ(A) \ar[dddd]_(0.45){\chi_{Q(A)}} \ar[rrrr]^{H(\theta^{011}_{A})} &&&& HQQ(A) \ar[dddd]^{\chi_{QQ(A)}} \\
GQ(A) \ar[dr]_(0.3){G(\theta^{010}_{A})} \ar@{-->}[rrrr]^-{G(\theta^{011}_{A})} &&&& GQQ(A) \ar@{-->}[dr]_(0.2){G(\theta^{010}_{Q(A)})} &&&\\
& GKH(A) \ar[dr]_(0.2){GK(\chi_{A})} \ar@{-->}[rrrr]_-{GK(\theta^{001}_{A})} &&&& GKHQ(A) \ar@{-->}[dr]_(0.3){GK(\chi_{Q(A)})} && \\
&& GKGQ(A) \ar@{-->}[dr]_(0.3){G(\varepsilon_{Q(A)})} \ar@{-->}[rrrr]_-{GKG(\theta^{011}_{A})} &&&& GKGQQ(A) \ar@{-->}[dr]_(0.3){G(\varepsilon_{QQ(A)})} & \\
&&& GQQ(A) \ar[rrrr]_{GQ(\theta^{011}_{A})} &&&& GQQQ(A) }\]
in which the top and bottom squares of the first one commute by (\ref{thetaA001}) and (\ref{thetaA111}) respectively. By horizontally and vertically composing the last two diagrams respectively, we get
\[\xymatrix@!=1pc@R1pc@C1pc{H(A) \ar[dddd]_{\chi_{A}} \ar[rrrr]^-{\theta^{001}_{A}} \ar[dr]^(0.7){\theta^{001}_{A}} &&&&
HQ(A) \ar@{-->}[dddd]^(0.6){\chi_{Q(A)}} \ar[dr]^(0.7){\theta^{001}_{Q(A)}} &\\
& HQ(A) \ar[dddd]_(0.45){\chi_{Q(A)}} \ar[rrrr]^{H(\theta^{011}_{A})} &&&& HQQ(A) \ar[dddd]^{\chi_{QQ(A)}} \\
&&&&& \\
&&&&& \\
GQ(A) \ar@{-->}[dr]_(0.3){G(\theta^{011}_{A})} \ar@{-->}[rrrr]^-{G(\theta^{011}_{A})} &&&& GQQ(A) \ar@{-->}[dr]_(0.3){G(\theta^{011}_{Q(A)})} & \\
& GQQ(A) \ar[rrrr]_{GQ(\theta^{011}_{A})} &&&& GQQQ(A) }\]
whose top and bottom squares commute by (\ref{thetaA011}) and (\ref{thetaA101}). 

\end{proof}
\end{proposition}

\begin{proposition}\label{KKleisliliftingprop}
There exists a lifing
\begin{equation}\label{KKleislilifting}
\begin{array}{c}
\xymatrix@!=2pc{\X_{C} \ar[d]_{U_{C}} \ar[r]^-{\widetilde{K}} & \A_{Q} \ar[d]^{U_{Q}} \\
\X \ar[r]_-{K} & \A }
\end{array}
\end{equation}
where the functor $\widetilde{K}_{G} \maps \X^{C} \to \A^{Q}$ acts as $K$ on objects and takes any Kleisli morphism $f \maps C(X) \to Y$ to
\begin{equation}\label{KKlifting}
\begin{array}{c}
\xymatrix@!=1pc{QK(X) \ar[rr]^-{Q(\theta^{110}_{X})} && QKC(X) \ar[rr]^-{QK(f)} && QK(Y) \ar[rr]^-{\zeta^{0}_{K(Y)}} && K(Y).}
\end{array}
\end{equation} 
\begin{proof}
We need to show that $\widetilde{K}$ is functorial. This means that 
\[\widetilde{K}(\widetilde{g} \circ_{Kl} \widetilde{f}) = \widetilde{K}(\widetilde{g}) \circ_{Kl} \widetilde{K}(\widetilde{f})\]
The left side of the equation is the composition
\[\xymatrix@!=1pc{QK(X) \ar[rr]^-{Q(\theta^{110}_{X})} && QKC(X) \ar[rr]^-{QK(\theta^{100}_{X})} && QKCC(X) \ar[rr]^-{QKC(f)} && QKC(Y)  \ar[rr]^-{QK(g)} && 
QK(Z) \ar[rr]^-{\zeta^{0}_{K(Z)}} && K(Z) }\]
and by (\ref{thetaX110}) and naturality this composition is equal to 
\[\xymatrix@!=1pc{QK(X) \ar[rr]^-{Q(\theta^{110}_{X})} && QKC(X) \ar[drr]_-{QK(f)} \ar[rr]^-{Q(\theta^{110}_{C(X)})} && QKCC(X) \ar[rr]^-{QKC(f)} && QKC(Y)  \ar[rr]^-{QK(g)} && 
QK(Z) \ar[rr]^-{\zeta^{0}_{K(Z)}} && K(Z) \\
&&&& QK(Y) \ar[urr]_-{Q(\theta^{110}_{Y})} &&&&&& }\]
The right side of the above equation is the composition
\[\hspace{-0.5cm}\xymatrix@!=0.75pc{QK(X) \ar[d]_{\theta^{011}_{K(X)}} && && && && && \\
QQK(X) \ar[rr]^-{QQ(\theta^{110}_{X})} && QQKC(X) \ar[rr]^-{QQK(f)} && QQK(Y) \ar[rr]^-{Q(\zeta^{0}_{K(Y)})} && QK(Y) \ar[rr]^-{Q(\theta^{110}_{Y})} && QKC(Y) \ar[rr]^-{QK(g)} && QK(Z) \ar[rr]^-{\zeta^{0}_{K(Z)}} && K(Z)}\]
and the statement will be proved if the composition of the first four terms of the last string is equal to $QK(f)Q(\theta^{110}_{X})$. This follows from a diagram
\[\xymatrix@!=0.5pc{HK(X) \ar[ddd]_{\chi_{K(X)}} \ar[rrr]^-{H(\theta^{110}_{X})}  \ar[dr]^(0.7){\theta^{001}_{K(X)}} &&& 
HKC(X) \ar@{-->}[ddd]^(0.6){\chi_{KC(X)}} \ar[dr]^(0.7){\theta^{001}_{KC(X)}} \ar[rrr]^-{HK(f)} &&& 
HK(Y) \ar@{-->}[ddd]_{\chi_{K(Y)}} \ar[dr]^(0.7){\theta^{001}_{K(X)}} \ar[rrr]^{1_{HK(Y)}} &&& HK(Y) \ar@{-->}[ddd]_{\chi_{K(Y)}}  \ar[dr]^-{1_{HK(Y)}}  &\\
& HQK(X) \ar[ddd]_{\chi_{QK(X)}} \ar[rrr]^-{HQ(\theta^{110}_{X})} &&& HQKC(X) \ar[ddd]_{\chi_{QKC(X)}} \ar[rrr]^{HQK(f)} &&& 
HQK(Y) \ar[ddd]_{\chi_{QK(Y)}}  \ar[rrr]^{H(\zeta^{0}_{K(Y)})} &&& HK(Y) \ar[ddd]^{\chi_{K(Y)}}   \\
&&&&&&& &&&  \\
GQK(X) \ar[dr]_(0.3){G(\theta^{011}_{K(X)})} \ar@{-->}[rrr]^-{GQ(\theta^{110}_{X})} &&& 
GQKC(X) \ar@{-->}[rrr]^-{GQK(f)} \ar@{-->}[dr]_(0.2){G(\theta^{011}_{KC(X)})} &&& 
GQK(Y) \ar@{-->}[rrr]^-{G(1_{QK(Y)})} \ar@{-->}[dr]_(0.3){GQ(\theta^{011}_{K(Y)})} &&& GQK(Y)  \ar@{-->}[dr]_(0.2){G(1_{QK(Y)})} & \\
& GQQK(X) \ar@{-->}[rrr]_-{GQQ(\theta^{110}_{X})} &&& GQQKC(X) \ar[rrr]_{GQQK(f)} &&& GQQK(Y)  \ar[rrr]_{GQ(\zeta^{0}_{K(Y)})} &&& GQK(Y)   }\]
in which the first two cubes commute by naturality and the top and bottom of the last one by (\ref{leftuniteqA1}) and (\ref{leftuniteqA2}) respectively.
\end{proof}
\end{proposition}

\begin{proposition}\label{HKleisliliftingprop}
There exists a lifing
\begin{equation}\label{Hlifting}
\begin{array}{c}
\xymatrix@!=2pc{\A_{Q} \ar[d]_{U_{Q}} \ar[r]^-{\widetilde{H}} & \X_{C} \ar[d]^{U_{C}} \\
\A \ar[r]_-{H} & \X }
\end{array}
\end{equation}
where the functor $\widetilde{H} \maps\A_{Q} \to \X_{C}$ acts as $H$ on objects and it sends $a \maps Q(A) \to B$ to
\begin{equation}\label{Hliftingdef}
\begin{array}{c}
\xymatrix@!=1pc{CH(A) \ar[rr]^-{C(\theta^{001}_{A})} && CHQ(A) \ar[rr]^-{CH(a)} && CH(B) \ar[rr]^-{\zeta^{1}_{H(B)}} && H(B).}
\end{array}
\end{equation}
\begin{proof}
We need to show that $\widetilde{H}$ is functorial. This means that 
\[\widetilde{H}(\widetilde{b} \circ_{Kl} \widetilde{a}) = \widetilde{H}(\widetilde{b}) \circ_{Kl} \widetilde{H}(\widetilde{a})\]
The left side of the equation is the composition of the string of morphisms
\[\xymatrix@!=1pc{CH(A) \ar[rr]^-{C(\theta^{001}_{A})} && CHQ(A) \ar[rr]^-{CH(\theta^{011}_{A})} && CHQQ(A) \ar[rr]^-{CHQ(a)} && CHQ(B) \ar[rr]^-{CH(b)} && 
CH(W) \ar[rr]^-{\zeta^{1}_{H(W)}} && H(W) }\]
and the right side of the equation is
\[\hspace{-0.5cm}\xymatrix@!=0.65pc{CH(A) \ar[d]_{\theta^{100}_{H(A)}} && && && && && \\
CCH(A) \ar[rr]^-{CC(\theta^{001}_{A})} && CCHQ(A) \ar[rr]^-{CCH(a)} && 
CCH(B) \ar[rr]^-{C(\zeta^{1}_{H(B)})} && CH(B) \ar[rr]^-{C(\theta^{001}_{B})} && CHQ(B) \ar[rr]^-{CH(b)} && CH(W) \ar[rr]^-{\zeta^{1}_{H(W)}} && H(W).}\]
The statement is proved by the following commutative diagram 
\[\xymatrix@!=4pc{CH(A) \ar[r]^-{C(\theta^{001}_{A})} \ar[d]_{\theta^{100}_{H(A)}} & CHQ(A) \ar[r]^-{CH(a)} \ar[d]^{\theta^{100}_{HQ(A)}} & 
CH(B) \ar[d]^{\theta^{100}_{H(B)}} \\
CCH(A) \ar[d]_{CC(\theta^{001}_{A})} \ar[r]^-{CC(\theta^{001}_{A})} & CCHQ(A) \ar[d]^-{CCH(a)} \ar[r]^-{CCH(a)} & CCH(B) \ar[d]^{C(\zeta^{1}_{H(B)})} \\
CCHQ(A) \ar[r]_-{CCH(a)} & CCH(B) \ar[r]_{C(\zeta^{1}_{H(B)})} & CH(B) }\]
whose right side compose to $1_{CH(B)}$ by (\ref{leftuniteqX1}) and therefore establishes an identity between the composition of the first four terms of the last string  and $CH(a)C(\theta^{001}_{A})$ which are the first two terms of the string
\[\xymatrix@!=1pc{CH(A) \ar[rr]^-{C(\theta^{001}_{A})} && CHQ(A) \ar[drr]_-{CH(a)} \ar[rr]^-{C(\theta^{001}_{Q(A)})} && CHQQ(A) \ar[rr]^-{CHQ(a)} && CHQ(B) \ar[rr]^-{CH(b)} && 
CH(W) \ar[rr]^-{\zeta^{1}_{H(W)}} && H(W) \\
&&&& CH(B) \ar[urr]_-{C(\theta^{001}_{B})} &&&&&& }\]
whose composition is by (\ref{thetaA011}) equal to the left side of the first equation of the proof.
\end{proof}
\end{proposition}
Although it follows from Proposition \ref{KKleisliliftingprop} and Proposition \ref{HKleisliliftingprop} that any colax $\D$-coalgebra $(G,\bm{F}_{G},\bm{\zeta},\bm{\theta})$ induces lifings (\ref{KKlifting}) and (\ref{Hlifting}) of functors $K$ and $H$ to corresponding Kleisli categories it does not necessarily induce the lifings to corresponding Eilenberg-Moore categories. We will now examine the conditions when this is the case. 

\begin{definition}\label{splitcolagebra}
We say that a colax $\D$-coalgebra $(G,\bm{F}_{G},\bm{\zeta},\bm{\theta})$ is split if satisfies the following conditions
\begin{equation}\label{splitrho}
K(\theta^{100}_{X}) = \theta^{110}_{C(X)},
\end{equation}
\begin{equation}\label{splitzeta}
\zeta^{0}_{K(X)} \epsilon_{K(X)} K(\varrho_{X}) = K(\zeta^{1}_{X}).
\end{equation}
\end{definition}

\begin{proposition}\label{splitrhofactor}
For any split colax $\D$-coalgebra $(G,\bm{F}_{G},\bm{\zeta},\bm{\theta})$ we have the following factorization
\begin{equation}\label{splitrhofact}
G(\zeta^{0}_{K(X)})\chi_{K(X)}\theta^{101}_{X} = \varrho_{X}.
\end{equation}
\begin{equation}\label{splitchifact}
G(\zeta^{0}_{Q(A)})\chi_{Q(A)}\theta^{001}_{A} = \chi_{A}.
\end{equation}
\begin{proof}
Consider the following commutative diagram in $(\X,G)$
\[\xymatrix@!=0.5pc{C(X) \ar[ddd]_{\varrho_{X}} \ar[rrr]^-{\theta^{100}_{X}} \ar[ddrr]^{1_{C(X)}} &&& 
CC(X) \ar@{-->}[ddd]_{\varrho_{C(X)}} \ar[dr]^{C(\varrho_{X})} \ar[rrr]^{\zeta^{1}_{C(X)}} \ar@{-->}[ddd]_{\varrho_{C(X)}} &&& 
C(X) \ar@{-->}[ddd]^{\varrho_{X}} \ar[dr]^{\varrho_{X}} &&\\
&&&& CGK(X) \ar@{-->}[ddd]_{\varrho_{GK(X)}} \ar[dr]^{\omega_{K(X)}}  \ar[rrr]^{\zeta^{1}_{GK(X)}} &&& GK(X) \ar@{-->}[ddd]_{1_{GK(X)}} \ar[dr]^{1_{GK(X)}} &\\
&& C(X) \ar[ddd]_{\varrho_{X}} \ar[rrr]^-{\theta^{101}_{X}} &&& HK(X) \ar[ddd]_{\chi_{K(X)}} \ar[rrr]^{G(\zeta^{0}_{K(X)})\chi_{K(X)}} &&& GK(X) \ar[ddd]^{1_{GK(X)}}    \\
GK(X) \ar[ddrr]_{G(1_{K(X)})} \ar@{-->}[rrr]^-{G(\theta^{110}_{X})} &&& GKC(X) \ar@{-->}[dr]_(0.3){GK(\varrho_{X})}  \ar@{-->}[rrr]^{GK(\zeta^{1}_{X})} &&& 
GK(X) \ar@{-->}[dr]_(0.3){G(1_{K(X)})}  && \\
&&&& GKGK(X) \ar@{-->}[dr]_(0.3){G(\varepsilon_{K(X)})}  \ar@{-->}[rrr]_{G(\zeta^{0}_{K(X)}\varepsilon_{K(X)})} &&& GK(X) \ar@{-->}[dr]_(0.3){G(1_{K(X)})} & \\
&& GK(X) \ar[rrr]_{G(\theta^{111}_{X})} &&& GQK(X) \ar[rrr]_{G(\zeta^{0}_{K(X)})} &&& GK(X) }\]
in which the lower right cube commutes by (\ref{colaxcoh2}). Its composition is equal to
\[\xymatrix@!=0.5pc{C(X) \ar[ddd]_{\varrho_{X}} \ar[dr]^{1_{C(X)}} \ar[rrr]^-{1_{C(X)}} &&& C(X) \ar[dr]^{\varrho_{X}} \ar[ddd]^(0.6){\varrho_{X}} & \\
& C(X) \ar[ddd]_(0.4){\varrho_{X}} \ar[rrr]_-{G(\zeta^{0}_{K(X)})\chi_{K(X)}\theta^{101}_{X}} &&& GK(X) \ar[ddd]^{1_{GK(X)}} \\
&&&&\\
GK(X)  \ar[dr]_-{G(1_{K(X)})} \ar@{-->}[rrr]^-{G(1_{K(X)})} &&& GK(X) \ar@{-->}[dr]_(0.4){G(1_{K(X)})} & \\
& GK(X) \ar[rrr]_{G(1_{K(X)})} &&& GK(X)}\]
since the horizontal compositions in the back square are identities by (\ref{rightuniteqX1}) and (\ref{leftuniteqX2}),  and the lower composition in the front square is identity by (\ref{rightuniteqX2}). The statement follows from the fact that the front square is well defined morphism in $(\X,G)$. The similar argument shows the validity of (\ref{splitchifact}).
\end{proof}
\end{proposition}
Split colax $\D$-coalgebras have remarkable self-generating structure as in a diagram 
\[\hspace{-1.5cm}\xymatrix@!=1pc@R1pc@C1pc{C(X) \ar[dddd]_{\theta^{101}_{X}} \ar[rrrr]^-{\theta^{100}_{X}} \ar[dr]^(0.7){\theta^{100}_{X}} &&&& 
CC(X) \ar@{-->}[dddd]^(0.6){\theta^{101}_{C(X)}} \ar[dr]^(0.7){\theta^{100}_{C(X)}} \ar[rrrr]^-{C(\varrho_{X})} &&&& 
CGK(X) \ar[rrrr]^-{\omega_{K(X)}} \ar[dr]^(0.7){\theta^{100}_{GK(X)}} \ar@{-->}[dddd]^(0.6){\theta^{101}_{GK(X)}} &&&& 
HK(X) \ar@{-->}[dddd]^(0.6){\theta^{001}_{K(X)}} \ar[dr]^(0.7){\theta^{000}_{K(X)}} &&&\\
& CC(X) \ar[dddd]_{\theta^{101}_{C(X)}} \ar[rrrr]^-{C(\theta^{100}_{X})} \ar[dr]^(0.7){C(\varrho_{X})} &&&& 
CCC(X) \ar@{-->}[dddd]^(0.6){\theta^{101}_{GQK(X)}} \ar[dr]^(0.7){C(\varrho_{C(X)})} 
\ar[rrrr]^-{CC(\varrho_{X})} &&&& CCGK(X) \ar@{-->}[dddd]^(0.6){\varrho_{GKC(X)}} \ar[rrrr]^-{C(\omega_{K(X)})} \ar[dr]^(0.7){C(\varrho_{GK(X)})} &&&& 
CHK(X)  \ar@{-->}[dddd]^(0.6){\theta^{101}_{HK(X)}} \ar[dr]^(0.7){C(\chi_{K(X)})} &&\\
&& CGK(X) \ar[dddd]_(0.6){\theta^{101}_{GK(X)}} \ar[dr]^(0.7){\omega_{K(X)}}  \ar[rrrr]^-{CG(\theta^{110}_{X})} &&&& 
CGKC(X) \ar@{-->}[dddd]_(0.6){\varrho_{GKC(X)}} \ar[dr]^(0.6){\omega_{KC(X)}} \ar[rrrr]^-{CGK(\varrho_{X})} &&&& 
CGKGK(X) \ar@{-->}[dddd]^{\varrho_{C(X)}} \ar[rrrr]^-{CG(\varepsilon_{K(X)})} \ar[dr]^(0.7){\omega_{KGK(X)}} &&&& 
CGQK(X)  \ar@{-->}[dddd]^{\theta^{101}_{GQK(X)}} \ar[dr]^(0.7){\omega_{QK(X)}} &\\
&&& HK(X) \ar[dddd]_(0.4){\theta^{001}_{K(X)}} \ar[rrrr]^{H(\theta^{110}_{X})} &&&& HKC(X) \ar[dddd]_(0.45){\theta^{001}_{KC(X)}} \ar[rrrr]^-{HK(\varrho_{X})} &&&& 
HKGK(X) \ar[rrrr]^-{H(\varepsilon_{K(X)})} \ar[dddd]_(0.45){\theta^{001}_{KGK(X)}} &&&& HQK(X) \ar[dddd]^{\theta^{001}_{QK(X)}} \\
HK(X) \ar[dddd]_{\chi_{K(X)}} \ar[dr]_(0.3){H(\theta^{110}_{X})} \ar@{-->}[rrrr]^-{H(\theta^{110}_{X})} &&&& 
HKC(X) \ar@{-->}[dddd]_{\chi_{KC(X)}} \ar@{-->}[rrrr]^-{HK(\varrho_{X})} \ar@{-->}[dr]_(0.2){H(\theta^{110}_{C(X)})} &&&& 
HKGK(X) \ar@{-->}[dddd]_{\chi_{KGK(X)}} \ar@{-->}[rrrr]^-{H(\varepsilon_{K(X)})} \ar@{-->}[dr]_(0.2){H(\theta^{110}_{GK(X)})} &&&& 
HQK(X) \ar@{-->}[dddd]_{\chi_{QK(X)}} \ar@{-->}[dr]_(0.2){H(\theta^{010}_{K(X)})} &&&\\
& HKC(X) \ar[dddd]_{\chi_{KC(X)}} \ar[dr]_(0.2){HK(\varrho_{X})} \ar@{-->}[rrrr]_-{HK(\theta^{100}_{X})} &&&& 
HKCC(X) \ar@{-->}[dddd]_{\chi_{KCC(X)}} \ar@{-->}[dr]_(0.3){HK(\varrho_{C(X)})} \ar@{-->}[rrrr]_-{HKC(\varrho_{X})} &&&& 
HKCHK(X) \ar@{-->}[dddd]_{\chi_{KCHK(X)}} \ar@{-->}[rrrr]_-{HK(\omega_{K(X)})}  \ar@{-->}[dr]_(0.3){HK(\varrho_{HK(X)})} &&&& 
HKHK(X) \ar@{-->}[dddd]_{\chi_{KHK(X)}} \ar@{-->}[dr]_(0.2){HK(\chi_{K(X)})} && \\
&& HKGK(X) \ar[dddd]_(0.6){\chi_{KGK(X)}} \ar@{-->}[dr]_(0.2){H(\varepsilon_{K(X)})} \ar@{-->}[rrrr]_-{HKH(\theta^{110}_{X})} &&&& 
HKHKC(X) \ar@{-->}[dddd]_{\chi_{KHKC(X)}} \ar@{-->}[dr]_(0.2){H(\varepsilon_{KC(X)})} \ar@{-->}[rrrr]_-{HKHK(\varrho_{C(X)})} &&&& 
HKHKHK(X) \ar@{-->}[dddd]_{\chi_{KHKHK(X)}} \ar@{-->}[rrrr]_-{HKH(\varepsilon_{K(X)})} \ar@{-->}[dr]_(0.2){H(\varepsilon_{KHK(X)})} &&&& 
HKGQK(X) \ar@{-->}[dddd]^{\chi_{KGQK(X)}} \ar@{-->}[dr]_(0.2){H(\varepsilon_{QK(X)})} & \\
&&& HQK(X) \ar[dddd]_{\chi_{QK(X)}} \ar[rrrr]_{HQ(\theta^{110}_{X})} &&&& HQKC(X) \ar[dddd]_(0.4){\chi_{QKC(X)}} \ar[rrrr]_-{GQK(\varrho_{C(X)})} &&&& 
HQKGK(X) \ar[dddd]_(0.4){\chi_{QKGK(X)}} \ar[rrrr]_-{HQ(\varepsilon_{K(X)})} &&&& HQQK(X) \ar[dddd]^{\chi_{QQK(X)}} \\
GQK(X) \ar[dddd]_{G(\zeta^{0}_{K(X)})} \ar[dr]_(0.3){GQ(\theta^{110}_{X})} \ar@{-->}[rrrr]^-{GQ(\theta^{110}_{X})} &&&& 
GQKC(X) \ar@{-->}[dddd]_{G(\zeta^{0}_{KC(X)})} \ar@{-->}[rrrr]^-{GQK(\varrho_{X})} \ar@{-->}[dr]_(0.2){GQ(\theta^{110}_{C(X)})} &&&& 
GQKGK(X) \ar@{-->}[dddd]_{G(\zeta^{0}_{KGK(X)})} \ar@{-->}[rrrr]^-{GQ(\varepsilon_{K(X)})} \ar@{-->}[dr]_(0.3){GQ(\theta^{110}_{GK(X)})} &&&& 
GQQK(X) \ar@{-->}[dddd]_{G(\zeta^{0}_{QK(X)})} \ar@{-->}[dr]_(0.2){GQ(\theta^{010}_{K(X)})} &&&\\
& GQKC(X) \ar[dddd]_(0.45){G(\zeta^{0}_{KC(X)})} \ar[dr]_(0.2){GQK(\varrho_{X})} \ar@{-->}[rrrr]_-{GQK(\theta^{100}_{X})} &&&& 
GQKCC(X) \ar@{-->}[dddd]_{} \ar@{-->}[dr]_(0.3){GQK(\varrho_{C(X)})} \ar@{-->}[rrrr]_-{GQKC(\varrho_{X})} &&&& 
GQKCGK(X) \ar@{-->}[dddd]_{}  \ar@{-->}[rrrr]_-{GQK(\omega_{K(X)})} \ar@{-->}[dr]_(0.3){GQK(\varrho_{GK(X)})} &&&& 
GQKHK(X) \ar@{-->}[dddd]_{} \ar@{-->}[dr]_(0.3){GQK(\chi_{K(X)})} && \\
&& GQKGK(X) \ar[dddd]_(0.3){G(\zeta^{0}_{KGK(X)})} \ar@{-->}[dr]_(0.3){GQ(\varepsilon_{K(X)})} \ar@{-->}[rrrr]_-{GQKG(\theta^{110}_{X})} &&&& 
GQKGKC(X) \ar@{-->}[dddd]_(0.4){G(\zeta^{0}_{KGKC(X)})} \ar@{-->}[dr]_(0.3){GQ(\varepsilon_{KC(X)})} \ar@{-->}[rrrr]_-{GQKGK(\varrho_{C(X)})} &&&& 
GQKGKGK(X) \ar@{-->}[dddd]^{G(\zeta^{0}_{KGKGK(X)})} \ar@{-->}[rrrr]_-{GQKG(\varepsilon_{K(X)})} \ar@{-->}[dr]_(0.3){GQ(\varepsilon_{KGK(X)})} &&&& 
GQKGQK(X) \ar@{-->}[dddd]^{G(\zeta^{0}_{KGQK(X)})} \ar@{-->}[dr]_(0.3){GQ(\varepsilon_{QK(X)})} & \\
&&& GQQK(X) \ar[dddd]_(0.4){G(\zeta^{0}_{QK(X)})} \ar[rrrr]_{GQQ(\theta^{110}_{X})} &&&& GQQKC(X) \ar[dddd]_{G(\zeta^{0}_{KC(X)})} \ar[rrrr]_-{GQQK(\varrho_{C(X)})} &&&& 
GQQKGK(X) \ar[dddd]_{} \ar[rrrr]_-{GQQ(\varepsilon_{K(X)})} &&&& GQQQK(X) \ar[dddd]^{G(\zeta^{0}_{QQK(X)})} \\
GK(X) \ar[dr]_(0.3){G(\theta^{110}_{X})} \ar@{-->}[rrrr]^-{G(\theta^{110}_{X})} &&&& 
GKC(X)\ar@{-->}[rrrr]^-{GK(\varrho_{X})} \ar@{-->}[dr]_(0.2){G(\theta^{110}_{C(X)})} &&&& 
GKGK(X) \ar@{-->}[rrrr]^-{G(\varepsilon_{K(X)})} \ar@{-->}[dr]_(0.3){G(\theta^{110}_{GK(X)})} &&&& GQK(X) \ar@{-->}[dr]_(0.2){G(\theta^{010}_{K(X)})} &&&\\
& GKC(X) \ar[dr]_(0.2){GK(\varrho_{X})} \ar@{-->}[rrrr]_-{GK(\theta^{100}_{X})} &&&& 
GKCC(X) \ar@{-->}[dr]_(0.3){GK(\varrho_{C(X)})} \ar@{-->}[rrrr]_-{GKC(\varrho_{X})} &&&& 
GKCGK(X) \ar@{-->}[rrrr]_-{GK(\omega_{K(X)})}  \ar@{-->}[dr]_(0.3){GK(\varrho_{GK(X)})} &&&& GKHK(X) \ar@{-->}[dr]_(0.3){GK(\chi_{K(X)})} && \\
&& GKGK(X) \ar@{-->}[dr]_(0.3){G(\varepsilon_{K(X)})} \ar@{-->}[rrrr]_-{GKG(\theta^{110}_{X})} &&&& GKGKC(X) \ar@{-->}[dr]_(0.3){G(\varepsilon_{KC(X)})} \ar@{-->}[rrrr]_-{GKGK(\varrho_{C(X)})} &&&& GKGKGK(X) \ar@{-->}[rrrr]_-{GKG(\varepsilon_{K(X)})} \ar@{-->}[dr]_(0.3){G(\varepsilon_{KGK(X)})} &&&& 
GKGQK(X) \ar@{-->}[dr]_(0.3){G(\varepsilon_{QK(X)})} & \\
&&& GQK(X) \ar[rrrr]_{GQ(\theta^{110}_{X})} &&&& GQKC(X) \ar[rrrr]_-{GQK(\varrho_{C(X)})} &&&& GQKGK(X) \ar[rrrr]_-{GQ(\varepsilon_{K(X)})} &&&& GQQK(X) }\]
which is obtained from the diagram in the statement of Proposition \ref{Xprop} by factorization of vertical edges dictated by equations (\ref{splitrhofact}) and (\ref{splitchifact}). 
\begin{remark}\label{3x3}
It cannot be a coincidence that such structures replicate themselves as Rubik's cubes. Note that in the above diagram we did not factorized the horizontal component of $\varrho$ in the same manner as we did with its vertical counterpart, but the reader can now easily imagine how this structures replicate themselves in a consistent and coherent manner. There is no any doubt in the mind of the author that a higher (combinatorial) principle is a driving force of such transformations but this remains to be verified in the future.
\end{remark}
It is well known (see \cite{Street1} for the dual statement) that there is a bijective correspondence between the liftings of $K$ to Eilenberg-Moore categories on the left
\begin{equation}\label{gendist}
\begin{array}{c}
\xymatrix@!=3pc{\X^{C} \ar[d]_{U^{C}} \ar[r]^-{\widetilde{K}} & \A^{Q} \ar[d]^{U^{Q}} \\
\X \ar[r]_-{K} & \A }
\hspace{2cm}
\xymatrix@!=3pc{\A \ar[r]^-{Q} & \A \dltwocell<\omit>{\kappa} \\
\X \ar[u]^{K} \ar[r]_-{C}& \X \ar[u]_{K}}
\end{array}
\end{equation}
and natural transformations on the right side which satisfy the following axioms:
\begin{itemize}
\item a commutative diagram
\begin{equation}\label{gendist1}
\begin{array}{c}
\xymatrix@!=3pc{KC \ar[d]_{K \zeta} \ar[r]^-{\kappa} & QK \ar[d]^{\zeta K} \\
K \ar@{=}[r]_-{} & K}
\end{array}
\end{equation}

\item a commutative diagram
\begin{equation}\label{gendist2}
\begin{array}{c}
\xymatrix@!=3pc{KC \ar[d]_{K \tau} \ar[rr]^-{\kappa} && QK \ar[d]^{\tau K} \\
KCC \ar[r]_-{\kappa C} & QKC \ar[r]_-{Q \kappa} & QQK}
\end{array}
\end{equation}
\end{itemize}
which we call {\it generalized distributive law} from $C$ to $Q$.  We have the following:

\begin{theorem}\label{kappadist}
Any split colax $\D$-coalgebra $(G,\bm{F}_{G},\bm{\zeta},\bm{\theta})$ induces a generalized distributive law from $C$ to $Q$ with a natural transformation $\kappa \maps KC \Rightarrow QK$ whose component indexed by an object $X$ in $\X$ is defined by a composition
\begin{equation}\label{kappacomp}
\begin{array}{c}
\xymatrix@!=1pc{KC(X) \ar[rr]^-{K(\varrho_{X})} && KGK(X) \ar[rr]^-{\varepsilon_{K(X)}} && QK(X).}
\end{array}
\end{equation}

\begin{proof}
The following diagram 
\[\xymatrix@!=4pc{KC(X) \ar[d]_-{K(\varrho_{X})} \ar[r]^-{K(\theta^{100}_{X})}_-{\theta^{110}_{C(X)}} & KCC(X) \ar[d]_-{KC(\varrho_{X})} \ar[r]^-{K(\varrho_{C(X)})} & KGKC(X) \ar[r]^-{\varepsilon_{KC(X)}} \ar[d]_-{KGK(\varrho_{X})} & QKC(X) \ar[d]^{QK(\varrho_{X})} \\
KGK(X) \ar[d]_-{\varepsilon_{K(X)}} \ar[r]_-{\theta^{110}_{GK(X)}} & KCGK(X) \ar[r]_-{K(\varrho_{GK(X)})} & KGKGK(X) \ar[r]_-{\varepsilon_{KGK(X)}} & QKGK(X) \ar[d]^-{Q(\varepsilon_{K(X)})} \\
QK(X) \ar[rrr]_-{\theta^{011}_{K(X)}} &&& QQK(X) }\]
shows that the component (\ref{kappacomp}) satisfies the axiom (\ref{gendist2}). The other condition is satisfied by (\ref{splitzeta}).
\end{proof}
\end{theorem}
Now we give an explicit construction of the lifing to Eilenberg-Moore categories.

\begin{theorem}\label{EMKlifing}
For any split colax $\D$-coalgebra $(G,\bm{F}_{G},\bm{\zeta},\bm{\theta})$ there exists a lifting
\[\xymatrix@!=3pc{\X^{C} \ar[d]_{U^{C}} \ar[r]^-{\widetilde{K}^{G}} & \A^{Q} \ar[d]^{U^{Q}} \\
\X \ar[r]_-{K} & \A }\]
where the functor $\widetilde{K}^{G}$ takes any $C$-coalgebra $(X,x)$ to $(K(X),\widetilde{K}^{G}(x))$ where $\widetilde{K}^{G}(x)$ is defined as a composition
\begin{equation}\label{EMKlifingdef}
\begin{array}{c}
\xymatrix@!=1pc{K(X) \ar[rr]^-{K(x)} && KC(X) \ar[rr]^-{K(\varrho_{X})} && KGK(X) \ar[rr]^-{\varepsilon_{K(X)}} && QK(X).}
\end{array}
\end{equation}

\begin{proof}
The following diagram
\[\hspace{-1cm}\xymatrix@!=1.15pc@R1pc@C1pc{C(X) \ar[dddd]_{\varrho_{X}} \ar[rrrr]^-{C(x)} \ar[dr]^(0.7){C(x)} &&&& CC(X) \ar@{-->}[dddd]^(0.6){\varrho_{C(X)}} \ar[dr]^(0.7){C(\theta^{100}_{X})} \ar[rrrr]^-{C(\theta^{101}_{X})} &&&& CHK(X) \ar[rrrr]^-{C(\chi_{K(X)})} \ar[dr]^(0.7){C(\theta^{001}_{K(X)})} \ar@{-->}[dddd]^(0.6){\varrho_{HK(X)}} &&&& 
CGQK(X) \ar@{-->}[dddd]^(0.6){\chi_{K(X)}} \ar[dr]^(0.7){\theta^{000}_{K(X)}} &&&\\
& CC(X) \ar[dddd]^{\varrho_{C(X)}} \ar[rrrr]^-{CC(x)} \ar[dr]^(0.7){C(\varrho_{X})} &&&& CCC(X) \ar@{-->}[dddd]^(0.6){\varrho_{CC(X)}} \ar[dr]^(0.7){C(\varrho_{C(X)})} 
\ar[rrrr]^-{CC(\theta^{101}_{X})} &&&& CCHK(X) \ar@{-->}[dddd]^(0.6){\varrho_{GKC(X)}} \ar[rrrr]^-{C(\omega_{K(X)})} \ar[dr]^(0.7){C(\varrho_{GK(X)})} &&&& 
CHK(X)  \ar@{-->}[dddd]^(0.6){\varrho_{HK(X)}} \ar[dr]^(0.7){C(\chi_{K(X)})} &&\\
&& CGK(X) \ar[dddd]_(0.6){\varrho_{GK(X)}} \ar[dr]^(0.7){\omega_{K(X)}}  \ar[rrrr]^-{CG(\theta^{110}_{X})} &&&& CGKC(X) \ar@{-->}[dddd]^(0.6){\varrho_{GKC(X)}} \ar[dr]^(0.6){\omega_{KC(X)}} \ar[rrrr]^-{CGK(\varrho_{X})} &&&& CGKGK(X) \ar@{-->}[dddd]^{\varrho_{C(X)}} \ar[rrrr]^-{CG(\varepsilon_{K(X)})} \ar[dr]^(0.7){\omega_{KGK(X)}} &&&& 
CGQK(X)  \ar@{-->}[dddd]^{\varrho_{GQK(X)}} \ar[dr]^(0.7){\omega_{QK(X)}} &\\
&&& HK(X) \ar[dddd]_(0.4){\chi_{K(X)}} \ar[rrrr]^{H(\theta^{110}_{X})} &&&& HKC(X) \ar[dddd]_(0.45){\chi_{KC(X)}} \ar[rrrr]^-{HK(\varrho_{X})} &&&& 
HKGK(X) \ar[rrrr]^-{H(\varepsilon_{K(X)})} \ar[dddd]_(0.45){\chi_{KGK(X)}} &&&& HQK(X) \ar[dddd]^{\chi_{QK(X)}} \\
GK(X) \ar[dr]_(0.3){GK(x)} \ar@{-->}[rrrr]^-{GK(x)} &&&& 
GKC(X) \ar@{-->}[rrrr]^-{GK(\varrho_{X})} \ar@{-->}[dr]_(0.2){GK(\theta^{100}_{X})}^(0.8){G(\theta^{110}_{C(X)})} &&&& 
GKGK(X) \ar@{-->}[rrrr]^-{G(\varepsilon_{K(X)})} \ar@{-->}[dr]^(0.8){G(\theta^{110}_{GK(X)})} &&&& GQK(X) \ar@{-->}[dr]_(0.2){G(\theta^{010}_{K(X)})} &&&\\
& GKC(X) \ar[dr]_(0.2){GK(\varrho_{X})} \ar@{-->}[rrrr]_-{GKC(x)} &&&& 
GKCC(X) \ar@{-->}[dr]_(0.3){GK(\varrho_{C(X)})} \ar@{-->}[rrrr]_-{GKC(\varrho_{X})} &&&& 
GKCGK(X) \ar@{-->}[rrrr]_-{GK(\omega_{K(X)})}  \ar@{-->}[dr]_(0.3){GK(\varrho_{GK(X)})} &&&& GKHK(X) \ar@{-->}[dr]_(0.3){GK(\chi_{K(X)})} && \\
&& GKGK(X) \ar@{-->}[dr]_(0.3){G(\varepsilon_{K(X)})} \ar@{-->}[rrrr]_-{GKG(x)} &&&& GKGKC(X) \ar@{-->}[dr]_(0.3){G(\varepsilon_{KC(X)})} \ar@{-->}[rrrr]_-{GKGK(\varrho_{X})} &&&& GKGKGK(X) \ar@{-->}[rrrr]_-{GKG(\varepsilon_{K(X)})} \ar@{-->}[dr]_(0.3){G(\varepsilon_{KGK(X)})} &&&& GQK(X) \ar@{-->}[dr]_(0.3){G(\varepsilon_{QK(X)})} & \\
&&& GQK(X) \ar[rrrr]_{GQK(x)} &&&& GQKC(X) \ar[rrrr]_-{GQK(\varrho_{X})} &&&& GQKGK(X) \ar[rrrr]_-{GQ(\varepsilon_{K(X)})} &&&& GQQK(X) }\]
commutes because it is identical to the diagram in the statement of Proposition \ref{Xprop} except for the first cube where we used a condition (\ref{splitrho}). Its base is a diagram
\[\xymatrix@!=4.5pc{K(X) \ar[r]^-{K(x)} \ar[d]_-{K(x)} & KC(X) \ar[d]^{\theta^{110}_{C(X)}}_{K(\theta^{100}_{X})} \ar[r]^-{K(\eta_{X})} & KGK(X) \ar[d]^{\theta^{110}_{GK(X)}} \ar[r]^-{\epsilon_{K(X)}} & 
QK(X) \ar[d]^{\theta^{010}_{K(X)}} \\
KC(X) \ar[d]_-{K(\eta_{X})} \ar[r]^-{KC(x)} & KCC(X) \ar[d]_-{K(\eta_{C(X)})} \ar[r]^-{KC(\eta_{X})}  & KCGK(X) \ar[r]^-{K(\omega_{K(X)})} \ar[d]_-{K(\eta_{GK(X)})} & KHK(X)  \ar[d]^{K(\chi_{K(X)})} \\
KGK(X) \ar[d]_{\epsilon_{K(X)}} \ar[r]_-{KGK(x)} & KGKC(X) \ar[d]_{\epsilon_{KC(X)}} \ar[r]_-{KGK(\eta_{X})} & KGKGK(X) \ar[d]_-{\epsilon_{KGK(X)}} \ar[r]_-{KG(\epsilon_{K(X)})} & 
KGQK(X) \ar[d]^{\epsilon_{QK(X)}} \\
QK(X) \ar[r]_-{QK(x)} & QKC(X) \ar[r]_-{QK(\eta_{X})} & QKGK(X)  \ar[r]_-{Q(\epsilon_{K(X)})} & QQK(X) }\]
which is exactly the condition which says that $(K(X),\widetilde{K}^{G}(x))$ is a $Q$-coalgebra.
\end{proof}
\end{theorem}
We refer the reader to definition of adjoint squares which appear as Definition I.6.7. of \cite{Gray}.
\begin{theorem}\label{EMKadjointsquare}
There exists an adjoint square
\begin{equation}\label{EMKadjointsquarediag}
\begin{array}{c}
\xymatrix@!=1pc@R1pc@C1pc{\X_{C} \ar@<-1ex>[ddd]_{U_{C}} \ar[rrr]^-{J_{C}} &&& \X^{C} \ar@<-1ex>[ddd]_{U_{C}} \\
& \vartheta^{11} & \vartheta^{12} & \\
& \vartheta^{21} & \vartheta^{22} & \\
\X \ar@<-1ex>[uuu]_{F_{C}}^{\dashv} \ar[rrr]_-{I_{\X}} &&& \X \ar@<-1ex>[uuu]_{F^{C}}^{\dashv} }
\end{array}
\end{equation}
where the comparison functor $J_{C} \maps \X_{C} \to \X^{C}$ takes any object $X$ in $\X_{C} $ to a free $C$-coalgebra $(C(X),\theta^{100}_{X})$ in $\X^{C}$  together with natural transformations
\[\begin{array}{c}
\vartheta^{11} \maps U^{C}J_{C} \Rightarrow I_{\X} U^{C}\\
\vartheta^{12} \maps U^{C}J_{C} F_{C} \Rightarrow I_{\X}\\
\vartheta^{21} \maps J_{C} \Rightarrow F^{C} I_{\X} U_{C}\\
\vartheta^{22} \maps J_{C} F_{C} \Rightarrow F^{C} I_{\X}.
\end{array}\]
defined by
\begin{equation}\label{EMKvarphi11}
\vartheta^{11}_{X} := 1_{C(X)},
\end{equation}
\begin{equation}\label{EMKvarphi12}
\vartheta^{12}_{X} := \zeta^{1}_{X},
\end{equation}
\begin{equation}\label{EMKvarphi21}
\vartheta^{21}_{X} := (\theta^{100}_{X},\theta^{100}_{X}), 
\end{equation}
\begin{equation}\label{EMKvarphi22}
\vartheta^{22}_{X} := (1_{C(X)},1_{C(X)}).
\end{equation}
\begin{proof}
The proof is a straighforward consequence of the axioms for a comonad $C$.

\end{proof}
\end{theorem}

\begin{theorem}\label{EMKadjointsquareC}
There exists an adjoint square
\begin{equation}\label{EMKadjointsquarediagC}
\begin{array}{c}
\xymatrix@!=1pc@R1pc@C1pc{\X_{C} \ar@<-1ex>[ddd]_{U_{C}} \ar[rrr]^-{J_{C}} &&& \X^{C} \ar@<-1ex>[ddd]_{U_{C}} \\
& \vartheta^{11} & \vartheta^{12} & \\
& \vartheta^{21} & \vartheta^{22} & \\
\X \ar@<-1ex>[uuu]_{F_{C}}^{\dashv} \ar[rrr]_-{C} &&& \X \ar@<-1ex>[uuu]_{F^{C}}^{\dashv} }
\end{array}
\end{equation}
where the comparison functor $J_{C} \maps \X_{C} \to \X^{C}$ takes any object $X$ in $\X_{C} $ to a free $C$-colagebra $(C(X),\theta^{100}_{X})$ in $\X^{C}$  together with natural transformations
\[\begin{array}{c}
\vartheta^{11} \maps U^{C}J_{C} \Rightarrow C U^{C}\\
\vartheta^{12} \maps U^{C}J_{C} F_{C} \Rightarrow C\\
\vartheta^{21} \maps J_{C} \Rightarrow F^{C} C U_{C}\\
\vartheta^{22} \maps J_{C} F_{C} \Rightarrow F^{C} C 
\end{array}\]
defined by
\begin{equation}\label{EMKC11}
\vartheta^{11}_{X} := \theta^{100}_{X},
\end{equation}
\begin{equation}\label{EMKC12}
\vartheta^{12}_{X} := 1_{C(X)},
\end{equation}
\begin{equation}\label{EMKC21}
\vartheta^{21}_{X} := (\theta^{100}_{C(X)}\theta^{100}_{X},C(\theta^{100}_{X})\theta^{100}_{X}),
\end{equation}
\begin{equation}\label{EMKC22}
\vartheta^{22}_{X} := (\theta^{100}_{X},\theta^{100}_{X}).
\end{equation}
\begin{proof}
The proof is a straighforward consequence of the axioms for a comonad $C$.

\end{proof}
\end{theorem}

\section{Normal colax coalgebras}
Let us investigate the conditions that a strictly normal colax $\D$-coalgebra needs to satisfy. By Definition \ref{colaxDcoalgebra} the 2-cell $\bm{\zeta}$ in (\ref{zeta}) must be identity. It follows that:
\begin{itemize}
\item [1)] The identity $\de_{0}\bm{F}_{0}=I_{\A}$ holds which means that the functor $\bm{F}_{0} \maps \A \to (\X,G)$ has the form $\bm{F}_{0}=(H,\chi,I_{\A})$ where $H \maps \A \to \X$ is a functor and $\chi \maps H \Rightarrow G$ is a natural transformation

\item [2)] The identity $\de_{1}\bm{F}_{1}=I_{\X}$ holds which means that the functor $\bm{F}_{1} \maps \X \to (\X,G)$ has the form $\bm{F}_{1}=(I_{\X},\varrho,K)$ where $K \maps \X \to \A$ is a functor and $\varrho \maps I_{\X} \Rightarrow GK$ is a natural transformation

\item [3)] As a consequence of the first two items the 1-cell $\bm{F}_{G}=(\bm{F}_{1},\bm{\varphi}, \bm{F}_{0}) \maps G \to \D(G)$ has the from $((I_{\X},\varrho,K),(\omega,\varepsilon), (H,\chi,I_{\A}))$ as in the following diagram
\begin{equation}\label{strcoalg}
\begin{array}{c}
\xymatrix@!=3pc{\A \ar[r]^-{(H,\chi,I_{\A})} \ar[d]_-{G} & (\X,G) \ar@{=}[d]^{} \\
\X \ar[r]_-{(I_{\X},\varrho,K)} & (\X,G) \ultwocell<\omit>{(\omega,\varepsilon)\,\,\,\,\,\,\,\,\,} }
\end{array}
\end{equation}
where $\omega \maps G \Rightarrow H$ and $\varepsilon \maps KG \Rightarrow I_{\A}$ are natural transformations such that a diagram
\begin{equation}\label{coh}
\begin{array}{c}
\xymatrix@!=3pc{G(A) \ar[d]_{\varrho_{G(A)}} \ar[r]^-{\omega_{A}} & H(A) \ar[d]^-{\chi_{A}}\\
GKG(A) \ar[r]_-{G(\varepsilon_{A})} & G(A)}
\end{array}
\end{equation}
commutes for any object $A$ in $\A$.

\item [4)] The horizontal composition $\del \circ \bm{\varphi}$ is a natural transformation whose component indexed by an object $A$ in $\A$ is defined by a commutative diagram
\[\xymatrix@!=0.75pc{\de_{1}(I_{\X},\varrho,K)G(A) \ar[rrrr]^-{\de_{1}((\omega,\varepsilon)_{A})} \ar[dd]_{\delta_{(C,\varrho,K)(A)}} &&&& \de_{1}(H,\chi,Q)G(A) \ar[dd]^-{\delta_{(H,\chi,I_{\A})(A)}} \\
&&&& \\
G\de_{0}(I_{\X},\varrho,K)G(A) \ar[rrrr]_-{G\de_{0}((\omega,\varepsilon)_{A})} &&&& G\de_{0}(H,\chi,I_{\A})(A) }\]
which is identical to (\ref{coh}). Then the normality condition says that (\ref{coh}) not only commutes but that its diagonal is equal to the identity $1_{G(A)}$ giving us the following two identities:
\begin{equation}\label{normalcoh}
G(\varepsilon_{A})\varrho_{G(A)} = 1_{G(A)},
\end{equation}
\begin{equation}\label{normalcoh2}
\chi_{A}\omega_{A} = 1_{G(A)}.
\end{equation}
\end{itemize}
Furthermore, the two equations (\ref{thetaleq1}) and (\ref{thetaleq0}) become the following two equations
\begin{equation}\label{normalthetaleq1}
\begin{array}{c}
\de_{1}\bm{\theta}^{1} =  \iota_{(I_{\X},\varrho,K)},
\end{array}
\end{equation}
\begin{equation}\label{normalthetaleq0}
\begin{array}{c}
\de_{0}\bm{\theta}^{0} =  \iota_{(H,\chi,I_{\A})}
\end{array}
\end{equation}
saying that
\[\xymatrix@!=0.25pc{X \ar[dd]_{\varrho_{X}} \ar[rr]^-{\theta^{100}_{X}} && X \ar[dd]^{\varrho_{X}} \\
&& \\
GK(X) \ar[rr]_-{G(\theta^{110}_{X})} && GK(X)  }
\hspace{1cm}
\xymatrix@!=0.25pc{H(A) \ar[dd]_{\chi_{A}} \ar[rr]^-{\theta^{001}_{A}} && H(A) \ar[dd]^{\chi_{A}}  \\
&& \\
G(A) \ar[rr]_-{G(\theta^{011}_{A})} && G(A) }\] 
are equal to identity morphisms $(1_{X},1_{K(X)}) \maps (X,\varrho_{X},K(X)) \to (X,\varrho_{X},K(X))$ and $(1_{H(A)},1_{A}) \maps (H(A),\chi_{A},A) \to (H(A),\chi_{A},A)$ in $(\X,G) $ respectively.  This gives us the equations
\begin{equation}\label{normalleftuniteqX1}
\theta^{100}_{X} = 1_{X},
\end{equation}
\begin{equation}\label{normalleftuniteqX2}
\theta^{110}_{X} = 1_{K(X)},
\end{equation}
\begin{equation}\label{normalleftuniteqA1}
\theta^{001}_{A} = 1_{H(A)},
\end{equation}
\begin{equation}\label{normalleftuniteqA2}
\theta^{011}_{A} = 1_{A}.
\end{equation}
The identities (\ref{thetareq1}) and (\ref{thetareq0}) reduce to the following two equations
\begin{equation}\label{normalthetareq1}
\begin{array}{c}
D(\del_{G})\bm{\theta}^{1} =  \iota_{(C,\varrho,K)},
\end{array}
\end{equation}
\begin{equation}\label{normalthetareq0}
\begin{array}{c}
D(\del_{G})\bm{\theta}^{0} =  \iota_{(H,\chi,Q)}
\end{array}
\end{equation}
meaning that
\[\xymatrix@!=0.5pc{X \ar[dd]_{\varrho_{X}} \ar[rr]^-{\theta^{100}_{X}} && X \ar[dd]^{\varrho_{X}}  \\
&& \\
GK(X) \ar[rr]_-{G(\theta^{111}_{X})} && GK(X)}
\hspace{1cm}
\xymatrix@!=0.5pc{H(A) \ar[dd]_{\chi_{A}} \ar[rr]^-{\theta^{000}_{A}} && H(A) \ar[dd]^{\chi_{Q(A)}}  \\
&& \\
G(A) \ar[rr]_-{G(\theta^{011}_{A})} && G(A) }\] 
are equal to identity morphisms $(1_{X},1_{K(X)}) \maps (X,\varrho_{X},K(X)) \to (X,\varrho_{X},K(X))$ and $(1_{H(A)},1_{A}) \maps (H(A),\chi_{A},A) \to (H(A),\chi_{A},A)$ in $(\X,G) $ respectively, giving the equations
\begin{equation}\label{normalrightuniteqX1}
\theta^{100}_{X} = 1_{X},
\end{equation}
\begin{equation}\label{normalrightuniteqX2}
\theta^{111}_{X} = 1_{K(X)},
\end{equation}
\begin{equation}\label{normalrightuniteqA1}
\theta^{000}_{A} = 1_{H(A)},
\end{equation}
\begin{equation}\label{normalrightuniteqA2}
\theta^{011}_{A} = 1_{A}
\end{equation}
which together with equations (\ref{normalleftuniteqX1}) - (\ref{normalleftuniteqA2}) show that components $\bm{\theta}^{1}_{X}$ and $\bm{\theta}^{0}_{A}$ are
\[\hspace{-1cm}\xymatrix@!=0.35pc{X \ar[ddd]_{\varrho_{X}} \ar[rrr]^-{1_{X}} \ar[ddrr]^{1_{X}} &&& X \ar@{-->}[ddd]_{\varrho_{X}} \ar[dr]^{\varrho_{X}} &&\\
&&&& GK(X) \ar@{-->}[ddd]^{\varrho_{GK(X)}} \ar[dr]^{\omega_{K(X)}} &\\
&& X \ar[ddd]_{\varrho_{X}} \ar[rrr]^-{\theta^{101}_{X}} &&& HK(X) \ar[ddd]^{\chi_{K(X)}}  \\
GK(X) \ar[ddrr]_{G(1_{K(X)})} \ar@{-->}[rrr]^-{G(1_{K(X)})} &&& GK(X) \ar@{-->}[dr]_(0.4){GK(\varrho_{X})} && \\
&&&& GKGK(X) \ar@{-->}[dr]_(0.4){G(\varepsilon_{K(X)})} & \\
&& GK(X) \ar[rrr]_{G(1_{K(X)})} &&& GK(X) }
\xymatrix@!=0.35pc{H(A) \ar[ddd]_{\chi_{A}} \ar[rrr]^-{1_{H(A)}} \ar[ddrr]^{1_{H(A)}} &&& H(A) \ar@{-->}[ddd]_{\varrho_{H(A)}} \ar[dr]^{\chi_{A}} &&\\
&&&& G(A) \ar@{-->}[ddd]^{\varrho_{G(A)}} \ar[dr]^{\omega_{A}} &\\
&& H(A) \ar[ddd]_{\chi_{A}} \ar[rrr]^-{1_{H(A)}} &&& H(A) \ar[ddd]^{\chi_{A}}  \\
G(A) \ar[ddrr]_{G(1_{A})} \ar@{-->}[rrr]^-{G(\theta^{010}_{A})} &&& GKH(A) \ar@{-->}[dr]_(0.4){GK(\chi_{A})} && \\
&&&& GKG(A) \ar@{-->}[dr]_(0.4){G(\varepsilon_{A})} & \\
&& G(A)  \ar[rrr]_{G(1_{A})} &&& G(A) }\]
respectively. From the left diagram representing $\bm{\theta}^{1}_{X}$  we obtain the following identities 
\begin{equation}\label{normalthetaeq1X}
\chi_{K(X)}\theta^{101}_{X} =  \varrho_{X},
\end{equation}
\begin{equation}\label{normalthetaeq2X}
\omega_{K(X)}\varrho_{X} = \theta^{101}_{X}
\end{equation}
\begin{equation}\label{normalthetaeq3X}
\varepsilon_{K(X)}K(\varrho_{X}) = 1_{K(X)}.
\end{equation}
The last equation in conjunction with (\ref{normalcoh}) says that $K$ is a left adjoint to $G$ with unit $\varrho$ and counit $\varepsilon$.
From the right diagram for $\bm{\theta}^{0}_{A}$ we obtain the following identities
\begin{equation}\label{normalthetaeq1A}
\varrho_{H(A)} =  G(\theta^{010}_{A})\chi_{A},
\end{equation}
\begin{equation}\label{normalthetaeq2A}
\omega_{A}\chi_{A}=1_{H(A)},
\end{equation}
\begin{equation}\label{normalthetaeq4A}
\varepsilon_{A}K(\chi_{A})\theta^{010}_{A} = 1_{A}.
\end{equation}
Then (\ref{normalthetaeq2A}) in conjunction with (\ref{normalcoh2}) implies that for every object $A$ in $\A$ we have
\begin{equation}\label{chiomega}
\omega_{A}^{-1}=\chi_{A}
\end{equation}
and diagrams
\[\xymatrix@!=0.25pc{G(A) \ar[ddd]_{\varrho_{G(A)}} \ar[rrr]^-{1_{G(A)}} \ar[ddrr]^{1_{G(A)}} &&& 
G(A) \ar@{-->}[ddd]_{\varrho_{G(A)}} \ar[dr]^{\varrho_{G(A)}} \ar[rrr]^-{\omega_{A}}  &&& 
H(A) \ar@{-->}[ddd]^{\varrho_{H(A)}} \ar[dr]^{\chi_{A}} &&\\
&&&& GKG(A) \ar@{-->}[ddd]_{\varrho_{GKG(A)}}  \ar[rrr]^-{G(\varepsilon_{A})} \ar[dr]^(0.6){\omega_{KG(A)}} &&& G(A) \ar@{-->}[ddd]^{\varrho_{G(A)}} \ar[dr]^{\omega_{A}} &\\
&& G(A) \ar[ddd]_(0.6){\varrho_{G(A)}} \ar[rrr]^{\theta^{101}_{G(A)}} &&& HKG(A) \ar[ddd] \ar[rrr]^{H(\varepsilon_{A})} &&& H(A) \ar[ddd]^{\chi_{A}}  \\
GKG(A) \ar[ddrr]_{GK(1_{G(A)})} \ar@{-->}[rrr]^-{G(1_{KG(A)})} &&& GKG(A) \ar@{-->}[dr]_(0.3){GK(\varrho_{G(A)})} \ar@{--}[rrr]^-{GK(\omega_{A})} &&& GKH(X) \ar@{-->}[dr]_(0.3){GK(\chi_{A})} && \\
&&&& GKGKG(A) \ar@{-->}[dr]_(0.3){G(\varepsilon_{KG(A)})} \ar@{-->}[rrr]_{GKG(\varepsilon_{A})} &&& 
GKG(A) \ar@{-->}[dr]_(0.3){G(\varepsilon_{A})} & \\
&& GKG(A)  \ar[rrr]_{G(1_{KG(A)})} &&& GKG(A)  \ar[rrr]_{G(\varepsilon_{A})} &&& G(A) }\]
\[\xymatrix@!=0.25pc{G(A) \ar[ddd]_{\varrho_{G(A)}} \ar[rrr]^-{\omega_{A}} \ar[ddrr]^{G(1_{A})} &&& H(A) \ar@{-->}[ddd]_{\chi_{A}} \ar[ddrr]^{H(1_{A})} \ar[rrr]^-{1_{H(A)}}  &&& 
H(A) \ar@{-->}[ddd]_{\varrho_{H(A)}} \ar[dr]^{\chi_{A}} &&\\
&&&&&&& G(A) \ar@{-->}[ddd]^{\varrho_{G(A)}} \ar[dr]^{\omega_{A}} &\\
&& G(A) \ar[ddd]_{\varrho_{G(A)}} \ar[rrr]^{\omega_{A}} &&& H(A) \ar[ddd]_{\chi_{A}} \ar[rrr]^{1_{H(A)}} &&& H(A) \ar[ddd]^{\chi_{A}}  \\
GKG(A) \ar[ddrr]_{GKG(1_{A})} \ar@{-->}[rrr]^-{G(\varepsilon_{A})} &&& G(A) \ar@{-->}[drdr]_{G(1_{A})} \ar@{-->}[rrr]^-{G(\theta^{010}_{A})} &&& GKH(A) \ar@{-->}[dr]_(0.3){GK(\chi_{A})} && \\
&&&&&&& GKG(A) \ar@{-->}[dr]_(0.3){G(\varepsilon_{A})} & \\
&& GKG(A)  \ar[rrr]_{G(\varepsilon_{A})} &&& G(A)  \ar[rrr]_{G(1_{A})} &&& G(A) }\]
have the same boundary by identity (\ref{theta2}) which implies the two nontrivial equations
\begin{equation}\label{normalthetacompeq1}
K(\omega_{A})=\theta^{010}_{A}\varepsilon_{A},
\end{equation}
\begin{equation}\label{normalthetacompeq2}
H(\varepsilon_{A})\theta^{101}_{G(A)}=\omega_{A},
\end{equation}
and (\ref{thetaX111}) becomes
\begin{equation}\label{normalthetaX111}
K(\theta^{101}_{X}) =  \theta^{010}_{K(X)} 
\end{equation}
while (\ref{thetaX000})-(\ref{thetaX110}) and (\ref{thetaA000})-(\ref{thetaA111}) become trivial except for (\ref{thetaA010}) which becomes
\begin{equation}\label{normalthetaA010}
H(\theta^{010}_{A}) =  \theta^{101}_{H(A)}. 
\end{equation}

Therefore we have proved the following
\begin{theorem}\label{strictnormalcolaxcoalgebra}
A functor $G \maps \A \to \X$ admits a strictly normal colax $\D$-coalgebra structure $(G,\bm{F}_{G},\bm{\theta})$ if and only if the components $\theta^{1}_{X}$ and $\theta^{0}_{A}$ have the following form 
\[\hspace{-1cm}\xymatrix@!=0.75pc@R1pc@C1pc{X \ar[dddd]_{\varrho_{X}} \ar[rrrr]^-{1_{X}} \ar[dddrrr]^{1_{X}} &&&& X \ar@{-->}[dddd]_{\varrho_{X}} \ar[dr]^{\varrho_{X}} &&&\\
&&&&& GK(X) \ar@{-->}[dddd]^(0.6){\varrho_{GK(X)}} \ar[dr]^(0.75){\theta^{101}_{GK(X)}} &&\\
&&&&&& HKGK(X) \ar@{-->}[dddd]^{\chi_{KGK(X)}}  \ar[dr]^(0.75){H(\varepsilon_{K(X)})} &\\
&&& X \ar[dddd]_{\varrho_{X}} \ar[rrrr]^(0.4){\theta^{101}_{X}} &&&& HK(X) \ar[dddd]^{\chi_{K(X)}}  \\
GK(X) \ar[dddrrr]_{G(1_{K(X)})} \ar@{-->}[rrrr]^-{G(1_{K(X)})} &&&& GK(X) \ar@{-->}[dr]_(0.3){GK(\varrho_{X})} &&& \\
&&&&& GKGK(X) \ar@{-->}[dr]_(0.3){G(1_{KGK(X)})} && \\
&&&&&& GKGK(X) \ar@{-->}[dr]_(0.3){G(\varepsilon_{K(X)})} &\\
&&& GK(X) \ar[rrrr]_{G(1_{K(X)})} &&&& GK(X) }
\xymatrix@!=0.75pc@R1pc@C1pc{H(A) \ar[dddd]_{\chi_{A}} \ar[rrrr]^-{1_{H(A)}} \ar[dddrrr]^{1_{H(A)}} &&&&  H(A) \ar@{-->}[dddd]_{\varrho_{H(A)}} \ar[dr]^{\chi_{A}} &&&\\
&&&&& G(A) \ar@{-->}[dddd]^(0.6){\varrho_{G(A)}} \ar[dr]^(0.75){\theta^{101}_{G(A)}} &&\\
&&&&&& HKG(A) \ar@{-->}[dddd]^{\chi_{KG(A)}}  \ar[dr]^(0.75){H(\varepsilon_{A})} &\\
&&& H(A) \ar[dddd]_{\chi_{A}} \ar[rrrr]^-{1_{H(A)}} &&&& H(A) \ar[dddd]^{\chi_{A}}   \\
G(A) \ar[dddrrr]_{G(1_{A})} \ar@{-->}[rrrr]^-{G(\theta^{010}_{A})} &&&& GKH(A) \ar@{-->}[dr]_(0.3){GK(\chi_{A})} &&& \\
&&&&& GKG(A) \ar@{-->}[dr]_(0.3){G(1_{KG(A)})} && \\
&&&&&& GKG(A) \ar@{-->}[dr]_(0.3){G(\varepsilon_{A})} &\\
&&& G(A) \ar[rrrr]_{G(1_{A})} &&&& G(A) }\]
respectively.  
\end{theorem}
Now let us show some notable examples of nortmal $\D$-coalgebras.
\begin{example}\label{Morita}
By taking $G=H$, and $\chi=\iota_{G}= \omega$ we obtain what Morita called in \cite{Morita} a \emph{strongly adjoint pair} consisting of an adjoint triple
\[G \dashv K \dashv G\]
where $G$ is simultaneously left and right adjoint of $K$. In the sequel to this paper \cite{Baka} we will follow a convention by Caenepeel, Militaru and Zhu \cite{CMZ} who suggested to refer to such $G$ as a \emph{Frobenius functor} followed by \cite{CIGTN}.
\end{example}
\begin{example}\label{Morita2}
By keeping $\omega$ and $\chi$ as mutually invertible natural transformations we end up with an ambidextrous adjunction
\[H \dashv K \dashv G\]
(or sometimes \emph{ambiadjunction} for short) which were pivotal in the work of Lauda \cite{Lauda} who showed that every Frobenius object $M$ in a monoidal category $\M$ arises from an ambijunction in some 2-category $\D$ into which $M$ fully and faithfully embeds. This result also shows that every 2D TQFT is obtained from an ambijunction in some 2-category since it is well known that a 2D topological quantum field theory is equivalent to a commutative Frobenius algebra.
\end{example}

The 2-category $\Cat$ of small categories, functors and natural transformations is a symmetric monoidal 2-category object in the cartesian closed 2-category $2\Cat$ of 2-categories, 2-functors and natural 2-transformations with respect to a coproduct in $\Cat$, i.e. a disjoint union of categories. Since any object $2\Cat$ is a comonoid with compultiplication given by a diagonal there exists a fibered 2-monad $\A r$ on $\Cat \times \Cat$ whose strict algebras are functors.  More precisely, there exists a canonical 2-bifibration $Pr_{1} \maps \Cat \times \Cat \to \Cat$ (from the existence of finite products in $2\Cat$) and the fibered 2-functor
\begin{equation}
\begin{array}{c}\label{Arr2fibered}
\xymatrix@!=1pc{\Cat \times \Cat \ar[rr]^-{\A r} \ar[dr]_-{Pr_{1}} && \Cat \times \Cat \ar[dl]^{Pr_{1}} \\
& \Cat &  }
\end{array}
\end{equation}
which is defined for any pair $(\X,\Y)$ of categories by $\A r(\X,\Y):= (\X, \X + \Y)$ where "+" denotes a coproduct in $\Cat$. For any pair of functors $F \maps \X \to \X'$ and $G \maps \Y \to \Y'$ we have $\A r(F,G):= (F, F + G)$ with obvious action on pairs of natural transformations. Then $\A r^{2}(\X,\Y):= \A r(\X, \X + \Y):= (\X, \X + \X + \Y)$ and there exist natural 2-transformations 
\begin{equation}
\begin{array}{c}\label{MN}
\M \maps \A r^{2} \Rightarrow \A r \\
\N \maps \I_{\Cat \times \Cat} \Rightarrow \A r
\end{array}
\end{equation}
with one component $(I_{\X}, \bigtriangledown_{\X} + I_{\Y}) \maps (\X, \X + \X + \Y) \to (\X, \X + \Y)$ and the other one $(I_{\X}, \iota_{\Y}) \maps (\X, \Y) \to (\X, \X + \Y)$ indexed by $(\X,\Y)$ respectively, as in the following diagrams
\[\hspace{-.5cm}
\xymatrix@!=4pc{(\X, \Y) \ar[r]^-{(I_{\X}, \iota_{\Y})} \ar[d]_-{(F, G)} & (\X, \X + \Y) \ar[d]^{(F,F + G)} \\
(\X', \Y') \ar[r]_-{(I_{\X'}, \iota_{\Y'})} & (\X', \X' + \Y')}
\hspace{.5cm}
\xymatrix@!=4pc{(\X, \X + \X +\Y) \ar[rr]^-{(I_{\X}, I_{\X} + I_{\X} + I_{\Y)}} \ar[d]_-{(F, F + F + G)} && (\X, \X + \Y) \ar[d]^{(F,F + G)} \\
(\X', \X' + \X' + \Y') \ar[rr]_-{(I_{\X'}, I_{\X'} + I_{\X'} + I_{\Y'})} && (\X', \X' + \Y')}\]
where $\iota_{\Y}$ is a canonical inclusion. It follows that we have components of natural 2-transformations
\begin{equation}
\begin{array}{c}\label{Str2alg}
\A r(\N_{(\X,\Y)})=\A r(I_{\X}, \iota_{\Y})=(I_{\X}, I_{\X} + \iota_{\Y})\\
\N_{\A r(\X,\Y)}= \N_{(\X, \X + \Y)} = (I_{\X}, I_{\X} +\iota_{\X + \Y})
\end{array}
\end{equation}
which send any object $(X,X')$ of ${(\X,\X + \Y)}$ where both $X$ and $X'$ are objects of $\X$ to an object of ${(\X,\X + \X + \Y)}$ which we also denote by $(X,X')$ where $X'$ now lies in the first or the second copy of the category $\X$ in the coproduct $\X + \X + \Y$ respectively. It follows that we have a fully faithful adjoint string
\begin{equation}
\begin{array}{c}\label{adjointstring}
\A r(\N_{(\X,\Y)}) \dashv \M_{(\X,\Y)} \dashv \N_{\A r(\X,\Y)}
\end{array}
\end{equation} 
together with a canonical iso-modification $\Phi$ whose component indexed by $(\X,\Y)$
\begin{equation}
\begin{array}{c}\label{Phi}
\Phi_{(\X,\Y)} \maps \A r(\N_{(\X,\Y)}) \Rightarrow \N_{\A r(\X,\Y)}
\end{array}
\end{equation} 
flips the two copes of the category $\X$ in the coproduct $\X + \X + \Y$.  Obviously we have $\Phi^{2}=Id$.
\begin{remark}\label{2monad1cell}
If one takes seriously a viewpoint which identifies $Mnd(\C)$ with the 2-category of lax algebras for the identity 2-monad on the 2-category then one could consider a diagram (\ref{Arr2fibered}) as a 1-cell
\begin{equation}
\begin{array}{c}\label{Arr2fibered2}
\xymatrix@!=1pc{\Cat \times \Cat \ar[rr]^-{\A r} \ar[d]_-{Pr_{1}} && \Cat \times \Cat \ar[d]^{Pr_{1}} \\
\Cat \ar[rr]^-{\I_{\Cat}} &&  \Cat}
\end{array}
\end{equation}
in the 2-category $2bi\Fib$ of 2-bifibrations, bicartesian 2-functors (which preserve both cartesian and cocartesian cells and cartesian transformations.  Then (\ref{Arr2fibered2}) is the underlying 1-cell of a monad in the 2-category $2bi\Fib$.  One can further show that the fibered 2-monad $\A r$ on $\Cat \times \Cat$ is monoidal monad.
\end{remark}

\begin{theorem}\label{Frobenius2monad}
The fibered 2-functor (\ref{Arr2fibered}) together with 2-natural transformations (\ref{MN}) and the iso-modification (\ref{Phi}) is a Frobenius fibred 2-monad.
\end{theorem}

\section{Conclusion and future directions}
In the sequels to this paper we will give a complete description of the Eilenberg-Moore 2-category of colax $\D$-coalgebras, colax morphisms between them and their transformations and we will show how many fundamental constructions in formal category theory like adjoint triples, distributive laws, comprehension structures, Frobenius functors etc. naturally fit in this context. Then we will proceed to describe various pseudo distributive laws between a comma 2-comonad and its cousins - the associated split fibration 2-monad and the associated split cofibration 2-monad. The former is an instance of a pseudodistributive law which Garner used \cite{Gar} in his description of Szabo’s polycategories, and the pseudoalgebras for the latter are Beck-Checalley fibrations. We will also show how this contexts is related to Bunge and Funk admissible 2-monads whose Eilenberg-Moore 2-category of algebras are characterised in terms of (co)completeness. Finally we will describe Kleisi 2-category of the associated split fibration 2-monad by means of its bifibrations which are defined by a Bunge and Funk in \cite{BF} by means of a certain bicomma object condition and the corresponding comprehensive factorization for those 1-cells which have an admissible domain.

One of the main motivations for the author to investigate the structures in this paper was his following observation:In the context of enhanced category theory which was initiated by Lack and Shulman in \cite{LSh} who developed a theory of what we could call 1-enhanced 2-categories, the authors considered the following situation:
\begin{itemize}
\item a {\it property} of 1-cells in a 2-category
\item a process of {\it assigning} properties to objects of the 2-category $\Cat$
\end{itemize}
I will show that besides its cartesian monoidal structure, the category $\F$ whose objects are injective on objects and fully faithful functors which Lack and Shulman called full embeddings has a much richer structure given by other closed (but not monoidal) structures.  One of the consequences of this fact is that the following notions are essentially equivalent:
\begin{itemize}
\item [(i)] A 2-category with a right ideal of 1-cells
\item [(ii)] A category enriched over the closed category $\F_{c}$ whose objects are functors that are fully faithful and injective on objects which I christened - enhanced categories
\end{itemize}
Any 2-category with Yoneda structure can be naturally seen as an enriched category in the sense of Eilenberg and Kelly \cite{EK} and the theory developed by Street and Walters in \cite{SW} has a natural interpretation in this context. 

\refs

\bibitem[Arkor,  Di Liberti, Loregian, 2024]{ADiLL}
N. Arkor,  I. Di Liberti, F. Loregian, Adjoint functor theorems for lax-idempotent pseudomonads, Theory and Applications of Categories, Vol. 41, No. 20 (2024), 667-685.

\bibitem[Arkor, Bourke, Ko, 2025]{ABK}
N. Arkor, J. Bourke, J. Ko, Enhanced 2-categorical structures, two-dimensional limit sketches and the symmetry of internalisation, preprint (2025), arXiv:2501.12510.

\bibitem[Bakovi\' c, 2025]{Baka}
I. Bakovi\' c, Comma 2-comonad II: Monads in Eilenberg-Moore 2-category, in preparation

\bibitem[Baez, Dolan, 1998]{BD}
J.C.  Baez, J.  Dolan,  Higher-Dimensional Algebra III.n-Categories and the Algebra of Opetopes, Advances in Mathematics Volume 135, Issue 2, (1998),  145-206.

\bibitem[Beck, 1969]{Beck}
J. Beck,  Distributive laws, In: Eckmann, B. (eds) Seminar on Triples and Categorical Homology Theory, Lecture Notes in Mathematics, vol 80 (1969), 119-140.

\bibitem[B\' enabou, 1967]{Be1}
J. B\' enabou, Introduction to bicategories, Reports of the Midwest Category Seminar, Lecture Notes in Math. 47 (1967), 1-77.

\bibitem[B\' enabou, 1973]{Be2}
J. B\'enabou, Les distributeurs, Institut de mathematiques pure and appliqu\' ee, Rapport No. 3 (1973),Universit\' e catholique de Louvain.

\bibitem[Blackwell,Kelly, Power,1989]{BKP}
R. Blackwell, G.M. Kelly, A.J. Power, Two-dimensional monad theory, Journal of Pure and Applied Algebra Volume 59, Issue 1 (1989), 1-41.

\bibitem[Bourke,2014]{Bo}
J. Bourke, Two-dimensional monadicity, Advances in Mathematics 252 (2014), 708-747.

\bibitem[Bourke, Lobbia,2023]{BL}
J. Bourke, G. Lobbia, A skew approach to enrichment for Gray-categories, Advances in Mathematics Volume 434 (2023), 109327.

\bibitem[Borceaux, Bourn,2004]{BB}
F. Borceaux, D. Bourn, Mal'cev, Protomodular, Homological and Semi-abelian categories, Mathematics and its Applications vol. ~566, Kluwer Academic Publishers (2004).

\bibitem[Bunge,1974]{Bu}
M. Bunge, Coherent Extensions and Relational Algebras, Transactions of the American Mathematical Society, Vol. 197 (1974), 355-390.

\bibitem[Bunge,  Carboni,, 1995]{BC}
M. Bunge,  A. Carboni, The symmetric topos, Journal of Pure and Applied Algebra 105 (1995), 233-249.

\bibitem[Bunge, Funk, 1999]{BF}
M. Bunge, J. Funk, On a bicomma object condition for KZ-doctrines, Journal of Pure and Applied Algebra 143 (1999), 69-105.

\bibitem[Burroni, 1971]{Burr}
A. Burroni, $T$-catégories (cat\' egories dans un triple), Cahiers de topologie et g\' eom\' etrie diff\' erentielle cat\' egoriques Volume 12 (1971) no. 3, pp. 215-321.

\bibitem[Casta\~{n}o Iglesias, G\' omez-Torrecillas, Nastasescu, 1998]{CIGTN}
F. Casta\~{n}o Iglesias, J. G\' omez-Torrecillas, C. Nastasescu, Frobenius functors: applications, Communications in Algebra vol.~27, Issue 10 (1998), 4879-4900.

\bibitem[Caenepeel, Militaru, Zhu,1997]{CMZ}
S. Caenepeel, G. Militaru, S. Zhu, Doi-Hopf modules, Yetter-Drinfel\' d modules and Frobenius type properties, Trans. Amer. Math. Soc. 349 (1997), 4311-4342.

\bibitem[Cheng, Gurski, Riehl, 2014]{CGR}
E. Cheng, N. Gurski, E. Riehl, Cyclic multicategories, multivariable adjunctions and mates,Journal of K-Theory 13 (2014), 337-396.

\bibitem[Chikladze, 2015]{Chik}
D. Chikladze, Lax formal theory of monads, monoidal approach to bicategorical structures and generalized operads,Theory and Applications of Categories, Vol. 30, No. 10, (2015), 332-386.

\bibitem[Cockett, Koslowski, Seely, Wood, 2003]{CKSW}
J.R.B. Cockett, J. Koslowski, R.A.G. Seely, R.J. Wood,  Modules, Theory and Applications of Categories, Vol. 11, No. 17 (2003), 375-396.

\bibitem[Crans, 1999]{Crans}
S. Crans, A tensor product for Gray-categories, Theory and Applications of Categories 5, No. 2 (1999),12-69.

\bibitem[Cruttwell, Shulman, 2010]{CS}
G.S.H. Cruttwell, M. Shulman, A unified framework for generalized multicategories, Theory and Applications of Categories 24, No. 21 (2010), 580-655.

\bibitem[Dawson, Par\' e,Pronk, 2003]{DPP}
R. Dawson, R. Par\' e, D.Pronk, Adjoining adjoints, Advances in Mathematics 178 (2003),  99-140.

\bibitem[Di Liberti, FLoregian, 2019]{DiLL}
I. Di Liberti, F. Loregian, On the unicity of formal category theories, preprint (2019), arXiv:1901.01594.

\bibitem[Duskin,2000]{Dus1}
J. Duskin, Simplicial matrices and the nerves of weak n-categories I: nerves of bicategories, Theory and Applications of Categories 9 (2000), 198-308.

\bibitem[Duskin, 2000]{Dus2}
J. Duskin, Simplicial Matrices and the Nerves of Weak n-Categories II : Nerves of Morphisms of Bicategories, incomplete draft.

\bibitem[Dyckhoff, Tholen, 1987]{DT2}
R. Dyckhoff, W. Tholen, Exponentiable morphisms, partial products and pullback complements, Journal of Pure Applied Algebra 49 (1987), 103-116.

\bibitem[Eilenberg, Kelly,  1966]{EK}
S.Eilenberg, G.M. Kelly,  Closed Categories, in: Eilenberg, Harrison,, MacLane, R\" ohrl, (eds.), Proceedings of the Conference on Categorical Algebra, (1966), 421-562.

\bibitem[Emmenegger, Mesiti, Rosolini, Streicher,2024]{EGMS}
J. Emmenegger, L. Mesiti, G. Rosolini, T.  Streicher, A comonad for Grothendieck fibrations, Theory and Applications of Categories Vol. 40, No. 13 (2024), 371-389.

\bibitem[Fiore, 1995]{Fi}
M. Fiore, Lifting as a KZ-doctrine, in: D. Pitt et al. (Eds.), Category Theory and Computer Science, Lecture Notes in Computer Science, vol. 953, Springer, Berlin (1995), 146-158.

\bibitem[Fiore, Gambino, Hyland, Wynskel,  2018]{FGHW}
M. Fiore, N. Gambino, M. Hyland, G. Wynskel, Relative pseudomonads, Kleisli bicategories, and substitution monoidal structures, Selecta Mathematica Volume 24 (2018), 2791-2830.

\bibitem[Fumex, Ghani, Johann, 2011]{FGJ}
C. Fumex, N. Ghani, P. Johann,  Indexed induction and coinduction, fibrationally, in Algebra and Coalgebra in Computer Science: 4th International Conference, CALCO 2011, Winchester, UK, August 30-September 2, 2011, Proceedings, volume 6859.

\bibitem[Garner, 2008]{Gar}
R. Garner, Polycategories via pseudo-distributive laws, Advances in Mathematics Volume 218, Issue 3 (2008), 781-827.

\bibitem[Garner, Shulman, 2016]{GS}
R. Garner, M. Shulman, Enriched categories as a free cocompletion, Advances in Mathematics 289 (2016), 1-94.

\bibitem[Gray,  1966]{Gray1}
J.W. Gray, Fibred and cofibred categories, in: Proceedings of the Conference on Categorical Algebra (La Jolla, California 1965), Springer,(1966), 21-83.

\bibitem[Gray,  1966]{Gray2}
J.W. Gray, The categorical comprehension scheme, In: P.J. Hilton (ed.) Category Theory, Homology Theory and Their Applications III, Lecture Notes in Mathematics, vol 99 (1966), 242-312.

\bibitem[Gray,  1974]{Gray}
J.W. Gray, Formal Category Theory: Adjointness for 2-Categories, Lecture Notes in Mathematics vol.~391, 1974.

\bibitem[Grothendieck et al.,1971]{Gr2}
A. Grothendieck, Catégories fibr\'ees et d\'escente, Expos\'e VI in Rev\^etements Etale et Groupe Fondamental, SGA 1, Lecture Notes in Math. 224, Springer-Verlag, Berlin, 1971.

\bibitem[Grothendieck, 2022]{Gr}
A. Grothendieck, \' A la poursuite des champs Volume I, \' edit\' e par G. Maltsiniotis, Documents Math\' ematiques Volume 20, Soci\' et\' e Math\' ematique de France, 2022.

\bibitem[Hermida, Jacobs, 1998]{HJ}
C. Hermida, B. Jacobs, Structural Induction and Coinduction in a Fibrational Setting, Information and Computation 145 (1998), 107-152.

\bibitem[Hermida, 2000]{Her1}
C. Hermida,  Representable multicategorires, Advances in Mathematics 151 (2000), 164-225.

\bibitem[Hermida, 2001]{Her2}
C. Hermida,  From coherent structures to universal properties,  Journal of Pure and Applied Algebra 165 (2001), 7-61.

\bibitem[Hyland, Power, 2002]{HP}
M.  Hyland, J.  Power,  Pseudo-commutative monads and pseudo-closed 2-categories,  Journal of Pure and Applied Algebra 175 (2002), 141-185.

\bibitem[Jacobs,1993]{Jac}
B. Jacobs, Comprehension categories and the semantics of type dependency, Theoretical Computer Science Volume 107, Issue 2 (1993), 169-207.

\bibitem[Johnstone,1993]{J1}
P.T. Johnstone, Fibrations and Partial Products in a 2-Category,  Applied Categorical Structures 1 (1993), 141-179.

\bibitem[Johnstone,1995]{J2}
P.T. Johnstone, Connected limits, familial representability and Artin glueing, Mathematical Structures in Computer Science vol. 5 (1995),  441-459.

\bibitem[Johnstone,2011]{J3}
P.T. Johnstone, Remarks on punctual local connectedness, Theory and Applications of Categories, Vol. ~25, No. 3 (2011), 51-63.

\bibitem[Kelly,1974]{Kelly4}
G.M. Kelly, Coherence theorems for lax algebras and for distributive laws, Category Seminar Sydney 1972/73, Lecture Notes in Mathematics, 420,(1974), 281-375.

\bibitem[Kelly, Lawvere,1989]{KL}
G.M. Kelly, F.W. Lawvere, On the complete lattice of essential localizations, Bulletin de la Soci\' et\' e Math\' ematique de Belgique, S\' erie A, v.ol.~ 41, no 2 (1989), 289-319.

\bibitem[Kock,1995]{Kock}
A. Kock, Monads for which structures are adjoint to units, Journal of Pure and Applied Algebra 104, Issue 1 (1995), 41-59.

\bibitem[Koudenburg, 2020]{Kou3}
S. R. Koudenburg, Augmented virtual double categories, Theory and Applications of Categories 35, No. 10 (2020), 261-325.

\bibitem[Korostenski,Tholen,]{KT}
M. Korostenski, W. Tholen, Factorization systems as Eilenberg-Moore algebras, Journal of Pure and Applied Algebra Volume 85, Issue 1 (1993), 57-72.

\bibitem[Lack, Street, 2002]{LSt}
S. Lack, R. Street, The formal theory of monads II, Journal of Pure and Applied Algebra vol.~175 (2002) 243--265.

\bibitem[Lack,Shulman,2012]{LSh}
S. Lack, M. Shulman, Enhanced 2-categories and limits for lax morphisms, Advances in Mathematics vol.~229 (2012), 294-356.

\bibitem[lack, Miranda,2024]{LM}
S. Lack, A. Miranda, What is the universal property of the 2-category of monads?, Theory and Applications of Categories, Vol. 42, No. 1 (2024), 2-18.

\bibitem[Lauda,2006]{Lauda}
A. D. Lauda, Frobenius algebras and ambidextrous adjunctions, Theory and Applications of Categories, Vol. 16,  No. 4 (2006),84-122.

\bibitem[Lawvere,  1963]{Law0}
F. W. Lawvere, Functorial Semantics of Algebraic Theories, Reprints in Theory and Applications of Categories, No. 5 (2004) 1-121.

\bibitem[Lawvere,  1965]{Law1}
F. W. Lawvere, The category of categories as a foundation for mathematics., Proceedings of the Conference on Categorical Algebra, La Jolla 1965,  S. Eilenberg, D. K. Harrison, S. MacLane, and H. R\" ohrl (eds.), Springer-Verlag New York Inc., New York (1966), 1-20.

\bibitem[Lawvere,  1969]{Law2}
F. W. Lawvere, Ordinal Sums and Equational Doctrines, Springer Lecture Notes in Mathematics No. 80, Springer-Verlag (1969), 141-155.

\bibitem[Lawvere,  1969]{Law3}
F.W. Lawvere, Adjointness in Foundations, Dialectica 23 (1969), 281-296.

\bibitem[Lawvere,  1970]{Law4}
F.W. Lawvere, Equality in hyperdoctrines and comprehension schema as an adjoint functor, Proceedings of the AMS Symposium on Pure Mathematics XVII (1970), 1-14.

\bibitem[Lawvere,  2007]{Law5}
F.W. Lawvere, Axiomatic cohesion, Theory and Applications of Categories Vol. 19, No. 3 (2007), 41-49.

\bibitem[Leinster,2004]{Le}
T. Leinster, Higher Operads, Higher Categories, Cambridge University Press, 2004.

\bibitem[Linton,1965]{Lin}
F.E.J. Linton,  Autonomous Categories and Duality of Functors, J. Alg. 2 (1965), p.315-349.

\bibitem[Maltsiniotis, 2005]{Maltsiniotis}
G. Maltsiniotis, Structures d'asph\' ericit\' e, foncteurs lisses, et fibrations, Annales math\' ematiques Blaise Pascal, Volume 12 (2005) no. 1,1-39.

\bibitem[Marmolejo,1997]{Ma}
F. Marmolejo,  Doctrines whose structure forms a fully faithful adjoint string, Theory and Applications of Categories, Vol. 3, No. 2 (1997), 23-44.

\bibitem[Marmolejo,2012]{Ma1}
F. Marmolejo,  Kan extensions and lax idempotent pseudomonads, Theory and Applications of Categories, Vol. 26, No. 1, 2012, pp. 1-29.

\bibitem[Melli\` es, 2012]{Mell1}
P.A. Melli\` es, Dialogue categories up to deformation, preprint (2012),359-412.

\bibitem[Melli\` es, 2016]{Mell2}
P.A. Melli\` es, Dialogue Categories and Chiralities, Publications of the Research Institute for Mathematical Sciences 52 (2016),359-412.

\bibitem[Melli\` es,Zeilberger, 2015]{MZ1}
P.A. Melli\` es, N. Zeilberger, Functors are Type Refinement Systems.,POPL '15: Proceedings of the 42nd Annual ACM SIGPLAN-SIGACT Symposium on Principles of Programming Languages (2015) 3-16.

\bibitem[Melli\` es,Zeilberger, 2016]{MZ}
P.A. Melli\` es, N. Zeilberger, A bifibrational reconstruction of Lawvere's presheaf hyperdoctrine, Proceedings of the 31st Annual ACM/IEEE Symposium on Logic in Computer Science,  LICS (2016), 555-564.

\bibitem[Melli\` es, Rolland, 2020]{MR}
P.A. Melli\` es, N. Rolland, Comprehension and quotient structures in the language of 2-categories, Proceedings of the 5th International Conference on Formal Structures for Computation and Deduction, Z. M. Ariola (ed.),  (2020), 6:1-6:18.

\bibitem[Mac Lane,1971]{McL}
S.  Mac Lane, Categories for the Working Mathematician,  Second Edition, Originally published by Springer-Verlag New York (1971).

\bibitem[Maranda, 1965]{Mara}
J.M. Maranda, Formal Categories, Canadian Journal of Mathematics vol.~ 17 (1965), 758-801.

\bibitem[Morita, 1965]{Morita}
K. Morita, Adjoint pairs of functors and Frobenius extensions, Science Reports of the Tokyo Kyoiku Daigaku, Section A 9 202/208 (1965) 40-71.

\bibitem[Osmond, 2012]{O1}
A. Osmond, On Diers theory of Spectrum I : Stable functors and right multi-adjoints, arXiv:2012.00853.

\bibitem[Palmquist, 1971]{Palm}
P.H.  Palmquist,  The double category of adjoint squares. In: Reports of the Midwest Category Seminar V, Lecture Notes in Mathematics, vol 195, (1971).

\bibitem[Pavlovi\' c,  1997]{Pav}
D. Pavlovi\' c, Chu I: cofree equivalences, dualities and *-autonomous categories, Mathematical Structures in Computer Science Volume 7,  Issue 1 (1997), 49-73.

\bibitem[Par\' e, 1971]{Pare1}
R. Par\' e, On absolute colimits, Journal  of Algebra 19 (1971), 80-95.

\bibitem[Par\' e,Rosebrugh, Wood, 1989]{PRW}
R. Par\' e, R. Rosebrugh, R.J. Wood, Idempotents in bicategories, Bulletin of the Australian Mathematical Society, Volume 39 , Issue 3 (1989), 421-434.

\bibitem[Power,  Watanabe, 2002]{PW}
A.J. Power,  H. Watanabe, Combining a monad and a comonad, Theoretical Computer Science Volume 280, Issues 1-2 (2002), 137-162.

\bibitem[Power, 2007]{Power3}
A.J. Power, Three dimensional monad theory, In: Categories in algebra, geometry and mathematical physics, Contemporary Mathematics Vol. 431 (2007), 405-426.

\bibitem[Rosicky, Tholen, 2007]{RT}
J. Rosicky,  W.Tholen, Factorization, fibration and torsion,  Journal of Homotopy and Related Structures Volume 2, Issue 2 (2007), 295-314.

\bibitem[Schanuel, Street,1986]{SS}
S. Schanuel, R. Street, The free adjunction, Cahiers de topologie et g\' eom\' etrie 7 (1986), 81-83.

\bibitem[Shulman, 2018]{Shu}
M. Shulman, Contravariance through enrichment, Theory and Applications of Categories, Vol. 33, No. 5 (2018), 95-130.

\bibitem[Street, 1972]{Street1}
R. Street, The formal theory of monads, Journal of Pure and Applied Algebra Volume 2, Issue 2 (1972), 149-168. 

\bibitem[Street, 1974]{Street2}
R. Street, Fibrations and Yoneda's lemma in a $2$-category, Category Seminar (Proc. Sem., Sydney, 1972/1973),  pp. 104-133. Lecture Notes in Mathematics, Vol. 420, Springer, Berlin, 1974.

\bibitem[Street,Walters,  1978]{SW}
R. Street, R. Walters, Yoneda structures on 2-categories, Journal of Algebra vol.~50, Issue 2 (1978), 350-379.

\bibitem[Street,Tholen, Wischnewsky, Wolff, 1980]{STWW}
R. Street, W. Tholen, M.B. Wischnewsky, H. Wolff, Semi-topological functors III: Lifting of monads and adjoint functors, Journal of Pure and Applied Algebra, vol.~16, Issue 3, (1980), 299-314.

\bibitem[Streicher, Weinberger, 2021]{StW}
T. Streicher, J. Weinberger, Simplicial sets inside cubical sets, Theory and Applications of Categories, Vol. 37 No. 10, (2021),  276-286.

\bibitem[Walker, 2018]{Wal0}
C. Walker, Yoneda structures and KZ doctrines, Journal of Pure and Applied Algebra 222, Issue 6 (2018), 1375-1387.

\bibitem[Walker, 2019]{Wal1}
C. Walker, Distributive Laws via Admissibility,  Applied Categorical Structures 27 (2019),567-617.

\bibitem[Wood, 1982]{Wo1}
R. J. Wood,  Abstract pro arrows I, Cahiers de Topologie et G\' eometrie Diff\' erentielle 23 (1982), 279-290.

\bibitem[Wood, 1985]{Wo2}
R. J. Wood,  Proarrows II, Cahiers de topologie et g\' eom\' etrie diff\' erentielle cat\' egoriques 26, no. 2 (1985), 135-168. 

\bibitem[Yanofsky, 2000]{Yan}
N. Yanofsky, The syntax of coherence, Cahiers de topologie et g\' eom\' etrie diff\' erentielle cat\' egoriques 41, no 4 (2000),  255-304.

\endrefs

\end{document}